\documentclass[a4paper,12pt]{article}
\usepackage{amsmath}
\usepackage{amscd}
\usepackage[xdvi]{graphicx}
\usepackage[dvips]{color}
\usepackage{times}
\usepackage{txfonts}

\newcommand{\bm}[1]{\mbox{\boldmath $#1$}}

\makeatletter

\@addtoreset{equation}{section}

\begin{document}

\begin{center}
\textbf{\LARGE{Extension and Unification of Singular Perturbation Methods for ODEs
Based on the Renormalization Group Method}}
\end{center}

\begin{center}
Department of Applied Mathematics and Physics

Kyoto University, Kyoto, 606-8501, Japan

Hayato CHIBA \footnote{E mail address : chiba@amp.i.kyoto-u.ac.jp}
\end{center}
\begin{center}
September 29, 2008; Revised May 1 2009
\end{center}

\setlength{\baselineskip}{14pt}

\textbf{Abstract}

The renormalization group (RG) method is one of the singular perturbation methods which
is used in search for asymptotic behavior of solutions of differential equations.
In this article, time-independent vector fields and time (almost) periodic vector fields are considered.
Theorems on
error estimates for approximate solutions, existence of approximate invariant manifolds and their stability,
inheritance of symmetries from those for the original equation to those for the RG equation, are proved.
Further it is proved that the RG method unifies traditional singular perturbation methods,
such as the averaging method, the multiple time scale method, the (hyper-) normal forms theory,
the center manifold reduction, the geometric singular perturbation method and the phase reduction.
A necessary and sufficient condition for the convergence of the infinite order RG equation is also investigated.

\section{Introduction}

Differential equations form a fundamental topic in mathematics and its application  to natural sciences.
In particular, perturbation methods occupy an important place in the theory of differential equations.
Although most of the differential equations can not be solved exactly, some of them are close
to solvable problems in some sense, so that 
perturbation methods, which provide techniques to handle such class of problems, have been long studied.

This article deals with a system of ordinary differential equations (ODEs) on a manifold $M$ of the form
\begin{equation}
\frac{dx}{dt} = \varepsilon g(t, x, \varepsilon ), \quad x\in M,
\end{equation}
which is almost periodic in $t$ with appropriate assumptions (see the assumption (A) in Sec.2.1),
where $\varepsilon \in \mathbf{R}$ or $\mathbf{C}$ is a small parameter.

Since $\varepsilon $ is small, it is natural to try to construct a solution of this system 
as a power series in $\varepsilon $
of the form
\begin{equation}
x = \hat{x}(t) = x_0(t) + \varepsilon x_1(t) + \varepsilon ^2 x_2(t) + \cdots .
\end{equation}
Substituting Eq.(1.2) into Eq.(1.1) yields a system of ODEs on $x_0, x_1, x_2, \cdots $.
The method to construct $\hat{x}(t)$ in this manner is called the \textit{regular perturbation method}.

It is known that if the function $g(t, x, \varepsilon )$ is analytic in $\varepsilon $,
the series (1.2) converges to an exact solution of (1.1), while if it is not analytic, (1.2) diverges
and no longer provides an exact solution.
However, the problem arising immediately is that one can not calculate infinite series like (1.2) in general whether it converges
or not, because it involves infinitely many ODEs on $x_0, x_1, x_2, \cdots $.
If the series is truncated at a finite-order term in $\varepsilon $, another problem arises.
For example, suppose that Eq.(1.1) admits an exact solution $x(t) = \sin (\varepsilon t)$,
and that we do not know the exact solution.
In this case, the regular perturbation method provides a series of the form
\begin{equation}
\hat{x}(t) = \varepsilon t - \frac{1}{3!}(\varepsilon t)^3 + \frac{1}{5!}(\varepsilon t)^5 + \cdots .
\end{equation}
If truncated, the series becomes a polynomial in $t$, which 
diverges as $t \to \infty$ although the exact solution is periodic in $t$.
Thus, the perturbation method fails to predict qualitative properties of the exact solution.
Methods which handle such a difficulty and provide acceptable approximate solutions are called
\textit{singular perturbation methods}.
Many singular perturbation methods have been proposed so far [1,3,6,38,39,43,44,47,50] and
many authors reported that some of them produced the same results though procedures
to construct approximate solutions were different from one another [9,38,43,44,47].

The renormalization group (RG) method is the relatively new method proposed by Chen, Goldenfeld and Oono [8,9],
which reduces a problem to a more simple equation called the \textit{RG equation}, based on an idea
of the renormalization group in quantum field theory.
In their papers, it is shown (without a proof) that the RG method unifies conventional singular perturbation
methods such as the multiple time scale method, the boundary layer technique, the WKB analysis and so on.

After their works, many studies of the RG method have been done [10-13,16,21,22,26-32,35,36,42,45,46,48,49,53,55-59].
Kunihiro [28,29] interpreted an approximate solution obtained by the RG method as the envelope of a family of 
regular perturbation solutions.
Nozaki \textit{et al.} [45,46] proposed the proto-RG method to derive the RG equation
effectively. Ziane [53] , DeVille \textit{et al.} [32] and Chiba [10] gave error estimates for approximate solutions.
Chiba [10] defined the higher order RG equation and the RG transformation to improve error estimates.
He also proved that the RG method could provide approximate vector fields and approximate invariant manifolds
as well as approximate solutions.
Ei, Fujii and Kunihiro [16] applied the RG method to obtain approximate center manifolds and their method was
rigorously formulated by Chiba [12].
DeVille \textit{et al.} [32] showed that lower order RG equations are equivalent to normal forms of vector
fields, and this fact was extended to higher order RG equations by Chiba [11].
Applications to partial differential equations are appeared in [16,31,42,45,46,55-59].

One of the purposes of this paper is to give basic theorems on the RG method extending
author's previous works [10-12], in which the RG method is discussed for more restricted problems than Eq.(1.1).
At first, definitions of the higher order RG equations for Eq.(1.1) are given
and properties of them are investigated.
It is proved that the RG method provides approximate vector fields (Thm.2.5)
and approximate solutions (Thm.2.7) along with error estimates.
Further, it is shown that if the RG equation has a normally hyperbolic invariant manifold $N$,
the original equation (1.1) also has an invariant manifold $N_\varepsilon $ which is diffeomorphic
to $N$ (Thms.2.9, 2.14).
The RG equation proves to have the same symmetries (action of Lie groups)
as those for the original equation (Thm.2.12).
In addition, if the original equation is an autonomous system, the RG equation is shown to have
an additional symmetry (Thm.2.15).
These facts imply that the RG equation is easier to analyze than the original equation.
An illustrative example to verify these theorems is also given (Sec.2.5).

The other purpose of this paper is to show that the RG method extends and unifies
other traditional singular perturbation methods, such as 
the averaging method (Sec.4.1), the multiple time scale method (Sec.4.2), the (hyper-) normal forms theory (Sec.4.3),
the center manifold reduction (Sec.3.2), the geometric singular perturbation method (Sec.3.3),
the phase reduction (Sec.3.4), and Kunihiro's method [28,29] based on envelopes (Sec.4.4).
A few of these results were partially obtained by many authors [9,38,43,44,47].
The present arguments will greatly reveal the relations among these methods.

Some properties of the infinite order RG equation are also investigated.
It is proved that the infinite order RG equation converges if and only if the original equation is invariant under 
an appropriate torus action (Thm.5.1).
This result extends Zung's theorem [54] which gives a necessary and sufficient condition for the convergence
of normal forms of infinite order.
The infinite RG equation for a time-dependent linear system proves to be convergent (Thm.5.6)
and be related to monodromy matrices in Floquet theory.

Throughout this paper, solutions of differential equations are supposed to be defined for all $t\in \mathbf{R}$.

%%%%%%%%%%%%%%%%%%%%%%%%%%%%%%%%%%%%%%%%%%%%%%%%%%%%%%%%%%%%%%%%%%%%%%%%%%%%%%%%%%%%%%%%%%%%%%%%%%%%%%
%%%%%%%%%%%%%%%%%%%%%%%%%%%%%%%%%%%%%%%%%%%%%%%%%%%%%%%%%%%%%%%%%%%%%%%%%%%%%%%%%%%%%%%%%%%%%%%%%%%%%%

\section{Renormalization group method}

In this section, we give the definition of the RG (renormalization group) equation and 
main theorems on the RG method, such as the existence of invariant manifolds and inheritance of 
symmetries.
An illustrative example and comments on symbolic computation of the RG equation are also provided.

\subsection{Setting, definitions and basic lemmas}

Let $M$ be an $n$ dimensional manifold and $U$ an open set in $M$ whose closure $\bar{U}$ is compact.
Let $g(t, \bm{\cdot}, \varepsilon )$ be a vector field on $U$ parameterized by 
$t\in \mathbf{R}$ and $\varepsilon \in \mathbf{C}$.
We consider the system of differential equations 
\begin{equation}
\frac{dx}{dt} = \dot{x} = \varepsilon g(t, x, \varepsilon ).
\end{equation}
For this system, we make the following assumption.
\\[0.2cm]
\textbf{(A)} \, The vector field $g(t,x,\varepsilon )$ is $C^1$ with respect to time $t \in \mathbf{R}$ and 
$C^\infty$ with respect to $x\in U$ and $\varepsilon \in \mathcal{E}$, 
where $\mathcal{E} \subset \mathbf{C}$ is a small neighborhood
of the origin. Further, $g$ is an almost periodic function with respect to $t$ uniformly in $x \in \bar{U}$ and 
$\varepsilon \in \bar{\mathcal{E}}$, the set of whose Fourier exponents has no accumulation points on $\mathbf{R}$.
\\[0.2cm]
\indent In general, 
a function $h(t,x)$ is called \textit{almost periodic with respect to $t$ uniformly in $x \in \bar{U}$} if the set
\begin{eqnarray*}
T(h, \delta ) := \{ \tau \,\, | \,\, || h(t+ \tau, x) - h(t,x) || < \delta, \quad
 \forall t\in \mathbf{R} ,\, \forall x\in \bar{U}\} \subset \mathbf{R}
\end{eqnarray*}
is relatively dense for any $\delta >0$; that is, 
there exists a positive number $L$ such that $[a, a+L] \cap T(h,\delta ) \neq \emptyset$
for all $a\in \mathbf{R}$. It is known that an almost periodic function is expanded in a Fourier series as 
$h(t,x) \sim \sum a_n(x) e^{i \lambda _n t},\,\, (i=\sqrt{-1})$, where $\lambda _n \in \mathbf{R}$ is called a \textit{Fourier exponent}.
See Fink [20] for basic facts on almost periodic functions.
The condition for Fourier exponents in the above assumption (A) is essentially used to prove Lemma 2.1 below.
We denote $\mathrm{Mod}(h)$ the smallest additive group of real numbers that contains the Fourier exponents $\lambda _n$ of 
an almost periodic function $h(t)$ and call it the \textit{module} of $h$.

Let $\sum^{{\infty}}_{k=1} \varepsilon ^k g_k(t,x)$ be the formal Taylor expansion of $\varepsilon g(t,x,\varepsilon )$
in $\varepsilon $ : 
\begin{equation}
\dot{x} = \varepsilon g_1(t,x) + \varepsilon ^2 g_2(t,x) + \cdots .
\end{equation}
By the assumption (A), we can show that $g_i(t,x) \,\, (i = 1,2, \cdots )$ are almost periodic functions
with respect to $t\in \mathbf{R}$ uniformly in $x\in \bar{U}$ such that $\mathrm{Mod}(g_i) \subset \mathrm{Mod}(g)$.

Though Eq.(2.1) is mainly considered in this paper,
we note here that Eqs.(2.3) and (2.5) below are reduced to Eq.(2.1):
Consider the system of the form
\begin{equation}
\dot{x} = f(t,x) + \varepsilon g(t,x, \varepsilon ),
\end{equation}
where $f(t, \bm{\cdot})$ is a $C^\infty$ vector field on $U$ and $g$ satisfies the assumption (A).
Let $\varphi _t$ be the flow of $f$ ; that is, $\varphi _t(x_0)$ is a solution of the equation $\dot{x} = f(t,x)$
whose initial value is $x_0$ at the initial time $t=0$.
For this system, changing the coordinates by $x = \varphi _t(X)$ provides
\begin{equation}
\dot{X} = \varepsilon \left( \frac{\partial \varphi _t}{\partial X}(X)\right) ^{-1}g(t, \varphi _t(X),\varepsilon )
:= \varepsilon \tilde{g}(t, X, \varepsilon ).
\end{equation}
We suppose that
\\[0.2cm]
\textbf{(B)} \, the vector field $g$ satisfies the assumption (A) and there exists an open set $W \subset U$
such that $\varphi _t(W) \subset U$ and $\varphi _t(x)$ is almost periodic with respect to $t$ uniformly in $x\in \bar{W}$,
the set of whose Fourier exponents has no accumulation points.
\\[0.2cm]
Under the assumption (B), we can show that the vector field $\tilde{g}(t, X, \varepsilon )$ in the right hand side of Eq.(2.4)
satisfies the assumption (A), in which $g$ is replaced by $\tilde{g}$.
Thus Eq.(2.3) is reduced to Eq.(2.1) by the transformation $ x \mapsto X$.

In many applications, Eq.(2.3) is of the form
\begin{eqnarray}
\dot{x} &=& Fx + \varepsilon g(x, \varepsilon ) \nonumber \\
&=& Fx + \varepsilon g_1(x) + \varepsilon ^2 g_2(x) + \cdots ,\,\, x\in \mathbf{C}^n,
\end{eqnarray}
where 
\\[0.2cm]
\textbf{(C1)}\, the matrix $F$ is a diagonalizable $n \times n$ constant matrix all of whose eigenvalues are on the 
imaginary axis,
\\
\textbf{(C2)}\, each $g_i (x)$ is a polynomial vector field on $\mathbf{C}^n$.
\\[0.2cm]
Then, the assumptions (C1) and (C2) imply the assumption (B) because $\varphi _t(x) = e^{Ft}x$ is almost periodic.
Therefore the coordinate transformation $x = e^{Ft}X$ brings Eq.(2.5) into the form of Eq.(2.1) :
$\dot{X} = \varepsilon e^{-Ft}g(e^{Ft}X, \varepsilon ) := \varepsilon \tilde{g}(t, X, \varepsilon )$.
In this case, $\mathrm{Mod}(\tilde{g})$ is generated by the absolute values of the eigenvalues of $F$.
Note that any equations $\dot{x} = f(x)$ with $C^\infty$ vector fields $f$ such that $f(0) = 0$
take the form (2.5) if we put $x \mapsto \varepsilon x$ and expand the equations in $\varepsilon $.

In what follows, we consider Eq.(2.1) with the assumption (A).
We suppose that the system (2.1) is defined on an open set $U$ on Euclidean space $M=\mathbf{C}^n$.
However, all results to be obtained below can be easily extended to those for a system on an arbitrary manifold
by taking local coordinates.
Let us substitute $x = x_0 + \varepsilon x_1 + \varepsilon ^2 x_2 + \cdots $ into the right hand side of Eq.(2.2) and expand
it with respect to $\varepsilon $. We write the resultant as
\begin{equation}
\sum^{\infty}_{k=1}\varepsilon ^k g_k(t, x_0 + \varepsilon x_1 + \varepsilon ^2 x_2 + \cdots)
 = \sum^{\infty}_{k=1}\varepsilon ^k G_k(t, x_0, x_1, \cdots , x_{k-1}).
\end{equation}
For instance, $G_1, G_2, G_3$ and $G_4$ are given by
\begin{eqnarray}
 G_1(t,x_0) & =&  g_1(t, x_0), \\
 G_2(t,x_0,x_1)& =&  \frac{\partial g_1}{\partial x}(t,x_0)x_1 + g_2(t, x_0), \\
 G_3(t,x_0,x_1,x_2)& =&  \frac{1}{2}\frac{\partial ^2g_1}{\partial x^2}(t,x_0)x_1^2
     + \frac{\partial g_1}{\partial x}(t, x_0)x_2 + \frac{\partial g_2}{\partial x}(t,x_0)x_1 + g_3(t, x_0), \\
 G_4(t,x_0,x_1,x_2,x_3) 
 &=&  \frac{1}{6}\frac{\partial ^3g_1}{\partial x^3}(t,x_0)x_1^3 + \frac{\partial ^2g_1}{\partial x^2}(t,x_0)x_1x_2
 + \frac{\partial g_1}{\partial x}(t,x_0)x_3 \nonumber \\
& & + \frac{1}{2} \frac{\partial ^2g_2}{\partial x^2}(t,x_0)x_1^2 + \frac{\partial g_2}{\partial x}(t,x_0)x_2
 + \frac{\partial g_3}{\partial x}(t,x_0)x_1 + g_4(t,x_0),
\end{eqnarray}
respectively.
Note that $G_i\,\, (i = 1,2, \cdots )$ are almost periodic functions with respect to $t$ uniformly in $x\in \bar{U}$
such that $\mathrm{Mod}(G_i) \subset \mathrm{Mod}(g)$.
With these $G_i$'s, we define the $C^\infty$ maps $R_i, u^{(i)}_t : U \to M$ to be
\begin{eqnarray}
& & R_1(y) = \lim_{t\to \infty}\frac{1}{t} \int ^t\! G_1(s,y)ds, \\
& & u^{(1)}_t(y) = \int^t \! \left( G_1(s,y) - R_1(y) \right) ds,
\end{eqnarray}
and 
\begin{eqnarray}
& & R_i(y)  
= \lim_{t\to \infty}\frac{1}{t} \int ^t\! \Bigl(
G_i(s, y, u^{(1)}_s(y),
 \cdots, u^{(i-1)}_s(y))- \sum^{i-1}_{k=1}\frac{\partial u_s^{(k)}}{\partial y}(y)R_{i-k}(y) \Bigr) ds , \\
& & u^{(i)}_t(y)= \int^t \! \Bigl( G_i(s, y,u^{(1)}_s(y)
, \cdots , u^{(i-1)}_s(y)) - \sum^{i-1}_{k=1}\frac{\partial u_s^{(k)}}{\partial y}(y) R_{i-k}(y) - R_i(y) \Bigr) ds,
\end{eqnarray}
for $i=2,3,\cdots$, respectively, where $\int^t$ denotes the indefinite integral, whose integral constants are 
fixed arbitrarily (see also Remark 2.4 and Section 2.4).
\\[0.2cm]
\textbf{Lemma 2.1.} \, (i)\, The maps $R_i \,\, (i = 1,2,\cdots )$ are well-defined (i.e. the limits exist).
\\
(ii) \, The maps $u^{(i)}_t(y)\,\, (i = 1,2,\cdots )$ are almost periodic functions with respect to $t$ 
uniformly in $y\in \bar{U}$ such that $\mathrm{Mod}(u^{(i)}) \subset \mathrm{Mod}(g)$.
In particular, $u_t^{(i)}$ are bounded in $t\in \mathbf{R}$.
\\[0.2cm]
\textbf{Proof.}\, We prove the lemma by induction.
Since $G_1(t,y) = g_1(t,y)$ is almost periodic, it is expanded in a Fourier series of the form
\begin{equation}
g_1(t,y) = \sum_{\lambda _n\in \mathrm{Mod} (g_1)} a_n(y) e^{i \lambda _n t}, \,\, \lambda _n\in \mathbf{R},
\end{equation} 
where $\lambda _0 = 0$. Clearly $R_1(y)$ coincides with $a_0(y)$. Thus $u_t^{(1)}(y)$ is written as
\begin{equation}
u_t^{(1)}(y) = \int^t \! \sum_{\lambda _n \neq 0} a_n(y) e^{i\lambda _n s} ds. 
\end{equation}
In general, it is known that the primitive function $\int \! h(t,y) dt$ 
of an uniformly almost periodic function $h(t,y)$ is also uniformly 
almost periodic if the set of Fourier exponents of $h(t,y)$ is bounded away from
zero (see Fink [20]). Since the set of Fourier exponents of $g_1(t,y) - R_1(y)$ is bounded away from
zero by the assumption (A), $u_t^{(1)}(y)$ is almost periodic and calculated as
\begin{equation}
u_t^{(1)}(y) 
  = \sum_{\lambda _n \neq 0} \frac{1}{i \lambda _n} a_n(y) e^{i\lambda _n t} + (\mathrm{integral\,\, constant}).
\end{equation}
This proves Lemma 2.1 for $i=1$.

Suppose that Lemma 2.1 holds for $i=1,2, \cdots ,k-1$.
Since $G_k(t, x_0, \cdots , x_{k-1})$ and $u_t^{(1)}(y) , \cdots , u_t^{(k-1)}(y)$ are uniformly almost periodic functions,
the composition $G_k(t,y, u_t^{(1)}(y), \cdots , u_t^{(k-1)}(y))$ is also an uniformly almost periodic function
whose module is included in $\mathrm{Mod}(g)$ (see Fink [20]).
Since the sum, the product and the derivative with respect to a parameter $y$ of uniformly almost periodic functions are 
also uniformly almost periodic (see Fink [20]), the integrand in Eq.(2.13) is an uniformly almost periodic function,
whose module is included in $\mathrm{Mod}(g)$.
The $R_k(y)$ coincides with its Fourier coefficient associated with the zero Fourier exponent.
By the assumption (A), the set of Fourier exponents of the integrand in Eq.(2.13) has no accumulation points.
Thus it turns out that the set of Fourier exponents of the integrand in Eq.(2.14) is bounded away from zero.
This proves that $u_t^{(k)}(y)$ is uniformly almost periodic and the proof of Lemma 2.1 is completed. \hfill $\blacksquare$ 
\\[-0.2cm]

Before introducing the RG equation, we want to explain how it is derived according to Chen, Goldenfeld and Oono [8,9].
The reader who is not interested in formal arguments can skip the next paragraph and go to Definition 2.2.

At first, let us try to construct a formal solution of Eq.(2.1) by the regular perturbation method;
that is, substitute Eq.(1.2) into Eq.(2.1).
Then we obtain a system of ODEs
\begin{eqnarray*}
\left\{ \begin{array}{l}
\dot{x}_0 = 0,  \\
\dot{x}_1 = G_1(t, x_0), \\
\quad \vdots  \\
\dot{x}_n = G_n(t, x_0,\cdots , x_{n-1}), \\
\quad \vdots 
\end{array} \right.
\end{eqnarray*}
Let $x_0(t) = y \in \mathbf{C}^n$ be a solution of the zero-th order equation.
Then, the first order equation is solved as
\begin{eqnarray*}
x_1(t) = \int^t \! G_1(s,y) ds
 = R_1(y)t + \int^t \! \left( G_1(s,y) - R_1(y) \right) ds = R_1(y)t + u^{(1)}_t(y),
\end{eqnarray*}
where we decompose $x_1(t)$ into the bounded term $u^{(1)}_t(y)$ and the divergence term $R_1(y)t$ 
called the \textit{secular term}.
In a similar manner, we solve the equations on $x_2, x_3, \cdots $ step by step.
We can show that solutions are expressed as
\begin{eqnarray*}
x_n(t) = u^{(n)}_t(y) + \left( R_n(y) + \sum^{n-1}_{k=1} \frac{\partial u^{(k)}}{\partial y}(y)R_{n-k}(y)\right) t + O(t^2),
\end{eqnarray*}
(see Chiba [10] for the proof).
In this way, we obtain a formal solution of the form
\begin{eqnarray*}
\hat{x}(t):=\hat{x}(t,y) = y + \sum^\infty_{n=1}\varepsilon ^n u^{(n)}_t(y)
 + \sum^\infty_{n=1} \varepsilon ^n \left( R_n(y) + \sum^{n-1}_{k=1} \frac{\partial u^{(k)}}{\partial y}(y)R_{n-k}(y)\right) t + O(t^2).
\end{eqnarray*}
Now we introduce a dummy parameter $\tau \in \mathbf{R}$ and replace polynomials $t^j$ in the above by $(t-\tau)^j$.
Next, we regard $y = y(\tau)$ as a function of $\tau$ to be determined so that we recover the formal solution $\hat{x}(t,y)$:
\begin{eqnarray*}
\hat{x}(t,y) = y(\tau) + \sum^\infty_{n=1}\varepsilon ^n u^{(n)}_t(y(\tau))
 + \sum^\infty_{n=1} \varepsilon ^n \left( R_n(y(\tau)) 
+ \sum^{n-1}_{k=1} \frac{\partial u^{(k)}}{\partial y}(y(\tau))R_{n-k}(y(\tau))\right) (t-\tau) + O((t-\tau)^2).
\end{eqnarray*}
Since $\hat{x}(t,y)$ has to be independent of the dummy parameter $\tau$, we impose the condition
\begin{eqnarray*}
\frac{d}{d\tau}\Bigl|_{\tau = t} \hat{x}(t,y) = 0,
\end{eqnarray*}
which is called the \textit{RG condition}.
This condition provides
\begin{eqnarray*}
0 &=& \frac{dy}{dt} + \sum^\infty_{n=1} \varepsilon ^n \frac{\partial u^{(n)}_t}{\partial y}(y)\frac{dy}{dt}
      - \sum^\infty_{n=1} \varepsilon ^n \left( R_n(y) + \sum^{n-1}_{k=1} \frac{\partial u^{(k)}_t}{\partial y}(y) R_{n-k}(y)\right) \\
&=& \left( id + \sum^\infty_{n=1} \varepsilon ^n  \frac{\partial u^{(n)}_t}{\partial y}(y) \right) \frac{dy}{dt}
   - \left( id + \sum^\infty_{n=1} \varepsilon ^n  \frac{\partial u^{(n)}_t}{\partial y}(y) \right) \sum^\infty_{k=1} \varepsilon ^k R_k(y).
\end{eqnarray*}
Thus we see that $y(t)$ has to satisfy the equation $dy/dt = \sum^\infty_{k=1} \varepsilon ^k R_k(y)$,
which gives the RG equation.
Motivated this formal argument, we define the RG equation as follows:
\\[0.2cm]
\textbf{Definition 2.2.}\, Along with $R_i$ and $u^{(i)}_t$, we define the 
\textit{$m$-th order RG equation} for Eq.(2.1) to be
\begin{equation}
\dot{y} = \varepsilon R_1(y) + \varepsilon ^2R_2(y) + \cdots  + \varepsilon ^mR_m(y),
\end{equation}
and the \textit{$m$-th order RG transformation} to be
\begin{equation}
\alpha ^{(m)}_t (y) = y + \varepsilon u^{(1)}_t(y) + \cdots + \varepsilon ^m u^{(m)}_t(y).
\end{equation}

Domains of Eq.(2.18) and the map $\alpha ^{(m)}_t$ are shown in the next lemma.
\\[0.2cm]
\textbf{Lemma 2.3.} \, If $|\varepsilon |$ is sufficiently small,
there exists an open set $V = V(\varepsilon ) \subset U$ such that $\alpha ^{(m)}_t(y)$ is a
diffeomorphism from $V$ into $U$, and the inverse $(\alpha^{(m)}_t)^{-1}(x) $ is also an almost periodic
function with respect to $t$ uniformly in $x$.
\\[0.2cm]
\textbf{Proof.} \, 
Since the vector field $g(t,x,\varepsilon )$ is $C^\infty$ with respect to $x$ and $\varepsilon $, 
so is the map $\alpha ^{(m)}_t$.
Since $\alpha ^{(m)}_t$ is close to the identity map if $|\varepsilon |$ is small,
there is an open set $V_t \subset U$ such that $\alpha ^{(m)}_t$ is a diffeomorphism on $V_t$.
Since $V_t$'s are $\varepsilon $-close to each other and since $\alpha ^{(m)}_t$ is almost periodic,
the set $\tilde{V} := \bigcap_{t\in \mathbf{R}} V_t$ is not empty.
We can take the subset $V \subset \tilde{V}$ if necessary so that $\alpha ^{(m)}_t(V) \subset U$.

Next thing to do is to prove that $(\alpha ^{(m)}_t)^{-1}$ is an uniformly almost periodic function.
Since $\alpha ^{(m)}_t$ is uniformly almost periodic, the set
\begin{equation}
T(\alpha ^{(m)}_t, \delta ) = 
\{ \tau \, | \, || \alpha ^{(m)}_{t+\tau}(y) - \alpha ^{(m)}_t(y) || < \delta, \quad
\forall t\in \mathbf{R},\, \forall y\in V \}
\end{equation}
is relatively dense for any small $\delta >0$.
For $y\in V$, put $x = \alpha ^{(m)}_t (y)$. Then
\begin{eqnarray}
|| (\alpha ^{(m)}_{t+\tau })^{-1}(x) - (\alpha ^{(m)}_t)^{-1}(x) ||
&=& || (\alpha ^{(m)}_{t+\tau })^{-1} (\alpha ^{(m)}_t (y)) 
    - (\alpha ^{(m)}_{t+\tau })^{-1}(\alpha ^{(m)}_{t+\tau } (y)) || \nonumber \\
&\leq & L_{t+ \tau} || \alpha ^{(m)}_t(y) - \alpha ^{(m)}_{t + \tau}(y) || < L_{t + \tau}\delta ,
\end{eqnarray}
if $\tau \in T(\alpha ^{(m)}_t, \delta ) $, where $L_t$ is the Lipschitz constant of the map $(\alpha ^{(m)}_t)^{-1}|_U$.
Since $\alpha ^{(m)}_t$ is almost periodic, we can prove that there exists the number $L:= \max_{t\in \mathbf{R}} L_t$.
Now the inequality
\begin{equation}
|| (\alpha ^{(m)}_{t+\tau })^{-1}(x) - (\alpha ^{(m)}_t)^{-1}(x) || < L \delta 
\end{equation}
holds for any small $\delta >0$, $\tau \in T(\alpha ^{(m)}_t, \delta ) $ and $x\in \alpha ^{(m)}_t(V)$.
This proves that $(\alpha ^{(m)}_t)^{-1}$ is an almost periodic function with respect to $t$ 
uniformly in $x\in \alpha ^{(m)}_t(V)$. \hfill $\blacksquare$
\\[-0.2cm]

In what follows, we suppose that the $m$-th order RG equation and the $m$-th order RG transformation
are defined on the set $V$ above.
Note that the smaller $|\varepsilon |$ is, the larger set $V$ we may take.
\\[0.2cm]
\textbf{Remark 2.4.}\, Since the integral constants in Eqs.(2.11) to (2.14) are left undetermined,
the $m$-th order RG equations and the $m$-th order RG transformations are not unique although $R_1(y)$
is uniquely determined.
However, the theorems described below hold for any choice of integral constants unless otherwise noted.
Good choices of integral constants simplify the RG equations and it will be studied in Section 2.4.

%%%%%%%%%%%%%%%%%%%%%%%%%%%%%%%%%%%%%%%%%%%%%%%%%%%%%%%%%%%%%%%%%%%%%%%%%%%%%%%%%%%%%%%%%%%%%%%%%%%%

\subsection{Main theorems}

Now we are in a position to state our main theorems.
\\[0.2cm]
\textbf{Theorem 2.5.}\, Let $\alpha ^{(m)}_t$ be the $m$-th order RG transformation for Eq.(2.1) defined on $V$ as Lemma 2.3.
If $|\varepsilon |$ is sufficiently small, there exists a vector field $S(t,y,\varepsilon )$ on $V$ parameterized by
$t$ and $\varepsilon $ such that
\\[0.2cm]
(i)\, by changing the coordinates as $x = \alpha ^{(m)}_t(y)$, Eq.(2.1) is transformed into the system
\begin{equation}
\dot{y} = \varepsilon R_1(y) + \varepsilon ^2 R_2(y) + \cdots  + \varepsilon ^mR_m(y) + \varepsilon ^{m+1}S(t, y, \varepsilon ),
\end{equation}
(ii) \, $S$ is an almost periodic function with respect to $t$ uniformly in $y\in V$ with $\mathrm{Mod} (S) \subset \mathrm{Mod} (g)$,
\\
(iii) \, $S(t,y,\varepsilon )$ is $C^1$ with respect to $t$ and $C^\infty$ with respect to $y$ and $\varepsilon $.
In particular, $S$ and its derivatives are bounded as $\varepsilon \to 0$ and $t\to \infty$.
\\[0.2cm]
\textbf{Proof.}\, The proof is done by simple calculation.
By putting $x = \alpha ^{(m)}_t(y)$, the left hand side of Eq.(2.1) is calculated as
\begin{eqnarray}
\frac{dx}{dt} &=& \frac{d}{dt}\alpha ^{(m)}_t(y) \nonumber \\
&=& \dot{y} + \sum^{m}_{k=1} \varepsilon ^k \frac{\partial u^{(k)}_t}{\partial y}(y) \dot{y}
 + \sum^{m}_{k=1}\varepsilon ^k \frac{\partial u^{(k)}_t}{\partial t}(y) \nonumber \\
&=& \left( id + \sum^{m}_{k=1} \varepsilon ^k \frac{\partial u^{(k)}_t}{\partial y}(y) \right) \dot{y}
 + \sum^{m}_{k=1}\varepsilon ^k \left( G_k(t,y, u^{(1)}_t , \cdots ,u^{(k-1)}_t) 
 - \sum^{k-1}_{j=1} \frac{\partial u^{(j)}_t}{\partial y} (y) R_{k-j}(y) - R_k(y) \right). \nonumber \\
\end{eqnarray}
On the other hand, the right hand side is calculated as
\begin{eqnarray}
\varepsilon g(t, \alpha ^{(m)}_t(y), \varepsilon )
 &=& \sum^{\infty}_{k=1} \varepsilon ^k g_k (t, y + \varepsilon  u^{(1)}_t(y) + \varepsilon ^2  u^{(2)}_t(y)+ \cdots  ) \nonumber \\
&=& \sum^{\infty}_{k=1} \varepsilon ^k G_k(t, y,  u^{(1)}_t(y) , \cdots ,  u^{(k-1)}_t(y)).
\end{eqnarray}
Thus Eq.(2.1) is transformed into
\begin{eqnarray}
\dot{y} &=& \left( id + \sum^{m}_{k=1} \varepsilon ^k \frac{\partial u^{(k)}_t}{\partial y}(y) \right)^{-1}
\sum^{m}_{k=1} \varepsilon ^k \left( R_k(y) + \sum^{k-1}_{j=1} \frac{\partial u^{(j)}_t}{\partial y} (y) R_{k-j}(y)\right) \nonumber \\
& & + \left( id + \sum^{m}_{k=1} \varepsilon ^k \frac{\partial u^{(k)}_t}{\partial y}(y) \right)^{-1}
\sum^{\infty}_{k=m+1} \varepsilon ^k G_k(t, y,  u^{(1)}_t(y) , \cdots ,  u^{(k-1)}_t(y)) \nonumber \\
&=& \left( id + \sum^{\infty}_{j=1}(-1)^j \left( \sum^{m}_{k=1} \varepsilon ^k \frac{\partial u^{(k)}_t}{\partial y}(y)\right)^j  \right)
\left( \sum^{m}_{k=1}\varepsilon ^k R_k(y) + \sum^{m}_{k=1} \varepsilon ^k \frac{\partial u^{(k)}_t}{\partial y} (y)
\sum^{m-k}_{j=1}\varepsilon ^j R_{j}(y)\right) \nonumber \\
& & + \left( id + \sum^{m}_{k=1} \varepsilon ^k \frac{\partial u^{(k)}_t}{\partial y}(y) \right)^{-1}
\sum^{\infty}_{k=m+1} \varepsilon ^k G_k(t, y,  u^{(1)}_t(y) , \cdots ,  u^{(k-1)}_t(y)) \nonumber \\
&=& \sum^{m}_{k=1}\varepsilon ^k R_k(y)
    + \sum^{\infty}_{j=1} (-1)^j \left( \sum^{m}_{k=1} \varepsilon ^k \frac{\partial u^{(k)}_t}{\partial y}(y)\right)^j
\sum^{m}_{i=m-k+1} \varepsilon ^i R_i(y) \nonumber \\
& & + \left( id + \sum^{m}_{k=1} \varepsilon ^k \frac{\partial u^{(k)}_t}{\partial y}(y) \right)^{-1}
 \sum^{\infty}_{k=m+1} \varepsilon ^k G_k(t, y,  u^{(1)}_t(y) , \cdots ,  u^{(k-1)}_t(y)).
\end{eqnarray}
The last two terms above are of order $O(\varepsilon ^{m+1})$ and almost periodic functions
because they consist of almost periodic functions $u^{(i)}_t$ and $G_i$. 
This proves Theorem 2.5. \hfill $\blacksquare$
\\[0.2cm]
\textbf{Remark 2.6.} \, To prove Theorem 2.5 (i),(iii), we do not need the assumption of almost periodicity
for $g(t,x,\varepsilon )$ as long as $R_i(y)$ are well-defined and $g, u_t^{(i)}$ and their derivatives are bounded in $t$
so that the last two terms in Eq.(2.26) are bounded.
In Chiba [10], Theorem 2.5 (i) and (iii) for $m=1$ are proved without the assumption (A) but assumptions
on boundedness of $g, u^{(i)}_t$ and their derivatives.

Thm.2.5 (iii) implies that we can use the $m$-th order RG equation to construct approximate solutions of Eq.(2.1).
Indeed, a curve $\alpha ^{(m)}_t(y(t))$, a solution of the RG equation transformed by the RG transformation,
gives an approximate solution of Eq.(2.1).
\\[0.2cm]
\textbf{Theorem 2.7 (Error estimate).} \, Let $y(t)$ be a solution of the $m$-th order RG equation and $\alpha ^{(m)}_t$
the $m$-th order RG transformation. There exist positive constants $\varepsilon _0, C$ and $T$ such that a solution $x(t)$ of
Eq.(2.1) with $x(0) = \alpha ^{(m)}_0(y(0))$ satisfies the inequality
\begin{equation}
|| x(t) - \alpha ^{(m)}_t (y(t)) || < C |\varepsilon |^m,
\end{equation}
as long as $|\varepsilon | < \varepsilon _0,\, y(t) \in V$ and $0 \leq t \leq T/|\varepsilon |$.
\\[0.2cm]
\textbf{Remark 2.8.}\, Since the velocity of $y(t)$ is of order $O(\varepsilon )$, $y(0) \in V$ implies
$y(t) \in V$ for $0 \leq t \leq T/|\varepsilon |$ unless $y(0)$ is $\varepsilon $-close to the boundary of $V$.
If we define $u^{(i)}_t$ so that the indefinite integrals in Eqs.(2.12, 14) are replaced by the definite integrals $\int^t_{0}$,
$\alpha ^{(m)}_0$ is the identity and $\alpha ^{(m)}_0(y(0)) = y(0)$.
\\[0.2cm]
\textbf{Proof of Thm.2.7.} Since $\alpha ^{(m)}_t$ is a diffeomorphism on $V$ and bounded in $t \in \mathbf{R}$,
it is sufficient to prove that a solution $y(t)$ of Eq.(2.18) and a solution $\tilde{y}(t)$ of Eq.(2.23) with
$y(0) = \tilde{y}(0)$ satisfy the inequality
\begin{equation}
|| \tilde{y}(t) - y(t) || < \tilde{C}|\varepsilon |^m,\,\, 0 \leq t \leq T/|\varepsilon |,
\end{equation}
for some positive constant $\tilde{C}$.

Let $L_1 > 0$ be the Lipschitz constant of the function $R_1(y) + \varepsilon R_2(y) + \cdots + \varepsilon ^{m-1}R_m(y)$
on $\bar{V}$ and $L_2 > 0$ a constant such that $\sup_{t\in \mathbf{R}, y\in \bar{V}} || S(t,y,\varepsilon ) || \leq L_2$.
Then, by Eq.(2.18) and Eq.(2.23), $y(t)$ and $\tilde{y}(t)$ prove to satisfy
\begin{equation}
|| \tilde{y}(t) - y(t) || \leq \varepsilon L_1 \int^t_{0} \! || \tilde{y}(s) - y(s) ||ds + L_2 \varepsilon ^{m+1} t. 
\end{equation}
Now the Gronwall inequality proves that
\begin{equation}
|| \tilde{y}(t) - y(t) || \leq \frac{L_2}{L_1} \varepsilon ^m (e^{\varepsilon L_1 t} - 1).
\end{equation}
The right hand side is of order $O(\varepsilon ^m)$ if $0 \leq t \leq T/\varepsilon $. \hfill $\blacksquare$
\\[-0.2cm]

In the same way as this proof, we can show that if $R_1(y) =\cdots  = R_{k}(y) = 0$ holds with $k \leq m$,
the inequality (2.27) holds for the longer time interval $0 \leq t \leq T/|\varepsilon |^{k+1}$.
This fact is proved by Murdock and Wang [41] for the case $k=1$ in terms of the multiple time scale method.

We can also detect existence of invariant manifolds.
Note that introducing the new variable $s$, we can rewrite Eq.(2.1) as the autonomous system
\begin{equation}
\left\{ \begin{array}{l}
\displaystyle \frac{dx}{dt} = \varepsilon g(s, x, \varepsilon ),  \\[0.2cm]
\displaystyle \frac{ds}{dt} = 1.  \\
\end{array} \right.
\end{equation}
Then we say that Eq.(2.31) is defined on the $(s,x)$ space.
\\[0.2cm]
\textbf{Theorem 2.9 (Existence of invariant manifolds).} \, Suppose that $R_1(y) = \cdots = R_{k-1}(y) = 0$ and
$\varepsilon ^k R_k(y)$ is the first non-zero term in the RG equation for Eq.(2.1).
If the vector field $R_k(y)$ has a boundaryless compact normally hyperbolic invariant manifold $N$,
then for sufficiently small $\varepsilon >0$, Eq.(2.31) has an invariant manifold $N_\varepsilon $ 
on the $(s,x)$ space which is diffeomorphic to $\mathbf{R} \times N$.
In particular, the stability of $N_\varepsilon $ coincides with that of $N$.
\\[-0.2cm]

To prove this theorem, we need Fenichel's theorem :
\\
\textbf{Theorem (Fenichel [18]).}\, 
Let $M$ be a $C^1$ manifold and $\mathcal{X}(M)$ the set of $C^1$ vector fields
on $M$ with the $C^1$ topology. 
Suppose that $f\in \mathcal{X}(M)$ has a boundaryless compact normally hyperbolic 
$f$-invariant manifold $N \subset M$. Then, the following holds:
\\[0.2cm]
(i) There is a neighborhood $\mathcal{U} \subset \mathcal{X}(M)$ of $f$ such that 
there exists a normally hyperbolic $g$-invariant manifold $N_g \subset M$ for any $g \in \mathcal{U}$.
The $N_g$ is diffeomorphic to $N$.
\\
(ii) If $|| f-g || \sim O(\varepsilon )$,
$N_g$ lies within an $O(\varepsilon )$ neighborhood of $N$ uniquely.
\\
(iii) The stability of $N_\varepsilon $ coincides with that of $N$.

Note that for the case of a compact normally hyperbolic invariant manifold with boundary,
Fenichel's theorem is modified as follows : If a vector field $f$ has a compact normally hyperbolic invariant manifold $N$
with boundary, then a vector field $g$, which is $C^1$ close to $f$, has a \textit{locally} invariant manifold $N_g$
which is diffeomorphic to $N$. In this case, an orbit of the flow of $g$ on $N_g$ may go out from $N_g$
through its boundary.
According to this theorem, Thm.2.9 has to be modified so that 
$N_\varepsilon $ is \textit{locally} invariant if $N$ has boundary.

See [18,24,51] for the proof of Fenichel's theorem and the definition of normal hyperbolicity.
\\[0.2cm]
\textbf{Proof of Thm.2.9.} \, Changing the time scale as $t\mapsto t/\varepsilon ^k$ and introducing the new variable $s$,
we rewrite the $k$-th order RG equation as
\begin{equation}
\left\{ \begin{array}{l}
\displaystyle \frac{dy}{dt} = R_k(y),  \\[0.2cm]
\displaystyle \frac{ds}{dt} = 1,  \\
\end{array} \right.
\end{equation}
and Eq.(2.23) as
\begin{equation}
\left\{ \begin{array}{l}
\displaystyle \frac{dy}{dt} = R_k(y) + \varepsilon R_{k+1}(y) + \cdots  + \varepsilon ^{m-k}R_m(y)
 + \varepsilon ^{m+1-k} S(s/\varepsilon ^k, y, \varepsilon ),  \\[0.2cm]
\displaystyle \frac{ds}{dt} = 1,  \\
\end{array} \right.
\end{equation}
respectively. Suppose that $m \geq 2k$.
Since $S$ is bounded in $s$ and since
\begin{equation}
\frac{\partial }{\partial y}\varepsilon ^{m+1-k} S(s/\varepsilon ^k, y, \varepsilon ) \sim O(\varepsilon^{k+1})
,\,\, \frac{\partial }{\partial s}\varepsilon ^{m+1-k} S(s/\varepsilon ^k, y, \varepsilon ) \sim O(\varepsilon),
\end{equation}
Eq.(2.33) is $\varepsilon $-close to Eq.(2.32) on the $(s,y)$ space in the $C^1$ topology.

By the assumption, Eq.(2.32) has a normally hyperbolic invariant manifold $\mathbf{R} \times N$ on the $(s,y)$ space.
At this time, Fenichel's theorem is not applicable because $\mathbf{R} \times N$ is not compact.
To handle this difficulty, we do as follows:

Since $S$ is almost periodic, the set
\begin{equation}
T(S, \delta ) 
:= \{ \tau \, | \, || S((s-\tau)/\varepsilon ^k, y, \varepsilon ) - S(s/\varepsilon ^k, y, \varepsilon ) || < \delta , 
\quad \forall s\in \mathbf{R}\}
\end{equation}
is relatively dense for any small $\delta >0$.
Let us fix $\delta $ so that it is sufficiently smaller than $\varepsilon $ and fix $\tau \in T(S, \delta )$ arbitrarily.
Then $W := [0, \tau] \times N$ is a compact locally invariant manifold of Eq.(2.32)
with boundaries $\{0\} \times N$ and $\{\tau\} \times N$ (see Fig.1).

Now Fenichel's theorem proves that Eq.(2.33) has a locally invariant manifold $W_\varepsilon $ which is diffeomorphic to $W$
and lies within an $O(\varepsilon )$ neighborhood of $W$ uniquely.

To extend $W_\varepsilon $ along the $s$ axis, consider the system
\begin{equation}
\left\{ \begin{array}{l}
\displaystyle \dot{y} = R_k(y) + \varepsilon R_{k+1}(y) + \cdots  + \varepsilon ^{m-k}R_m(y)
 + \varepsilon ^{m+1-k} S((s-\tau)/\varepsilon ^k, y, \varepsilon ),  \\[0.2cm]
\displaystyle \dot{s} = 1.  \\
\end{array} \right.
\end{equation}
Since the above system is $\delta $-close to Eq.(2.33), it has a locally invariant manifold $W_{\varepsilon ,\delta }$,
which is diffeomorphic to $W _\varepsilon $.
By putting $\tilde{s} = s- \tau$, Eq.(2.36) is rewritten as
\begin{equation}
\left\{ \begin{array}{l}
\displaystyle \dot{y} = R_k(y) + \varepsilon R_{k+1}(y) + \cdots  + \varepsilon ^{m-k}R_m(y)
 + \varepsilon ^{m+1-k} S(\tilde{s}/\varepsilon ^k, y, \varepsilon ),  \\[0.2cm]
\displaystyle \dot{\tilde{s}} = 1,  \\
\end{array} \right.
\end{equation}
and it takes the same form as Eq.(2.33).
This means that the set
\begin{eqnarray*}
K:= \{ (s,y) \, | \, (s-\tau, y) \in W_{\varepsilon ,\delta }\}
\end{eqnarray*}
is a locally invariant manifold of Eq.(2.33).
Since $W_{\varepsilon , \delta }$ is $\delta $-close to $W_\varepsilon $ and since $\delta  \ll \varepsilon $,
both of $W_\varepsilon \cap \{ s= \tau \}$ and $K \cap \{ s= \tau \}$ are $\varepsilon $-close to $W$.
Since an invariant manifold of Eq.(2.33) which lies within an $O(\varepsilon )$ neighborhood of $W$ is unique
by Fenichel's theorem, $K \cap \{ s= \tau\}$ has to coincide with $W_{\varepsilon} \cap \{ s= \tau\}$.
This proves that $K$ is connected to $W_\varepsilon $ and $K \cup W_\varepsilon $
gives a locally invariant manifold of Eq.(2.33).

This procedure is done for any $\tau \in T(S, \delta )$.
Thus it turns out that $W_\varepsilon $ is extended along the $s$ axis
and it gives an invariant manifold $\tilde{N}_\varepsilon \simeq \mathbf{R} \times N$ of Eq.(2.33).
An invariant manifold $N_\varepsilon $ of Eq.(2.1) is obtained by transforming $\tilde{N}_\varepsilon $ by $\alpha ^{(m)}_t$.

Note that by the construction, projections of the sets $\tilde{N}_{\varepsilon } \cap \{ s = \tau\},\, \tau \in T(S, \delta )$
on to the $y$ space are $\delta$-close to each other. This fact is used to prove the next corollary. \hfill $\blacksquare$

\begin{figure}[h]
\begin{center}
\includegraphics[scale=1.0]{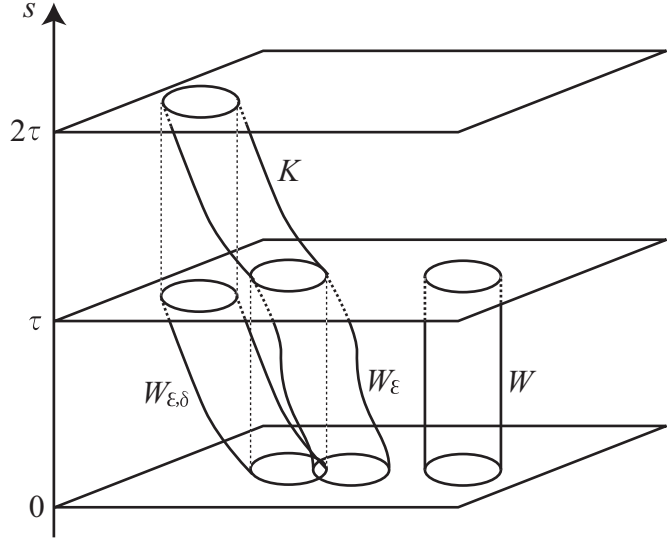}
\end{center}
\caption{A schematic view of the proof for the case that $N$ is a circle.
The $W_{\varepsilon }$ is $\varepsilon $-close to $W$ and $W_{\varepsilon , \delta }$ is $\delta $-close to $W_\varepsilon $.
The $K$ is the ``copy" of $W_{\varepsilon ,\delta }$.}
\end{figure}

The next corollary (ii) and (iii) for $k=1$ are proved in Bogoliubov, Mitropolsky [6] and Fink [20]
and immediately follow from Thm.2.9.
\\[0.2cm]
\textbf{Corollary 2.10.}\, Suppose that $R_1(y) = \cdots = R_{k-1}(y) = 0$ and
$\varepsilon ^k R_k(y)$ is the first non-zero term in the RG equation for Eq.(2.1).
For sufficiently small $\varepsilon >0$,
\\[0.2cm]
(i) \, if the vector field $R_k(y)$ has a hyperbolic periodic orbit $\gamma _0(t)$, then Eq.(2.1)
has an almost periodic solution with the same stability as $\gamma _0(t)$,
\\
(ii) if the vector field $R_k(y)$ has a hyperbolic fixed point $\gamma _0$, then Eq.(2.1)
has an almost periodic solution $\gamma _\varepsilon (t)$ with the same stability as $\gamma _0$ such that 
$\mathrm{Mod} (\gamma _\varepsilon ) \subset \mathrm{Mod} (g)$,
\\
(iii) if the vector field $R_k(y)$ has a hyperbolic fixed point $\gamma _0$ and if $g$ is periodic in $t$
with a period $T$, then Eq.(2.1)
has a periodic solution $\gamma _\varepsilon (t)$ with the same stability as $\gamma _0$ and 
the period $T$ (it need not be the least period).
\\[0.2cm]
\textbf{Proof.}\, If $R_k(y)$ has a periodic orbit, 
Eq.(2.33) has an invariant cylinder $\tilde{N}_{\varepsilon }$ on the $(s,y)$ space as is represented in Fig.1.
To prove Corollary 2.10 (i), at first we suppose that $g(t,x,\varepsilon )$ is periodic with a period $T$.
In this case, since $S(t,y,\varepsilon )$ is a periodic function with the period $T$,
$\tilde{N}_\varepsilon $ is periodic along the $s$ axis in the sense that 
the projections $S^1 := \tilde{N}_{\varepsilon } \cap \{ s = mT\} $
give the same circle for all integers $m$. Let $y= \gamma (t),\, s=t$ be a solution of Eq.(2.33) on the cylinder.
Then $\gamma (mT),\, m= 0,1,\cdots $ gives a discrete dynamics on $S^1$.
If $\gamma (mT)$ converges to a fixed point or a periodic orbit as $m\to \infty$, 
then $\gamma (t)$ converges to a periodic function as $t\to \infty$.
Otherwise, the orbit of $\gamma (mT)$ is dense on $S^1$ and in this case $\gamma (t)$ is an almost periodic function.
A solution of Eq.(2.1) is obtained by transforming $\gamma (t)$ by the almost periodic map $\alpha ^{(m)}_t$.
This proves (i) of Corollary 2.10 for the case that $g$ is periodic.

If $g$ is almost periodic, the sets $\tilde{N}_{\varepsilon } \cap \{ s = \tau \} $
give circles for any $\tau \in T(S, \delta )$ and they are $\delta $-close to each other as is mentioned
in the end of the proof of Thm.2.9.
In this case, there exists a coordinate transformation $Y = \varphi (y, t)$
such that the cylinder $\tilde{N}_{\varepsilon }$ is straightened along the $s$ axis.
The function $\varphi $ is almost periodic in $t$ because $|| \varphi (y,t+ \tau) - \varphi (y, t)||$ is of order $O(\delta )$
for any $\tau \in T(S, \delta )$.
Now the proof is reduced to the case that $g$ is periodic.

The proofs of (ii) and (iii) of Corollary 2.10 are done in the same way as (i), details of which are left to the reader. \hfill $\blacksquare$
\\[0.2cm]
\textbf{Remark 2.11.} \, Suppose that the first order RG equation $\varepsilon R_1(y) \neq 0$ 
does not have normally hyperbolic invariant manifolds but
the second order RG equation $\varepsilon R_1(y) + \varepsilon ^2 R_2(y)$ does.
Then can we conclude that the original system (2.1) has an invariant manifold with the same stability as that of the second order 
RG equation?
Unfortunately, it is not true in general.
For example, suppose that the RG equation for some system is a linear equation of the form
\begin{eqnarray}
\dot{y}/\varepsilon  = \left(
\begin{array}{@{\,}cc@{\,}}
0 & 1 \\
0 & 0 
\end{array}
\right) y - \varepsilon \left(
\begin{array}{@{\,}cc@{\,}}
1 & 0 \\
0 & 1 
\end{array}
\right) y+ \varepsilon ^2 \left(
\begin{array}{@{\,}cc@{\,}}
0 & 0 \\
4 & 0
\end{array}
\right) y+ \cdots ,\,\,\, y\in \mathbf{R}^2.
\end{eqnarray}
The origin is a fixed point of this system, however, the first term has zero eigenvalues and
we can not determine the stability up to the first order RG equation.
If we calculate up to the second order,
the eigenvalues of the matrix
\begin{eqnarray}
\left(
\begin{array}{@{\,}cc@{\,}}
0 & 1 \\
0 & 0 
\end{array}
\right)  - \varepsilon \left(
\begin{array}{@{\,}cc@{\,}}
1 & 0 \\
0 & 1 
\end{array}
\right)
\end{eqnarray}
are $-\varepsilon $ (double root), so that $y=0$ is a stable fixed point of the second order RG equation
if $\varepsilon >0$. 
Unlike Corollary 2.10 (ii),
this \textit{does not} prove that the original system has a stable almost periodic solution.
Indeed, if we calculate the third order RG equation, the eigenvalues of the matrix in the right hand side of Eq.(2.38)
are $3\varepsilon$ and $-\varepsilon $. Therefore the origin is an unstable fixed point
of the third order RG equation.
This example shows that if we truncate higher order terms of the RG equation, stability
of an invariant manifold may change and we can not use $\varepsilon R_1(y) + \varepsilon ^2 R_2(y)$
to investigate stability of an invariant manifold as long as $R_1(y) \neq 0$.
This is because Fenichel's theorem does not hold if the vector field $f$ 
in his theorem depends on the parameter $\varepsilon $.
\\[-0.2cm]

Theorems 2.7 and 2.9 mean that the RG equation is useful to understand the properties of the flow of the system (2.1).
Since the RG equation is an autonomous system while Eq.(2.1) is not, it seems that the RG equation
is easier to analyze than the original system (2.1).
Actually, we can show that the RG equation does not lose symmetries the system (2.1) has.

Recall that integral constants in Eqs.(2.12, 14) are left undetermined and they can depend on $y$ (see Remark 2.4).
To express the integral constants $B _i(y)$ in Eqs.(2.12, 14) explicitly, we rewrite them as
\begin{eqnarray*}
& & u^{(1)}_t(y) = B_1(y) + \int^t \! \left( G_1(s,y) - R_1(y) \right) ds,
\end{eqnarray*}
and 
\begin{eqnarray*}
& & u^{(i)}_t(y)=  B_i(y) + \int^t \! \Bigl( G_i(s, y,u^{(1)}_s(y)
, \cdots , u^{(i-1)}_s(y)) - \sum^{i-1}_{k=1}\frac{\partial u_s^{(k)}}{\partial y}(y) R_{i-k}(y) - R_i(y) \Bigr) ds,
\end{eqnarray*}
for $i= 2,3, \cdots $, where integral constants of the indefinite integrals in the above formulas are chosen to be zero.
\\[0.2cm]
\textbf{Theorem 2.12 (Inheritance of symmetries).}\, 
Suppose that an $\varepsilon $-independent Lie group $H$ acts on $U \subset M$.
If the vector field $g$ and integral constants $B_i(y),\, i= 1, \cdots ,m-1$ in Eqs.(2.12, 14)
are invariant under the action of $H$;
that is, they satisfy
\begin{equation}
g(t, hy, \varepsilon ) = \frac{\partial h}{\partial y}(y)g(t,y,\varepsilon ), \quad
B_i(hy) = \frac{\partial h}{\partial y}(y) B_i(y),
\end{equation}
for any $h\in H, y\in U, t\in \mathbf{R}$ and $\varepsilon $, 
then the $m$-th order RG equation for Eq.(2.1) is also invariant under the action of $H$.
\\[0.2cm]
\textbf{Proof.}\, Since $h\in H$ is independent of $\varepsilon $, Eq.(2.40) implies
\begin{equation}
g_i(t, hy) = \frac{\partial h}{\partial y}(y)g_i(t,y),
\end{equation}
for $i= 1,2,\cdots $.
We prove by induction that $R_i(y)$ and $u^{(i)}_t(y) ,\, i=1,2,\cdots ,$ are invariant under the action of $H$.
At first, $R_1(hy),\, h\in H$ is calculated as
\begin{eqnarray}
R_1(hy) &=& \lim_{t\to \infty} \frac{1}{t}\int^t \! G_1(t,hy) ds \nonumber \\
&=& \lim_{t\to \infty} \frac{1}{t}\int^t \! \frac{\partial h}{\partial y}(y)G_1(t,y) ds = \frac{\partial h}{\partial y}(y)R_1(y).
\end{eqnarray}
Next, $u^{(1)}_t$ is calculated in a similar way:
\begin{eqnarray*}
u^{(1)}_t(hy) &=& B_1(hy) + \int^t \! \left( G_1(s,hy) - R_1(hy) \right) ds \nonumber \\
&=& \frac{\partial h}{\partial y}(y) B_1(y)
    + \frac{\partial h}{\partial y}(y) \int^t \! \left( G_1(s,y) - R_1(y) \right) ds \nonumber \\
&=& \frac{\partial h}{\partial y}(y)u^{(1)}_t(y).
\end{eqnarray*}
Suppose that $R_k$ and $u^{(k)}_t$ are invariant under the action of $H$ for $k=1,2, \cdots , i-1$.
Then, it is easy to verify that
\begin{eqnarray}
& & \frac{\partial u^{(k)}_t}{\partial y}(hy) = \frac{\partial h}{\partial y}(y)\frac{\partial u^{(k)}_t}{\partial y}(y)
\left( \frac{\partial h}{\partial y}(y) \right)^{-1}, \\
& & G_k(hy, u^{(1)}_t(hy) , \cdots , u^{(k-1)}_t(hy))
 = \frac{\partial h}{\partial y}(y) G_k(y, u^{(1)}_t(y) , \cdots , u^{(k-1)}_t(y)),
\end{eqnarray}
for $k= 1,2, \cdots ,i-1$.
These equalities and Eqs.(2.13), (2.14) prove Theorem 2.12 by a similar calculation to Eq.(2.42). \hfill $\blacksquare$

%%%%%%%%%%%%%%%%%%%%%%%%%%%%%%%%%%%%%%%%%%%%%%%%%%%%%%%%%%%%%%%%%%%%%%%%%%%%%%%%%%%%%%%%%%%%%%%%

\subsection{Main theorems for autonomous systems}

In this subsection, we consider an autonomous system of the form
\begin{eqnarray}
\dot{x} &=& f(x) + \varepsilon g(x, \varepsilon ) \nonumber \\
&=& f(x) + \varepsilon g_1(x) + \varepsilon ^2 g_2(x) + \cdots ,\,\, x\in U \subset M,
\end{eqnarray}
where the flow $\varphi _t$ of $f$ is assumed to be almost periodic due to the assumption (B) so that
Eq.(2.45) is transformed into the system of the form of (2.1).
For this system, 
we restate definitions and theorems obtained so far in the present notation for convenience.
We also show a few additional theorems.
\\[0.2cm]
\textbf{Definition 2.13.} \, Let $\varphi _t$ be the flow of the vector field $f$.
For Eq.(2.45), define the $C^\infty$ maps $R_i, h^{(i)}_t : U \to M $ to be
\begin{eqnarray}
& & R_1(y) = \lim_{t\to \infty} \frac{1}{t}
\int^t \! \left( D\varphi _s \right)_y^{-1}G_1(s, \varphi _s(y))ds, \\
& & h^{(1)}_{t}(y) = (D\varphi _t)_y \int^t \! 
\left( (D\varphi _s)_y^{-1}G_1(s, \varphi _s(y)) - R_1(y)\right) ds, 
\end{eqnarray}
and 
\begin{eqnarray}
 R_i(y) &=&  \lim_{t\to \infty} \frac{1}{t} \int^t \! 
\Bigl( (D\varphi _s)_y^{-1} G_i(s, \varphi _s(y), h^{(1)}_s(y) , \cdots , h^{(i-1)}_s(y)) \nonumber \\
& & \quad \quad \quad -(D\varphi _s)_y^{-1}\! \sum^{i-1}_{k=1}(Dh^{(k)}_s)_y R_{i-k}(y) \Bigr) ds, \\
 h^{(i)}_t(y) &=& (D\varphi _t)_y \int^t \! 
\Bigl( (D\varphi _s)_y^{-1} G_i(s, \varphi _s(y), h^{(1)}_s(y) , \cdots , h^{(i-1)}_s(y)) \nonumber \\
& & \quad \quad \quad - (D\varphi _s)_y^{-1}\! \sum^{i-1}_{k=1}(Dh^{(k)}_s)_y R_{i-k}(y) - R_i(y) \Bigr) ds ,
\end{eqnarray}
for $i=2,3, \cdots $, respectively, where $(Dh^{(k)}_t)_y$ is the derivative of 
$h^{(k)}_t(y)$ with respect to $y$, $(D\varphi _t)_y$ is the derivative of 
$\varphi _t(y)$ with respect to $y$, and where $G_i$ are defined through Eq.(2.6).
With these $R_i$ and $h^{(i)}_t$, define the 
\textit{$m$-th order RG equation} for Eq.(2.45) to be
\begin{equation}
\dot{y} = \varepsilon R_1(y) + \varepsilon ^2R_2(y) + \cdots  + \varepsilon ^mR_m(y),
\end{equation}
and define the \textit{$m$-th order RG transformation} to be
\begin{equation}
\alpha _t^{(m)}(y) = \varphi _t(y) + \varepsilon h^{(1)}_t(y) + \cdots  + \varepsilon ^mh^{(m)}_t(y),
\end{equation}
respectively. 
\\[-0.2cm]

In the present notation, Theorems 2.5 and 2.7 are true though the relation $\mathrm{Mod} (S) \subset \mathrm{Mod} (g)$
in Thm.2.5 (ii) is replaced by $\mathrm{Mod} (S) \subset \mathrm{Mod} (\varphi _t)$.
Note that even if Eq.(2.45) is autonomous, the function $S$ depends on $t$ as long as the flow $\varphi _t$ depends on $t$.

Theorem 2.9 is refined as follows:
\\[0.2cm]
\textbf{Theorem 2.14 (Existence of invariant manifolds).}\, \, Suppose that $R_1(y) = \cdots = R_{k-1}(y) = 0$ and
$\varepsilon ^k R_k(y)$ is the first non-zero term in the RG equation for Eq.(2.45).
If the vector field $R_k(y)$ has a boundaryless compact normally hyperbolic invariant manifold $N$,
then for sufficiently small $\varepsilon >0$, 
Eq.(2.45) has an invariant manifold $N_\varepsilon $, which is diffeomorphic to $N$.
In particular, the stability of $N_\varepsilon $ coincides with that of $N$.
\\[-0.2cm]

Note that unlike Thm.2.9, we need not prepare the $(s,x)$ space, 
and the invariant manifold $N_\varepsilon $ lies on $M$ not $\mathbf{R} \times M$.
This theorem immediately follows from the proof of Thm.2.9.
Indeed, Eq.(2.33) has an invariant manifold $\tilde{N}_\varepsilon \simeq \mathbf{R} \times N$ 
on the $(s,y)$ space as is shown in the proof of Thm.2.9.
An invariant manifold of Eq.(2.45) on the $(s,x)$ is obtained by transforming $\tilde{N}_\varepsilon $
by the RG transformation. However, it has to be straight along the $s$ axis because Eq.(2.45) is autonomous.
Thus its projection onto the $x$ space gives the invariant manifold $N_\varepsilon $ of Eq.(2.45) (see Fig.2).

\begin{figure}[h]
\begin{center}
\includegraphics[scale=1.0]{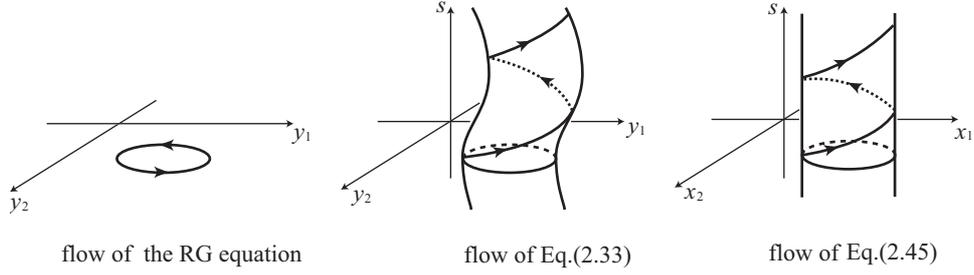}
\end{center}
\caption{A schematic view of the proof for the case that $N$ is a circle. The projection of the straight cylinder
on to the $x$ space gives an invariant manifold of Eq.(2.45).}
\end{figure}

For the case of Eq.(2.1), the RG equation is simpler than the original system (2.1) in the sense that
it has the same symmetries as (2.1) and further it is an autonomous system while (2.1) is not.
In the present situation of (2.45), Theorem 2.12 of inheritance of symmetries still holds as long as the assumption 
for $g$ is replaced as ``the vector field $f$ and $g$ are invariant under the action of a Lie group $H$".
However, since Eq.(2.45) is originally an autonomous system, it is not clear that the RG equation for Eq.(2.45)
is easier to analyze than the original system (2.45).
The next theorem shows that the RG equation for Eq.(2.45) has larger symmetries than Eq.(2.45).

To express the integral constants $B _i(y)$ in Eqs.(2.47, 49) explicitly, we rewrite them as
\begin{eqnarray*}
& & h^{(1)}_{t}(y) = (D\varphi _t)_yB_1(y) + (D\varphi _t)_y \int^t \! 
\left( (D\varphi _s)_y^{-1}G_1(s, \varphi _s(y)) - R_1(y)\right) ds, 
\end{eqnarray*}
and 
\begin{eqnarray*}
h^{(i)}_t(y) &=& (D\varphi _t)_yB_i(y) + (D\varphi _t)_y \int^t \! 
\Bigl( (D\varphi _s)_y^{-1} G_i(s, \varphi _s(y), h^{(1)}_s(y) , \cdots , h^{(i-1)}_s(y)) \nonumber \\
& & \quad \quad \quad - (D\varphi _s)_y^{-1}\! \sum^{i-1}_{k=1}(Dh^{(k)}_s)_y R_{i-k}(y) - R_i(y) \Bigr) ds ,
\end{eqnarray*}
for $i=2,3, \cdots $, where integral constants of the indefinite integrals in the above formulas are chosen to be zero.
\\[0.2cm]
\textbf{Theorem 2.15 (Additional symmetry).}\, Let $\varphi _t$ be the flow of the vector field $f$ defined on $U\subset M$.
If the integral constants $B_i$ in Eqs.(2.47, 49) are chosen so that they are 
invariant under the action of the one-parameter group $\{ \varphi _t : U\to M \, | \, t\in \mathbf{R}\}$,
then the RG equation for Eq.(2.45) is also invariant under the action of the group.
In other words, $R_i$ satisfies the equality
\begin{equation}
R_i(\varphi _t(y)) = (D\varphi _t)_y R_i(y),
\end{equation}
for $i=1,2, \cdots $.
\\[0.2cm]
\textbf{Proof.}\, Since Eq.(2.45) is autonomous, the function $G_k(t, x_0, \cdots ,x_{k-1})$ defined through Eq.(2.6)
is independent of $t$ and we write it as $G_k(x_0, \cdots ,x_{k-1})$.
We prove by induction that equalities $R_i(\varphi _t(y)) = (D\varphi _t)_yR_i(y)$
 and $h^{(i)}_t(\varphi _{t'}(y)) = h^{(i)}_{t+t'}(y)$ hold for $i=1,2,\cdots $.
For all $s' \in \mathbf{R}$, $R_1(\varphi _{s'}(y))$ takes the form
\begin{eqnarray*}
R_1(\varphi _{s'}(y)) 
 &=& \lim_{t\to \infty}\frac{1}{t}\int^t \! (D\varphi _s)_{\varphi _{s'}(y)}^{-1}G_1(\varphi _s \circ \varphi _{s'}(y))ds \\
 &=& (D\varphi _{s'})_y\lim_{t\to \infty}\frac{1}{t}\int^t \! (D\varphi _{s+s'})_y^{-1} G_1(\varphi _{s+s'}(y))ds.  
\end{eqnarray*}
Putting $s+s' = s''$, we verify that
\begin{eqnarray*}
R_1(\varphi _{s'}(y)) 
 &=& (D\varphi _{s'})_y \lim_{t\to \infty}\frac{1}{t}\int^{t+s'} \! (D\varphi _{s''})_y^{-1}G_1(\varphi _{s''}(y))ds'' \\
 &=& (D\varphi _{s'})_yR_1(y) 
  + (D\varphi _{s'})_y \lim_{t\to \infty}\frac{1}{t}\int^{t+s'}_t \! (D\varphi _{s''})_y^{-1}G_1(\varphi _{s''}(y))ds''\\
&=& (D\varphi _{s'})_yR_1(y).
\end{eqnarray*}
The $h^{(1)}_t(\varphi _{s'}(y))$ is calculated in a similar way as
\begin{eqnarray*}
h^{(1)}_t(\varphi _{s'}(y)) &=&  (D\varphi _{t})_y B_1(\varphi _{s'}(y))
 + (D\varphi _{t})_{\varphi _{s'}(y)}\int^t \! \left( 
   (D\varphi _{s})_{\varphi _{s'}(y)}^{-1}G_1(\varphi _s\circ \varphi _{s'}(y)) - R_1(\varphi _{s'}(y)) \right) ds \\
&=& (D\varphi _{t+s'})_y B_1(y)
 + (D\varphi _{t+s'})_y \int^t \! \left( (D\varphi _{s+s'})_y^{-1}G_1(\varphi _{s+s'}(y)) - R_1(y) \right) ds.
\end{eqnarray*}
Putting $s+s' = s''$ provides
\begin{eqnarray}
h^{(1)}_t(\varphi _{s'}(y)) &=& (D\varphi _{t+s'})_y B_1(y)
+ (D\varphi _{t+s'})_y\int^{t+s'} \! 
  \left( (D\varphi _{s''})_y^{-1}G_1(\varphi _{s''}(y)) - R_1(y) \right) ds'' \nonumber \\
 &=& h^{(1)}_{t+s'}(y).
\end{eqnarray}
Suppose that $R_k(\varphi _t(y)) = (D\varphi _{t})_yR_k(y)$ and 
$h^{(k)}_t(\varphi _{t'}(y)) = h^{(k)}_{t+t'}(y)$ hold for $k=1,2,\cdots ,i-1$.
Then, $R_i(\varphi _{s'}(y))$ is calculated as
\begin{eqnarray*}
R_i(\varphi _{s'}(y)) &=& \lim_{t\to \infty} \frac{1}{t} \int^t \! 
\Bigl( (D\varphi _{s})_{\varphi _{s'}(y)}^{-1} G_i(\varphi _s \circ \varphi _{s'}(y), h^{(1)}_s(\varphi _{s'}(y)) ,
   \cdots , h^{(i-1)}_s(\varphi _{s'}(y))) \nonumber \\
& & \quad \quad \quad - (D\varphi _{s})_{\varphi _{s'}(y)}^{-1}\! 
   \sum^{i-1}_{k=1}(Dh^{(k)}_s)_{\varphi _{s'}(y)} R_{i-k}(\varphi _{s'}(y)) \Bigr) ds \\
&=& (D\varphi _{s'})_y\lim_{t\to \infty} \frac{1}{t} \int^t \! 
\Bigl( (D\varphi _{s+s'})_y^{-1} G_i(\varphi _{s+s'}(y), h^{(1)}_{s+s'}(y) , \cdots , h^{(i-1)}_{s+s'}(y)) \nonumber \\
& & \quad \quad \quad - (D\varphi _{s+s'})_y^{-1}\! \sum^{i-1}_{k=1}(Dh^{(k)}_{s+s'})_{y} R_{i-k}(y) \Bigr) ds.
\end{eqnarray*}
Putting $s+s' = s''$ provides
\begin{eqnarray*}
R_i(\varphi _{s'}(y)) &=&(D\varphi _{s'})_y\lim_{t\to \infty} \frac{1}{t} \int^{t+s'} \! 
\Bigl( (D\varphi _{s''})_y^{-1} G_i(\varphi _{s''}(y), h^{(1)}_{s''}(y) , \cdots , h^{(i-1)}_{s''}(y)) \nonumber \\
& & \quad \quad \quad - (D\varphi _{s''})_y^{-1}\! \sum^{i-1}_{k=1}(Dh^{(k)}_{s''})_{y} R_{i-k}(y) \Bigr) ds'' \\
&=& (D\varphi _{s'})_yR_i(y) + (D\varphi _{s'})_y\lim_{t\to \infty} \frac{1}{t} \int^{t+s'}_t \! 
\Bigl( (D\varphi _{s''})_y^{-1} G_i(\varphi _{s''}(y), h^{(1)}_{s''}(y) , \cdots , h^{(i-1)}_{s''}(y)) \nonumber \\
& & \quad \quad \quad - (D\varphi _{s''})_y^{-1}\! \sum^{i-1}_{k=1}(Dh^{(k)}_{s''})_{y} R_{i-k}(y) \Bigr) ds''\\
&=& (D\varphi _{s'})_yR_i(y).
\end{eqnarray*}
We can show the equality $h^{(i)}_t(\varphi _{t'}(y)) = h^{(i)}_{t+t'}(y)$ in a similar way. \hfill $\blacksquare$

%%%%%%%%%%%%%%%%%%%%%%%%%%%%%%%%%%%%%%%%%%%%%%%%%%%%%%%%%%%%%%%%%%%%%%%%%%%%%%%%%%%%%%%%%%%%%%%%%%%%%%%%%%%%%%%

\subsection{Simplified RG equation}

Recall that the definitions of $R_i$ and $u^{(i)}_t$ given in Eqs.(2.11) to (2.14) include the indefinite
integrals and we have left the integral constants undetermined (see Rem.2.4).
In this subsection, we investigate how different choices of the integral constants change the forms of RG equations
and RG transformations.

Let us fix definitions of $R_i$ and $u^{(i)}_t,\, i=1,2,\cdots $ by fixing
integral constants in Eqs.(2.12, 14) arbitrarily
(note that $R_i$ and $u^{(i)}_t$ are independent of the integral constants in Eqs.(2.11, 13)).
For these $R_i$ and $u^{(i)}_t$, we define $\widetilde{R}_i$ and $\widetilde{u}^{(i)}_t$ to be
\begin{eqnarray}
& & \widetilde{R}_1(y) = \lim_{t\to \infty}\frac{1}{t} \int ^t\! G_1(s,y)ds, \label{2_4_1} \\
& & \widetilde{u}^{(1)}_t(y) = B_1(y) + \int^t \! \left( G_1(s,y) - \widetilde{R}_1(y) \right) ds, \label{2_4_2}
\end{eqnarray}
and
\begin{eqnarray}
& & \widetilde{R}_i(y)  
= \lim_{t\to \infty}\frac{1}{t} \int ^t\! \Bigl(
G_i(s, y, \widetilde{u}^{(1)}_s(y),  \cdots, \widetilde{u}^{(i-1)}_s(y))- 
\sum^{i-1}_{k=1}\frac{\partial \widetilde{u}_s^{(k)}}{\partial y}(y)\widetilde{R}_{i-k}(y) \Bigr) ds , \label{2_4_3} \\
& & \widetilde{u}^{(i)}_t(y)=  B_i(y) + \int^t \! \Bigl( G_i(s, y,\widetilde{u}^{(1)}_s(y)
, \cdots , \widetilde{u}^{(i-1)}_s(y))
 - \sum^{i-1}_{k=1}\frac{\partial \widetilde{u}_s^{(k)}}{\partial y}(y) \widetilde{R}_{i-k}(y)
   - \widetilde{R}_i(y) \Bigr) ds \label{2_4_4},
\end{eqnarray}
for $i=2,3,\cdots$, respectively, where integral constants in the definitions of $\widetilde{u}_t^{(i)}$ are the same 
as those of $u_t^{(i)}$ and where $B_i,\, i= 1,2, \cdots $ are arbitrary vector fields which imply other choice of 
integral constants. Along with these functions, We define the $m$-th order RG equation and the 
$m$-th order RG transformation to be
\begin{equation}
\dot{y} = \varepsilon \widetilde{R}_1(y) + \cdots  + \varepsilon ^m \widetilde{R}_m(y),\label{2_4_5}
\end{equation}
and 
\begin{equation}
\widetilde{\alpha} ^{(m)}_t (y) 
   = y + \varepsilon \widetilde{u}_t^{(1)}(y) + \cdots  + \varepsilon ^m \widetilde{u}_t^{(m)}(y), \label{2_4_6}
\end{equation}
respectively. 
Now we have two pairs of RG equations-transformations, Eqs.(2.18),(2.19) and Eqs.(\ref{2_4_5}),(\ref{2_4_6}).
Main theorems described so far hold for both of them for any choices of
$B_i(y)$'s except to Thms. 2.11 and 2.15, in which we need additional assumptions for $B_i(y)$'s as was stated.

Let us examine relations between $R_i,\, u_t^{(i)}$ and $\widetilde{R}_i,\, \widetilde{u}_t^{(i)}$.
Clearly $\widetilde{R}_1$ coincides with $R_1$.
Thus $\widetilde{u}_t^{(1)}$ is given as $\widetilde{u}_t^{(1)}(y) =u_t^{(1)}(y) + B_1(y)$.
According to the definition of $G_2$ (Eq.(2.8)), $\widetilde{R}_2$ is calculated as
\begin{eqnarray}
\widetilde{R}_2(y) &=& \lim_{t\to \infty} \frac{1}{t} \int^t \! 
\left( \frac{\partial g_1}{\partial y}(s,y) \widetilde{u}_s^{(1)}(y) + g_2(s,y)
 - \frac{\partial \widetilde{u}_s^{(1)}}{\partial y} \widetilde{R}_1(y) \right) ds \nonumber \\
&=&  \lim_{t\to \infty} \frac{1}{t} \int^t \! 
\left( \frac{\partial g_1}{\partial y}(s,y) u_s^{(1)}(y) + g_2(s,y) 
- \frac{\partial u_s^{(1)}}{\partial y} R_1(y) \right) ds \nonumber \\
& & \quad + \lim_{t\to \infty} \frac{1}{t} \int^t \! 
\left( \frac{\partial g_1}{\partial y}(s,y) B_1(y) - \frac{\partial B_1}{\partial y}(y) R_1(y) \right) ds \nonumber \\
&=& R_2(y) + \frac{\partial R_1}{\partial y}(y) B_1(y) - \frac{\partial B_1}{\partial y}(y) R_1(y). 
\label{2_4_7}
\end{eqnarray}
If we define the commutator $[\, \bm{\cdot} \,,\, \bm{\cdot} \,]$ of vector fields to be
\begin{equation}
[B_1, R_1] (y) = \frac{\partial B_1}{\partial y}(y) R_1(y) - \frac{\partial R_1}{\partial y}(y) B_1(y),
\label{2_4_8}
\end{equation}
Eq.(\ref{2_4_7}) is rewritten as
\begin{equation}
\widetilde{R}_2(y) = R_2(y) - [B_1, R_1] (y). 
\label{2_4_9}
\end{equation}
Similar calculation proves that $\widetilde{R}_i,\, i= 2,3, \cdots $ are expressed as
\begin{equation}
\widetilde{R}_i(y) = R_i(y) + P_i(R_1, \cdots , R_{i-1}, B_1, \cdots , B_{i-2})(y) - [ B_{i-1}, R_1](y),
\label{2_4_10}
\end{equation}
where $P_i$ is a function of $R_1, \cdots , R_{i-1}$ and $ B_1, \cdots , B_{i-2}$.
See Chiba[11] for the proof.
Thus appropriate choices of vector fields $B_i(y),\, i=1,2, \cdots $ may simplify $\widetilde{R}_i ,\, i= 2,3, \cdots $
through Eq.(\ref{2_4_10}).

Suppose that $R_i,\, i= 1,2, \cdots $ are elements of some finite dimensional vector space $V$ and the commutator 
$[\, \bm{\cdot} \, , R_1]$
defines the linear map on $V$. For example if the function $g(t, x, \varepsilon )$ in Eq.(2.1) is polynomial in $x$,
$V$ is the space of polynomial vector fields. If Eq.(2.1) is an $n$-dimensional linear equation,
then $V$ is the space of all $n \times n$ constant matrices.
Let us take a complementary subspace $\mathcal{C}$ to $\mathrm{Im}\, [\, \bm{\cdot} \, , R_1]$ into $V$ arbitrarily:
$V = \mathrm{Im}\, [\, \bm{\cdot} \, , R_1] \bigoplus \mathcal{C}$.
Then, there exist $B_i \in V,\, i = 1, 2, \cdots , m-1$ such that $\widetilde{R}_i \in \mathcal{C}$ 
for $i= 2,3, \cdots , m$ because of Eq.(\ref{2_4_10}).
If $\widetilde{R}_i \in \mathcal{C}$ for $i= 2,3, \cdots , m$, we call Eq.(\ref{2_4_5}) the \textit{$m$-th order simplified RG equation}.
See Chiba[11] for explicit forms of the simplified RG equation for the cases that Eq.(2.1) is polynomial in $x$ or 
a linear system. In particular, they are quite related to the simplified normal forms theory (hyper-normal forms
theory) [2,39,40]. See also Section 4.3. 

Since integral constants $B_i$'s in Eqs.(\ref{2_4_2}) and (\ref{2_4_4}) are independent of $t$,
we can show the next claim, which will be used to prove Thm.5.1.
\\[0.2cm]
\textbf{Claim 2.16.} \, RG equations and RG transformations are not unique in general because of
undetermined integral constants in Eqs.(2.12) and (2.14).
Let $\alpha ^{(m)}_t$ and $\widetilde{\alpha }_t^{(m)}$ be two different RG transformations
for a given system (2.1).
Then, there exists a \textit{time-independent} transformation $\phi (y, \varepsilon )$,
which is $C^\infty$ with respect to $y$ and $\varepsilon $, such that 
$\widetilde{\alpha }_t^{(m)}(y) = \alpha ^{(m)}_t \circ \phi (y, \varepsilon )$.
Conversely, $\alpha ^{(m)}_t \circ \phi (y, \varepsilon )$ gives one of the RG transformations
for any $C^\infty$ maps $\phi (y, \varepsilon )$.
\\[-0.2cm]

Note that the map $\phi (y, \varepsilon )$ is independent of $t$ because it brings one of the RG equations
into the other RG equation, both of which are autonomous systems.
According to Claim 2.16, one of the simplest way to achieve simplified RG equations is as follows:
At first, we calculate $R_i$ and $u^{(i)}_t$ by fixing integral constants arbitrarily and obtain the RG equation.
It may be convenient in practice to choose zeros as integral constants (see Prop.2.18 below).
Then, any other RG equations are given by transforming the present RG equation by time-independent $C^\infty$ maps.

%%%%%%%%%%%%%%%%%%%%%%%%%%%%%%%%%%%%%%%%%%%%%%%%%%%%%%%%%%%%%%%%%%%%%%%%%%%%%%%%%%%%%%%%%%%%%%%%%%%%%%%%%%%%%%%%

\subsection{An example}

In this subsection, we give an example to verify the main theorems.
See Chiba[10,11] for more examples.
Consider the system on $\mathbf{R}^2$
\begin{equation}
\left\{ \begin{array}{l}
\dot{X}_1 = X_2 + X_2^2 + \varepsilon^2 k \sin (\omega t),  \\
\dot{X}_2 = -X_1 + \varepsilon^2 X_2 - X_1X_2 + X_2^2 , \\
\end{array} \right.
\label{ex1}
\end{equation}
where $\varepsilon >0,\, k \geq 0$ and $\omega >0$ are parameters.
Changing the coordinates by $(X_1, X_2) = (\varepsilon x_1, \varepsilon x_2)$ yields
\begin{equation}
\left\{ \begin{array}{l}
\dot{x}_1 = x_2 + \varepsilon x_2^2 + \varepsilon k\sin (\omega t),  \\
\dot{x}_2 = -x_1 + \varepsilon (x_2^2 - x_1x_2) + \varepsilon ^2 x_2 . \\
\end{array} \right.
\label{ex2}
\end{equation}
Diagonalizing the unperturbed term by introducing the complex variable $z$ as
$x_1 = z + \overline{z},\, x_2 = i(z- \overline{z})$ may simplify our calculation :
\begin{equation}
\left\{ \begin{array}{l}
\displaystyle \dot{z} = iz + \frac{\varepsilon }{2}\left( i(z-\overline{z})^2 -2z^2 + 2z\overline{z} + k \sin (\omega t)\right)
 + \frac{\varepsilon ^2}{2}(z- \overline{z}),  \\
\displaystyle \dot{\overline{z}} = -i\overline{z} 
   + \frac{\varepsilon }{2}\left( -i(z-\overline{z})^2 -2\overline{z}^2 + 2z\overline{z}  + k \sin (\omega t)\right)
 - \frac{\varepsilon ^2}{2}(z- \overline{z}),  \\
\end{array} \right.
\label{ex3}
\end{equation}
where $i = \sqrt{-1}$.
Let us calculate the RG equation for the system (\ref{ex2}) or (\ref{ex3}).
In this example, all integral constants in Eqs.(2.47, 49) are chosen to be zero.
\\[0.2cm]
\textbf{(i)} When $\omega \neq 1, 2$, the second order RG equation for Eq.(\ref{ex3}) is given as 
\begin{equation}
\left\{ \begin{array}{l}
\displaystyle \dot{y}_1 = \frac{1}{2}\varepsilon ^2 (y_1 - 3y_1^2 y_2 - \frac{16i}{3}y_1^2y_2),  \\[0.2cm]
\displaystyle \dot{y}_2 = \frac{1}{2}\varepsilon ^2 (y_2 - 3y_1 y_2^2 + \frac{16i}{3}y_1y_2^2),  \\
\end{array} \right.
\label{ex4}
\end{equation}
where the first order RG equation $R_1$ vanishes.
Note that it is independent of the time periodic external force $k \sin (\omega t)$.
Thus this RG equation coincides with that of the autonomous system obtained by putting $k=0$ in Eq.(\ref{ex3}) and
Theorem 2.15 is applicable to Eq.(\ref{ex4}).
Indeed, since the above RG equation is invariant under the rotation group $(y_1, y_2) \mapsto (e^{i\tau} y_1, e^{-i\tau}y_2)$,
putting $y_1 = re^{i\theta },\, y_2 = re^{-i\theta }$ results in
\begin{equation}
\left\{ \begin{array}{l}
\displaystyle \dot{r} = \frac{1}{2}\varepsilon ^2r (1-3r^2),  \\[0.1cm]
\displaystyle \dot{\theta } = -\frac{8}{3}\varepsilon ^2 r^2,  \\
\end{array} \right.
\label{ex5}
\end{equation}
and it is easily solved.
We can verify that this RG equation has a stable periodic orbit $r = \sqrt{1/3}$ if $\varepsilon >0$.
Now Corollary 2.10 (i) proves that the original system (\ref{ex2})
has a stable almost periodic solution if $\varepsilon >0$ is small (see Fig.3).
If $k= 0$ and Eq.(\ref{ex2}) is autonomous, then Thm.2.14 is applied to conclude that Eq.(\ref{ex2})
has a stable periodic orbit.

\begin{figure}[h]
\begin{center}
\includegraphics[scale=1.0]{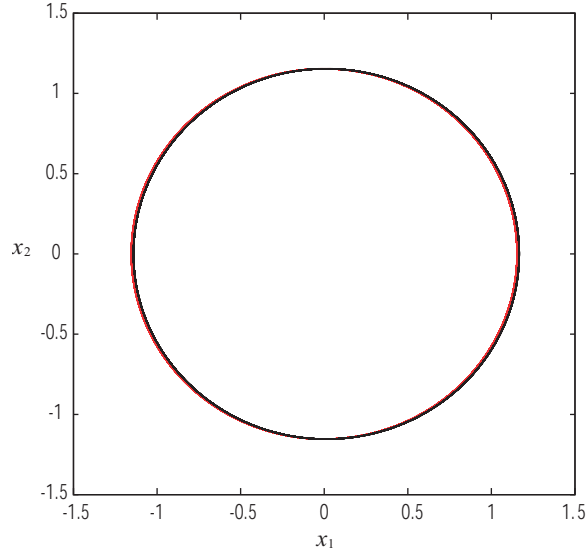}
\end{center}
\caption{Numerical results of the system (\ref{ex2}) and its RG equation (\ref{ex5}) for $\omega = 3,\, k = 1.8$
and $\varepsilon =0.01$.
The red curve denotes the stable periodic orbit of the RG equation and the black curve denotes the 
almost periodic solution of Eq.(\ref{ex2}). They almost overlap with one another.}
\end{figure}

\vspace*{0.2cm}
\noindent \textbf{(ii)} When $\omega = 2$, we prefer Eq.(\ref{ex2}) to Eq.(\ref{ex3}) 
to avoid complex numbers. The second order RG equation for Eq.(\ref{ex2}) is given by
\begin{equation}
\left\{ \begin{array}{l}
\displaystyle \dot{y}_1 
  = \frac{\varepsilon ^2}{24} \left( 12y_1 - 9 y_1^3 - 16 y_1^2 y_2 - 9 y_1 y_2^2 - 16 y_2^3 - k(6y_1 + 4y_2) \right) , \\[0.2cm]
\displaystyle \dot{y}_2 
  = \frac{\varepsilon ^2}{24} \left( 12y_2 - 9 y_2^3 + 16 y_1 y_2^2 - 9 y_1^2 y_2 + 16 y_1^3 - k(4y_1 - 6y_2) \right).  \\
  \end{array} \right.
\label{ex6}
\end{equation}
Since it depends on the time periodic term $k \sin (\omega ^2 t)$, Thm.2.15 is no longer applicable and
to analyze this RG equation is rather difficult.
However, numerical simulation shows that Eq.(\ref{ex6}) undergoes a typical homoclinic bifurcations 
(see Chow, Li and Wang [14]) whose phase portraits are given as Fig.4.

\begin{figure}[h]
\begin{center}
\includegraphics[scale=1.0]{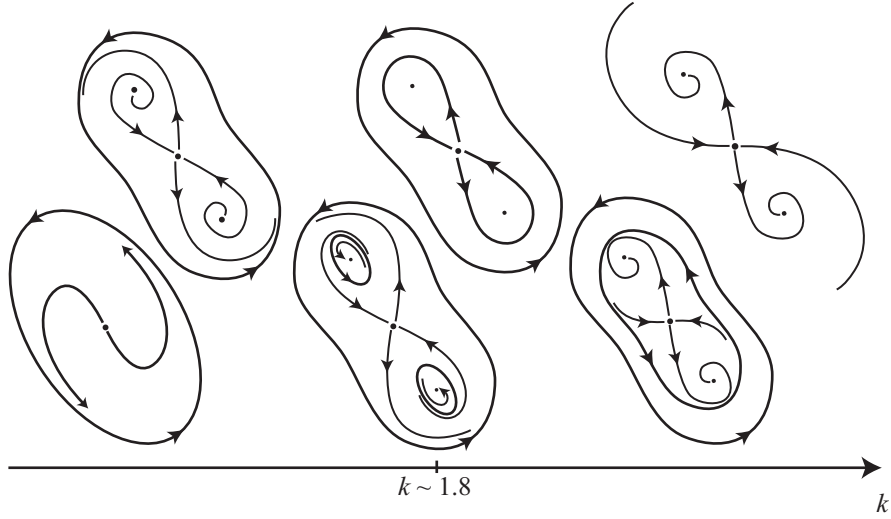}
\end{center}
\caption{Phase portraits of the RG equation (\ref{ex6}).}
\end{figure}

Let us consider the case $k= 1.8$.
In this case, the RG equation has one stable periodic orbit and two stable fixed points.
Thus Corollary 2.10 proves that the original system (\ref{ex2}) has one almost periodic solution $\gamma (t)$ and
two periodic solutions with the period $\pi$ (see Fig.5; in the arXiv version, Fig.5 is omitted because of the limitation of size).

\vspace*{0.2cm}
\noindent \textbf{(iii)} When $\omega = 1$, the second RG equation for Eq.(\ref{ex2}) is given by
\begin{equation}
\left\{ \begin{array}{l}
\displaystyle \dot{y}_1 = \frac{\varepsilon ^2}{24} \left( 12y_1 - 9 y_1^3 - 16 y_1^2 y_2 - 9 y_1 y_2^2 - 16 y_2^3 \right) , \\[0.2cm]
\displaystyle \dot{y}_2 
  = \frac{\varepsilon k}{2} + \frac{\varepsilon ^2}{24} \left( 12y_2 - 9 y_2^3 + 16 y_1 y_2^2 - 9 y_1^2 y_2 + 16 y_1^3 \right).  \\
\end{array} \right.
\label{ex7}
\end{equation}
In this case, the first order RG equation does not vanish.
For small $\varepsilon $, Eq.(\ref{ex7}) has a stable fixed point $\bm{y}_0 = \bm{y}_0 (k, \varepsilon )$, 
which tends to infinity as $\varepsilon  \to 0$. For example if $k= 1.8$ and $\varepsilon = 0.01$,
$\bm{y}_0$ is given by $\bm{y}_0 \sim (-4.35, 2.31)$ (see Fig.6).
Therefore, the original system has a stable periodic orbit whose radius tends to infinity as $\varepsilon \to 0$.  

\setcounter{figure}{5}
\begin{figure}[h]
\begin{center}
\includegraphics[scale=1.0]{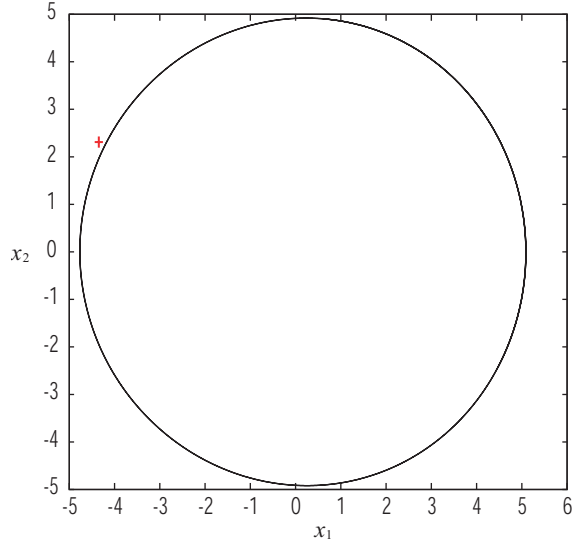}
\end{center}
\caption{Numerical results of the system (\ref{ex2}) and its RG equation (\ref{ex7}) for $\omega = 1,\, k = 1.8$
and $\varepsilon  = 0.01$.
The red cross point denotes the
stable fixed point of the RG equation. 
The black curve denotes the periodic solution of Eq.(\ref{ex2}).}
\end{figure}

%%%%%%%%%%%%%%%%%%%%%%%%%%%%%%%%%%%%%%%%%%%%%%%%%%%%%%%%%%%%%%%%%%%%%%%%%%%%%%%%%%%%%%%%%%%%%%%%%%%%%%%%%%%%%%%

\subsection{A few remarks on symbolic computation}

One can calculate the RG equation and the RG transformation by using symbolic computation softwares, such as
Mathematica and Maple, with formulas (2.11) to (2.14).
In this subsection, we provide a few remarks which may be convenient for symbolic computation.

It is not comfortable to compute the limits in formulas (2.11) and (2.13) directly by symbolic computation
softwares because it takes too much time.
The next proposition is useful to compute them.
\\[0.2cm]
\textbf{Proposition 2.17.} \, Suppose that $F(t,y)$ and its primitive function are almost periodic functions in $t$.
Then, $\lim_{t\to \infty} \frac{1}{t} \int^t \! F(t,y) dt$ gives a coefficient of a linear term of $\int^t \! F(t,y) dt$
with respect to $t$.

This proposition is easily proved because $F$ is expanded in a Fourier series as $F(t,y) \sim \sum a_n(y) e^{i\lambda _n t}$.
Thus, to obtain $R_i(y)$'s, we compute integrals in Eqs.(2.11,13) and extract linear terms 
with respect to $t$. To do so, for example in Mathematica, the command 
\texttt{Coefficient[Integrate[F[t,y],t],t]} is available, where \texttt{F[t,y]} is the integrand in Eqs.(2.11, 13).

When computing the integrals in Eqs.(2.11) to (2.14), we recommend that the integrands are expressed
by exponential functions with respect to $t$ such as $e^{i\lambda t}$ if they include trigonometric functions such as 
$\sin \lambda t,\, \cos \lambda t$ (in Mathematica, it is done by using the command \texttt{TrigToExp}).
It is because if they include trigonometric functions, softwares may choose unexpected 
integral constants while if they consist of exponential functions,
then zeros are chosen as integral constants.
If all integral constants in Eqs.(2.12), (2.14) are zeros, the second term in the integrand in Eq.(2.13)
does not include a constant term with respect to $t$.
Thus we obtain the next proposition, which reduces the amount of calculation.
\\[0.2cm]
\textbf{Proposition 2.18.}\, If we choose zeros as integral constants in the formulas (2.11) to (2.14)
for $i = 1, \cdots  , k-1$, then Eq.(2.13) for $i=k$ is written as
\begin{eqnarray}
& & R_k(y)  
= \lim_{t\to \infty}\frac{1}{t} \int ^t\! G_k(s, y, u^{(1)}_s(y), \cdots, u^{(k-1)}_s(y)) ds.
\end{eqnarray}

If a system for which we calculate the RG equation includes parameters such as $k$ and $\omega $ in Eq.(\ref{ex2}),
then softwares automatically assume that it is in a generic case.
For example if we compute the RG equation for Eq.(\ref{ex2}) by using Mathematica, it is assumed that
$\omega \neq 0, 1, 2$.
We can verify that the RG equation (\ref{ex4}) obtained by softwares is invalid when $\omega = 0, 1, 2$ 
by computing the second order RG transformation. Indeed, it
includes the factors $\omega , \omega -1, \omega - 2$ in denominators.
The explicit form of the RG transformation for Eq.(\ref{ex2}) is too complicated to show here.
To obtain the RG equation for $\omega = 0, 1, 2$, substitute them into the original system (\ref{ex2}) and compute
the RG equation again for each case.

%%%%%%%%%%%%%%%%%%%%%%%%%%%%%%%%%%%%%%%%%%%%%%%%%%%%%%%%%%%%%%%%%%%%%%%%%%%%%%%%%%%%%%%%%%%%%

\section{Restricted RG method}

Suppose that Eq.(2.3) defined on $U \subset M$ does not satisfy the assumption (B) on $U$ but satisfies it on a
submanifold $N_0 \subset U$. Then, the RG method discussed in the previous section is still valid if domains of the 
RG equation and the RG transformation are restricted to $N_0$.
This method to construct an approximate flow on some restricted region is called the \textit{restricted RG method}
and it gives extension of the center manifold theory, the geometric singular perturbation method
and the phase reduction.

\subsection{Main results of the restricted RG method}

Consider the system of the form
\begin{equation}
\dot{x} = f(x) + \varepsilon g(t, x, \varepsilon ),
\end{equation}
defined on $U \subset M$. For this system, we suppose that
\\[0.2cm]
\textbf{(D1)}\, the vector field $f$ is $C^\infty$ and it has a compact attracting normally hyperbolic invariant manifold
$N_0 \subset U$. The flow $\varphi _t (x)$ of $f$ on $N_0$ 
is an almost periodic function with respect to $t \in \mathbf{R}$
uniformly in $x\in N_0$, the set of whose Fourier exponents has no accumulation points.
\\
\textbf{(D2)}\, there exists an open set $V \supset N_0$ in $U$ such that the vector field $g$ is $C^1$ in $t\in \mathbf{R}$,
$C^\infty$ in $x\in V$ and small $\varepsilon $, and that $g$ is an almost periodic function with respect to $t\in \mathbf{R}$
uniformly in $x\in V$ and small $\varepsilon $, the set of whose Fourier exponents has no accumulation points
(i.e. the assumption (A) is satisfied on $V$).
\\[0.2cm]
\indent If $g$ is independent of $t$ and Eq.(3.1) is autonomous, Fenichel's theorem proves that Eq.(3.1)
has an attracting invariant manifold $N_\varepsilon $ near $N_0$.
If $g$ depends on $t$, we rewrite Eq.(3.1) as
\begin{equation}
\left\{ \begin{array}{l}
\dot{x} = f(x) + \varepsilon g(s, x, \varepsilon ), \\
\dot{s} = 1, \\
\end{array} \right.
\end{equation}
so that the unperturbed term $(f(x),1)$ has an attracting normally hyperbolic invariant manifold $\mathbf{R} \times N_0$
on the $(s,x)$ space. Then, a similar argument to the proof of Thm.2.9 proves that Eq.(3.2) has an attracting
invariant manifold $N_\varepsilon $ on the $(s,x)$ space, which is diffeomorphic to $\mathbf{R} \times N_0$.

In both cases, since $N_\varepsilon $ is attracting, long time behavior of the flow of Eq.(3.1) is well described
by the flow on $N_\varepsilon $.
Thus a central issue is to construct $N_\varepsilon $ and the flow on it approximately.
To do so, we establish the restricted RG method on $N_0$.
\\[0.2cm]
\textbf{Definition 3.1.} \, For Eq.(3.1), we define $C^\infty$ maps $R_i,\, h^{(i)}_t : N_0 \to M$ to be
\begin{eqnarray}
& & R_1(y) = \lim_{t\to -\infty} \frac{1}{t}
\int^t \! \left( D\varphi _s \right)_y^{-1}G_1(s, \varphi _s(y))ds, \\
& & h^{(1)}_{t}(y) = (D\varphi _t)_y \int^t_{-\infty} \! 
\left( (D\varphi _s)_y^{-1}G_1(s, \varphi _s(y)) - R_1(y)\right) ds, 
\end{eqnarray}
and 
\begin{eqnarray}
 R_i(y) &=&  \lim_{t\to -\infty} \frac{1}{t} \int^t \! 
\Bigl( (D\varphi _s)_y^{-1}G_i(s, \varphi _s(y), h^{(1)}_s(y) , \cdots , h^{(i-1)}_s(y)) \nonumber \\
& & \quad \quad \quad -(D\varphi _s)_y^{-1}\! \sum^{i-1}_{k=1}(Dh^{(k)}_s)_y R_{i-k}(y) \Bigr) ds, \\
 h^{(i)}_t(y) &=& (D\varphi _t)_y \int^t_{-\infty} \! 
\Bigl( (D\varphi _s)_y^{-1} G_i(s, \varphi _s(y), h^{(1)}_s(y) , \cdots , h^{(i-1)}_s(y)) \nonumber \\
& & \quad \quad \quad - (D\varphi _s)_y^{-1}\! \sum^{i-1}_{k=1}(Dh^{(k)}_s)_y R_{i-k}(y) - R_i(y) \Bigr) ds ,
\end{eqnarray}
for $i=2,3, \cdots $, respectively.
Note that $\lim_{t\to \infty}$
in Eqs.(2.46, 48) are replaced by $\lim_{t\to -\infty}$ and the indefinite integrals in Eqs.(2.47, 49)
are replaced by the definite integrals.
With these $R_i$ and $h^{(i)}_t$, define the 
\textit{restricted $m$-th order RG equation} for Eq.(3.1) to be
\begin{equation}
\dot{y} = \varepsilon R_1(y) + \varepsilon ^2R_2(y) + \cdots  + \varepsilon ^mR_m(y),\,\, y\in N_0,
\end{equation}
and define the \textit{restricted $m$-th order RG transformation} to be
\begin{equation}
\alpha _t(y) = \varphi _t(y) + \varepsilon h^{(1)}_t(y) + \cdots  + \varepsilon ^mh^{(m)}_t(y),\,\, y\in N_0,
\end{equation}
respectively. 

Note that the domains of them are restricted to $N_0$.
To see that they make sense, we prove the next lemma, which corresponds to Lemma 2.1.
\\[0.2cm]
\textbf{Lemma 3.2.}\, (i) The maps $R_i\,\, (i = 1,2, \cdots )$ are well-defined and $R_i(y) \in T_yN_0$ for any $y\in N_0$.
\\
(ii) The maps $h^{(i)}_t\,\, (i = 1,2, \cdots )$ are almost periodic functions with respect to $t$ uniformly in $y\in N_0$.
In particular, $h^{(i)}_t$ are bounded in $t\in \mathbf{R}$ if $y\in N_0$.
\\[0.2cm]
\textbf{Proof.}\, Let $\pi_s$ be the projection from $T_yM$ to the stable subspace $\bm{E}_s$ of $N_0$ 
and $\pi_{N_0}$ the projection from $T_yM$ to the tangent space $T_yN_0$.
Note that $\pi_s + \pi_{N_0} = id$.
By the definition of an attracting hyperbolic invariant manifold, there exist positive constants $C_1$ and $\alpha $
such that
\begin{equation}
|| \pi_s (D\varphi_t )_y^{-1}v || < C_1 e^{\alpha t} || v ||,
\end{equation}
for any $t <0, y\in N_0$ and $v\in T_yM$.
Since $\varphi _t(y) \in N_0$ for all $t \in \mathbf{R}$ and since $G_1$ is almost periodic in $t$,
there exists a positive constant $C_2$ such that $|| G_1(t, \varphi _t(y)) || < C_2$.
Then, $\pi_s R_1(y)$ proves to satisfy
\begin{eqnarray}
|| \pi_s R_1(y) || &\leq & \lim_{t\to -\infty} \frac{1}{t}\int^t \! || \pi_s (D\varphi _s)^{-1}_y G_1(s, \varphi _s(y))|| ds \nonumber \\
 & < & \lim_{t\to -\infty} \frac{1}{t}\int^t \! C_1C_2 e^{\alpha s} ds = 0.
\end{eqnarray}
This means that $R_1(y) \in T_yN_0$.

To prove (ii) of Lemma 3.2, note that the set
\begin{equation}
T(\delta ) = \{ \tau \, | \, || G_1(s+\tau , y) - G_1(s,y) || < \delta,\,
 || \varphi _{s+\tau}(y) - \varphi _s(y) ||< \delta  ,\quad
\forall s \in \mathbf{R},\, \forall y\in N_0\}
\end{equation}
is relatively dense. For $\tau \in T(\delta )$, $h^{(1)}_t (\varphi _\tau (y))$ is calculated as
\begin{eqnarray}
h^{(1)}_t (\varphi _\tau (y))
 &=& (D\varphi _t)_{\varphi _\tau (y)} \int^t_{-\infty} \! \left( 
(D\varphi _s)_{\varphi _\tau (y)}^{-1}G_1(s, \varphi _s\circ \varphi _\tau (y))
 - R_1(\varphi _\tau (y))\right) ds \nonumber \\
&=& \int^t_{-\infty} \! \left( (D\varphi _{s-t})_{y}^{-1}G_1(s, \varphi _{s+\tau}(y))
 - (D\varphi _{t})_{\varphi _{\tau} (y)}R_1(\varphi _\tau (y))\right) ds.
\end{eqnarray}
Putting $s' = s+ \tau$ yields
\begin{equation}
h^{(1)}_t (\varphi _\tau (y))= 
\int^{t+\tau}_{-\infty} \! \left( (D\varphi _{s'-(t+\tau)})_{y}^{-1}G_1(s'-\tau , \varphi _{s'}(y))
 - (D\varphi _{t})_{\varphi _{\tau} (y)}R_1(\varphi _\tau (y))\right) ds'.
\end{equation}
Since the space $T_yN_0$ is $(D\varphi _t)_y$-invariant, $\pi_s (D\varphi _t)_y R_1(y) = 0$.
This and Eq.(3.13) provide that
\begin{eqnarray}
|| \pi_sh^{(1)}_{t+\tau}(y) - \pi_s h^{(1)}_t (\varphi _\tau (y)) ||
& \leq & \int^{t+\tau}_{-\infty} \! || \pi_s  (D\varphi _{s-(t+\tau)})_{y}^{-1} || \cdot
|| G_1(s, \varphi _s (y)) - G_1(s-\tau , \varphi _{s}(y)) || ds \nonumber \\
&\leq & \int^{t+\tau}_{-\infty} \! \delta C_1e^{\alpha (s-(t + \tau))}ds = \delta C_1/\alpha .
\end{eqnarray}
Thus we obtain
\begin{eqnarray*}
||  \pi_sh^{(1)}_{t+\tau}(y) -  \pi_sh^{(1)}_{t}(y) ||
&\leq & ||  \pi_sh^{(1)}_{t+\tau}(y) -  \pi_sh^{(1)}_{t}( \varphi _\tau (y)) ||
 + || \pi_sh^{(1)}_{t}( \varphi _\tau (y))- \pi_sh^{(1)}_{t}(y)  || \\
& \leq & (C_1/\alpha  + L_t) \delta ,
\end{eqnarray*}
where $L_t$, which is bounded in $t$, is the Lipschitz constant of $\pi_s h^{(1)}_t|_{N_0}$.
This proves that $\pi_s h^{(1)}_t $ is an almost periodic function with respect to $t$.
On the other hand, $\pi_{N_0} h^{(1)}_t$ is written as
\begin{equation}
\pi_{N_0} h^{(1)}_{t}(y) =  \pi_{N_0}(D\varphi _t)_y \int^t_{-\infty} \! 
\left( \pi_{N_0}(D\varphi _{s})_y^{-1}G_1(s, \varphi _s(y)) - R_1(y)\right) ds.
\end{equation}
Since $\pi_{N_0}(D\varphi _t)_y$ is almost periodic, we can show that $\pi_{N_0} h^{(1)}_{t}(y)$ is almost
periodic in $t$ in the same way as the proof of Lemma 2.1 (ii).
This proves that $h^{(1)}_{t} = \pi_{s} h^{(1)}_{t} + \pi_{N_0} h^{(1)}_{t}$ is an almost periodic function in $t$.
The proof of Lemma 3.2 for $i=2,3,\cdots $ is done in a similar way by induction and we omit it here. \hfill $\blacksquare$
\\[0.2cm]
\indent Since $R_i(y) \in T_yN_0$, Eq.(3.7) defines a $\mathrm{dim} N_0$-dimensional differential equation on $N_0$.
\\[0.2cm]
\textbf{Remark 3.3.} Even if $h^{(i)}_t$ are defined by using indefinite integrals as Eqs.(2.47, 49),
we can show that $h^{(i)}_t$ is bounded as $t \to \infty$, though it is not bounded as $t\to -\infty$.
In this case, the theorems listed below are true for large $t$.

Now we are in a position to state main theorems of the restricted RG method, all proofs of which are the same as before
and omitted.
\\[0.2cm]
\textbf{Theorem 3.4.}\, Let $\alpha ^{(m)}_t$ be the restricted $m$-th order RG transformation for Eq.(3.1).
If $|\varepsilon |$ is sufficiently small, there exists a function $S(t,y,\varepsilon ),\, y\in N_0$ such that
\\[0.2cm]
(i)\, by changing the coordinates as $x = \alpha ^{(m)}_t(y)$, Eq.(3.1) is transformed into the system
\begin{equation}
\dot{y} = \varepsilon R_1(y) + \varepsilon ^2 R_2(y) + \cdots  + \varepsilon ^mR_m(y) + \varepsilon ^{m+1}S(t, y, \varepsilon ),
\end{equation}
(ii) \, $S$ is an almost periodic function with respect to $t$ uniformly in $y\in N_0$,
\\
(iii) \, $S(t,y,\varepsilon )$ is $C^1$ with respect to $t$ and $C^\infty$ with respect to $y\in N_0$ and $\varepsilon $.
\\[0.2cm]
\textbf{Theorem 3.5 (Error estimate).} \, Let $y(t)$ be a solution of the restricted $m$-th order RG equation and $\alpha ^{(m)}_t$
the restricted $m$-th order RG transformation. There exist positive constants $\varepsilon _0, C$ and $T$ such that a solution $x(t)$ of
Eq.(3.1) with $x(0) = \alpha ^{(m)}_0(y(0)) \in \alpha _0^{(m)} (N_0)$ satisfies the inequality
\begin{equation}
|| x(t) - \alpha ^{(m)}_t (y(t)) || < C |\varepsilon |^m,
\end{equation}
for $|\varepsilon | < \varepsilon _0$ and $0 \leq t \leq T/|\varepsilon |$.
\\[0.2cm]
\textbf{Theorem 3.6 (Existence of invariant manifolds).} \, Suppose that $R_1(y) = \cdots = R_{k-1}(y) = 0$ and
$\varepsilon ^k R_k(y)$ is the first non-zero term in the restricted RG equation for Eq.(3.1).
If the vector field $R_k(y)$ has a boundaryless compact normally hyperbolic invariant manifold $L \subset N_0$,
then for sufficiently small $\varepsilon >0$,
the system (3.2) has an invariant manifold $L_\varepsilon $ on the $(s,x)$ space which is diffeomorphic to $\mathbf{R} \times L$.
In particular, the stability of $L_\varepsilon $ coincides with that of $L$.
\\[0.2cm]
\textbf{Theorem 3.7 (Inheritance of symmetries).}\, 
Suppose that an $\varepsilon $-independent Lie group $H$ acts on $N_0$.
If the vector field $f$ and $g$ are invariant under the action of $H$,
then the restricted $m$-th order RG equation for Eq.(3.1) is also invariant under the action of $H$.
\\[0.2cm]
\indent Recall that our purpose is to construct the invariant manifold $N_\varepsilon $ of Eq.(3.1) and the flow
on $N_\varepsilon $ approximately. The flow on $N_\varepsilon $ is well understood by Theorems 3.4 to 3.6,
and $N_\varepsilon $ is given by the next theorem.
\\[0.2cm]
\textbf{Theorem 3.8.}\,  Let $\alpha ^{(m)}_t$ be the restricted $m$-th order RG transformation for Eq.(3.1).
Then, the set $\{ (t, x) \, | \, x\in \alpha ^{(m)}_t (N_0)\}$ lies within an $O(\varepsilon ^{m+1})$ neighborhood of
the attracting invariant manifold $N_\varepsilon $ of Eq.(3.2).
\\[0.2cm]
\textbf{Proof.} \, Though the maps $R_i(y)$ and $\alpha ^{(m)}_t$ are
defined on $N_0$, we can extend them to the maps defined on $V \supset N_0$ so that Eq.(3.7) is $C^1$
close to Eq.(3.16) on $V$ and that $N_0$ is an attracting normally hyperbolic invariant manifold of Eq.(3.7).
Then the same argument as the proof of Thm.2.9 proves Theorem 3.8. \hfill $\blacksquare$
\\

If the vector field $g$ is independent of $t$ and Eq.(3.1) is autonomous, we can prove the next theorems.
\\[0.2cm]
\textbf{Theorem 3.9 (Existence of invariant manifolds).}\, \, Suppose that $R_1(y) = \cdots = R_{k-1}(y) = 0$ and
$\varepsilon ^k R_k(y)$ is the first non-zero term in the restricted RG equation for Eq.(3.1) with $t$-independent $g$.
If the vector field $R_k(y)$ has a boundaryless compact normally hyperbolic invariant manifold $L$,
then for sufficiently small $\varepsilon >0$,
Eq.(3.1) has an invariant manifold $L_\varepsilon $, which is diffeomorphic to $L$.
In particular, the stability of $L_\varepsilon $ coincides with that of $L$.
\\[0.2cm]
\textbf{Theorem 3.10 (Additional symmetry).}\, 
The restricted RG equation for Eq.(3.1) with $t$-independent $g$ is invariant under the action of the 
one-parameter group $\{ \varphi _t : N_0\to N_0 \, | \, t\in \mathbf{R}\}$.
In other words, $R_i$ satisfies the equality
\begin{equation}
R_i(\varphi _t(y)) = (D\varphi _t)_y R_i(y),\,\, y\in N_0,
\end{equation}
for $i=1,2, \cdots $.
\\[-0.2cm]

For autonomous systems, Thm.3.8 is restated as follows:
Recall that if the function $g$ depends on $t$, the attracting invariant manifold $N_\varepsilon $ of Eq.(3.1)
and the approximate invariant manifold described in Thm.3.8 depend on $t$ in the sense that
they lie on the $(s,x)$ space.
If Eq.(3.1) is autonomous, its attracting invariant manifold $N_\varepsilon $ lies on $M$
and is independent of $t$. Thus we want to construct an approximate invariant manifold of $N_\varepsilon $
so that it is also independent of $t$.
\\[0.2cm]
\textbf{Theorem 3.11.}\,  Let $\alpha ^{(m)}_t$ be the restricted $m$-th order RG transformation for Eq.(3.1).
If $g$ is independent of $t$, 
the set $\alpha ^{(m)}_t (N_0) = \{ \alpha ^{(m)}_t(y) \, | \, y \in N_0\}$ is independent of $t$
and lies within an $O(\varepsilon ^{m+1})$ neighborhood of
the attracting invariant manifold $N_\varepsilon $ of Eq.(3.1).
\\[0.2cm]
\textbf{Proof.} We have to show that the set $\alpha ^{(m)}_t (N_0)$ is independent of $t$.
Indeed, we know the equality $h^{(i)}_t(\varphi _{t'}(y)) = h^{(i)}_{t+t'}(y)$ as is shown in the proof of
the Thm.2.15. This proves that
\begin{equation}
\alpha ^{(m)}_{t+t'} (N_0) = \alpha ^{(m)}_t (\varphi _{t'}(N_0)) = \alpha ^{(m)}_t (N_0).
\end{equation}
The rest of the proof is the same as the proofs of Thm.2.14 and Thm.3.8. \hfill $\blacksquare$

%%%%%%%%%%%%%%%%%%%%%%%%%%%%%%%%%%%%%%%%%%%%%%%%%%%%%%%%%%%%%%%%%%%%%%%%%%%%%%%%%%%%%%%%%%%%%%%%%%%%%%%%%%%%%%

\subsection{Center manifold reduction}

The restricted RG method recovers the approximation theory of center manifolds (Carr [7]).
Consider a system of the form
\begin{eqnarray}
\dot{x} &=& Fx + \varepsilon g(x, \varepsilon ) \nonumber \\
&=& Fx + \varepsilon g_1(x) + \varepsilon ^2 g_2(x) + \cdots ,\,\, x\in \mathbf{R}^n,
\end{eqnarray}
where unperturbed term $Fx$ is linear.
For this system, we suppose that
\\[0.2cm]
\textbf{(E1)} all eigenvalues of the $n \times n$ constant matrix $F$ are on the imaginary axis or the left half plane.
The Jordan block corresponding to eigenvalues on the imaginary axis is diagonalizable.
\\
\textbf{(E2)} $g$ is $C^\infty$ with respect to $x$ and $\varepsilon $ such that $g(0, \varepsilon ) = 0$.
\\[0.2cm]
If all eigenvalues of $F$ are on the left half plane, the origin is a stable fixed point and the flow near the origin is trivial.
In what follows, we suppose that at least one eigenvalue is on the imaginary axis.
In this case, Eq.(3.20) has a center manifold which is tangent to the center subspace $N_0$ at the origin.
The center subspace $N_0$, which is spanned by eigenvectors associated with eigenvalues on the imaginary axis,
is an attracting normally hyperbolic invariant manifold of the unperturbed term $Fx$, and 
the flow of $Fx$ on $N_0$ is almost periodic. However, since $N_0$ is not compact,
we take an $n$-dimensional closed ball $K$ including the origin and consider $N_0 \cap K$.
Then, we obtain the next theorem as a corollary of Thm.3.9 and Thm.3.11.
\\[0.2cm]
\textbf{Theorem 3.12 (Approximation of Center Manifolds, [12]).}\, 
Let $\alpha ^{(m)}_t$ be the restricted $m$-th order RG transformation for Eq.(3.20) 
and $K$ a small compact neighborhood of the origin. 
Then, the set $\alpha _t^{(m)} (K \cap N_0)$ lies within an $O(\varepsilon ^{m+1})$ neighborhood of the center manifold of Eq.(3.20).
The flow of Eq.(3.20) on the center manifold is well approximated by those of the restricted RG equation.
In particular,
suppose that $R_1(y) = \cdots = R_{k-1}(y) = 0$ and
$\varepsilon ^k R_k(y)$ is the first non-zero term in the restricted RG equation.
If the vector field $R_k(y)$ has a boundaryless compact normally hyperbolic invariant manifold $L$,
then for sufficiently small $\varepsilon >0$, 
Eq.(3.20) has an invariant manifold $L_\varepsilon $ on the center manifold, which is diffeomorphic to $L$.
The stability of $L_\varepsilon $ coincides with that of $L$.
\\[-0.2cm]

See Chiba [12] for the detail of the proof and examples.

%%%%%%%%%%%%%%%%%%%%%%%%%%%%%%%%%%%%%%%%%%%%%%%%%%%%%%%%%%%%%%%%%%%%%%%%%%%%%%%%%%%%%%%%%%%%%%%%%%%%%%%%

\subsection{Geometric singular perturbation method}

The restricted RG method can also recover the geometric singular perturbation method proposed by Fenichel [19].

Consider the autonomous system 
\begin{equation}
\dot{x} = f(x) + \varepsilon g(x, \varepsilon ),\,\, x\in \mathbf{R}^n
\end{equation}
on $\mathbf{R}^n$ with the assumption that
\\[0.2cm]
\textbf{(F)}\, suppose that $f$ and $g$ are $C^\infty$ with respect to $x$ and $\varepsilon $,
and that $f$ has an $m$-dimensional attracting normally hyperbolic invariant manifold $N_0$ which consists of fixed points of $f$,
where $m < n$.
\\[-0.2cm]

Note that this system satisfies the assumptions (D1) and (D2), so that Thm.3.9 to Thm.3.11 hold.
The invariant manifold $N_0$ consisting of fixed points of the unperturbed system 
is called the \textit{critical manifold} (see [1]).
For this system, Fenichel[19] proved that there exist local coordinates $(u, v)$ such that the system (3.21)
is expressed as
\begin{equation}
\left\{ \begin{array}{ll}
\dot{u} = \varepsilon \hat{g}_1(u,v, \varepsilon ), & u\in \mathbf{R}^m,  \\
\dot{v} = \hat{f}(u,v) + \varepsilon \hat{g}_2(u,v, \varepsilon ), & v\in \mathbf{R}^{n-m}, \\
\end{array} \right.
\end{equation}
where $\hat{f}(u,0) = 0$ for any $u\in \mathbf{R}^m$.
In this coordinate, the critical manifold $N_0$ is locally given as the $u$-plane.
Further he proved the next theorem.
\\[0.5cm]
\textbf{Theorem 3.13 (Fenichel [19]).} 
Suppose that the system $\dot{u} = \varepsilon \hat{g}_1(u,0,0)$ has a compact normally hyperbolic invariant manifold $L$.
If $\varepsilon >0$ is sufficiently small, the system (3.22) has an invariant manifold $L_\varepsilon $ which is diffeomorphic to $L$.
\\

This method to obtain an invariant manifold of (3.22) is called the \textit{geometric singular perturbation method}.
By using the fact that $\varphi _t(u) = u$ for $u \in N_0$, it is easy to verify that
the system $\dot{u} = \varepsilon  \hat{g}_1(u,0,0)$ described above is just the restricted first order RG equation for Eq.(3.22).
Thus Thm.3.13 immediately follows from Thm.3.9.
Note that in our method, we need not change the coordinates so that Eq.(3.21) is transformed into the form of Eq.(3.22).
\\[0.5cm]
\textbf{Example 3.14.}@Consider the system on $\mathbf{R}^2$
\begin{equation}
\left\{ \begin{array}{l}
\dot{x}_1 = -x_1 + (x_1 + c)x_2,  \\
\varepsilon \dot{x}_2 = x_1 - (x_1 + 1)x_2,  \\
\end{array} \right.
\end{equation}
where $0< c <1$ is a constant.
This system arises from a model of the kinetics of enzyme reactions (see Carr [7]).
Set $t = \varepsilon s$ and denote differentiation with respect to $s$ by $\,'\,$.
Then, the above system is rewritten as
\begin{equation}
\left\{ \begin{array}{l}
x_1' = \varepsilon (-x_1 + cx_2 + x_1x_2),  \\
x_2' = x_1 - x_2 - x_1x_2.  \\
\end{array} \right.
\end{equation}
The attracting critical manifold $N_0$ of this system is expressed as the graph of the function
\begin{equation}
x_2 = h(x_1) := \frac{x_1}{1 + x_1}.
\end{equation}
Since the restricted first order RG transformation for Eq.(3.24) is given by
\begin{equation}
\alpha ^{(1)}_t(y_1) = \left(
\begin{array}{@{\,}c@{\,}}
y_1 \\
y_1/(1 + y_1)
\end{array}
\right) + \varepsilon  \left(
\begin{array}{@{\,}c@{\,}}
0 \\
-(c-1)y_1/ (1+y_1)^4
\end{array}
\right),
\end{equation}
Theorem 3.11 proves that the attracting invariant manifold of Eq.(3.24) is given as the graph of
\begin{equation}
x_2 = \frac{x_1}{1+x_1} - \varepsilon \frac{(c-1)x_1}{(1+x_1)^4} + O(\varepsilon ^2).
\end{equation}
If $|x_1|$ is sufficiently small, it is expanded as 
\begin{eqnarray}
x_2 &=& x_1(1 - x_1) - \varepsilon (c-1)x_1(1 - 4x_1) + O(x_1^3, \varepsilon ^2) \nonumber  \\
&=& (1-\varepsilon (c-1))x_1 - (1-4\varepsilon (c-1))x_1^2 + O(x_1^3, \varepsilon ^2).
\end{eqnarray}
This result coincides with the result obtained by the local center manifold theory (see Carr [7]).
The restricted first order RG equation on $N_0$ is given by
\begin{eqnarray}
y_1' = \varepsilon \frac{(c-1)y_1}{1+y_1}.
\end{eqnarray}
This RG equation describes a motion on the invariant manifold (3.27) approximately.
Since it has the stable fixed point $y_1 =0$ if $c<1$, 
the system (3.24) also has a stable fixed point $(x_1,x_2) = (0,0)$ by virtue of Theorem 3.9.

%%%%%%%%%%%%%%%%%%%%%%%%%%%%%%%%%%%%%%%%%%%%%%%%%%%%%%%%%%%%%%%%%%%%%%%%%%%%%%%%%%%%%%%%%%%%%%%%%%%%

\subsection{Phase reduction}

Consider a system of the form
\begin{equation}
\dot{x} = f(t,x) + \varepsilon g_1(t,x),\, x\in \mathbf{R}^n.
\end{equation}
For this system, we suppose that
\\[0.2cm]
\textbf{(G1)} \, the vector fields $f$ and $g$ are $C^\infty$ in $x$, $C^1$ in $t$,
and $T$-periodic in $t$. It need not be the least period.
In particular, $f$ and $g$ are allowed to be independent of $t$.
\\
\textbf{(G2)} \, the unperturbed system $\dot{x} = f(t, x)$ has a $k$-parameter family of $T$-periodic solutions
which constructs an attracting invariant torus $\mathbb{T}^k \subset \mathbf{R}^n$. 

Let $\alpha  = (\alpha _1 , \cdots , \alpha _k)$ be coordinates on $\mathbb{T}^k$, so-called \textit{phase variables}.
It is called the \textit{phase reduction} to derive equations on $\alpha $ which govern the dynamics of Eq.(3.30)
on the invariant torus. The phase reduction was first introduced by Malkin [33,34] and rediscovered by many authors.
In this subsection, we show that the RG method can recover the phase reduction method.

The next theorem is due to Malkin [33,34].
See also Hahn [23], Blekhman [4], and Hoppensteadt and Izhikevich [25].
\\
\textbf{Theorem 3.15 (Malkin [33,34]).} \, Consider the system (3.30) with the assumptions (G1) and (G2).
Let $\alpha = (\alpha _1, \cdots , \alpha _k)$ be phase variables and $U(t; \alpha )$ the periodic solutions
of the unperturbed system parameterized by $\alpha $.
Suppose that the adjoint equation
\begin{equation}
\frac{dQ_i}{dt}= - \left( \frac{\partial f}{\partial x} (t, U(t; \alpha ))\right) ^{\mathsf{T}}Q_i
\end{equation}
has exactly $k$ independent $T$-periodic solutions $Q_1(t; \alpha ), \cdots , Q_k(t; \alpha )$,
where $A^{\mathsf{T}}$ denotes the transpose matrix of a matrix $A$.
Let $Q = Q(t; \alpha )$ be the $k \times n$ matrix whose columns are these solutions such that
\begin{equation}
Q^{\mathsf{T}} \frac{\partial U}{\partial \alpha }(t; \alpha ) = id.
\end{equation}
Then, Eq.(3.30) has a solution of the form
\begin{equation}
x(t) = U(t, \alpha (t)) + O(\varepsilon ),
\end{equation}
where $\alpha (t)$ is a solution of the system
\begin{equation}
\frac{d \alpha }{dt} = \frac{\varepsilon }{T} \int^T_{0} \! Q(s; \alpha ) ^{\mathsf{T}} g_1(s, U(s; \alpha ))ds.
\end{equation}

Now we show that the system (3.34) of the phase variables is just the first order RG equation.
Note that the system (3.30) satisfies the assumptions (D1) and (D2) with $N_0 = \mathbb{T}^k$
and the RG method is applicable.
The restricted first order RG equation for Eq.(3.30) is given by
\begin{equation}
\dot{y} = \varepsilon \lim_{t\to -\infty} \int^t_{0} \!
\left( \frac{\partial \varphi _s}{\partial y} (y)\right) ^{-1} g_1(s, \varphi _s(y)) ds, \quad y\in \mathbb{T}^k.
\end{equation}
Let us change the coordinates by using the $k$-parameter family of periodic solutions as $y = U(0;\alpha )$.
Then, Eq.(3.35) is rewritten as
\begin{equation}
\frac{\partial U}{\partial \alpha } (0; \alpha ) \dot{\alpha } =
\varepsilon \lim_{t\to -\infty} \int^t_{0} \!
\left( \frac{\partial \varphi _s}{\partial y} (U(0; \alpha ))\right) ^{-1} g_1(s, U(s; \alpha )) ds.
\end{equation}
Since $U(t; \alpha ) = \varphi _t(U(0; \alpha ))$, the equality
\begin{equation}
\frac{\partial U}{\partial \alpha }(t; \alpha ) = \frac{\partial \varphi _t}{\partial y}(U(0 ; \alpha )) 
\frac{\partial U}{\partial \alpha }(0; \alpha )
\end{equation}
holds. Then, Eqs.(3.32), (3.36) and (3.37) are put together to obtain
\begin{eqnarray}
\dot{\alpha } &=& \varepsilon \lim_{t\to -\infty} \int^t_{0} \! \left( \frac{\partial U}{\partial \alpha }(0; \alpha ) \right) ^{-1}
\left( \frac{\partial \varphi _s}{\partial y} (U(0; \alpha ))\right) ^{-1} g_1(s, U(s; \alpha )) ds \nonumber \\
&=& \varepsilon \lim_{t\to -\infty} \int^t_{0} \! 
\left( \frac{\partial U}{\partial y} (s; \alpha )\right) ^{-1} g_1(s, U(s; \alpha )) ds \nonumber \\
&=& \varepsilon \lim_{t\to -\infty} \int^t_{0} \! 
Q(s; \alpha ) ^{\mathsf{T}} g_1(s, U(s; \alpha )) ds.
\end{eqnarray} 
Since the integrand in the above equation is $T$-periodic, Eq.(3.38) is reduced to Eq.(3.34).

%%%%%%%%%%%%%%%%%%%%%%%%%%%%%%%%%%%%%%%%%%%%%%%%%%%%%%%%%%%%%%%%%%%%%%%%%%%%%%%%%%%%%%%%%%%%%%%%%%%%%%%%%%%%%
%%%%%%%%%%%%%%%%%%%%%%%%%%%%%%%%%%%%%%%%%%%%%%%%%%%%%%%%%%%%%%%%%%%%%%%%%%%%%%%%%%%%%%%%%%%%%%%%%%%%%%%%%%%%%

\section{Relation to other singular perturbation methods}

In the previous section, we have seen that the restricted RG method unifies the center manifold reduction,
the geometric singular perturbation method and the phase reduction.
In this section, we show that the RG method described in Section 2 unifies the traditional singular perturbation
methods, such as the averaging method, the multiple time scale method and the normal forms theory.
We will also give explicit formulas for regular perturbation solutions and show how the RG equation is derived from
the regular perturbation method through the envelope theory.

%%%%%%%%%%%%%%%%%%%%%%%%%%%%%%%%%%%%%%%%%%%%%%%%%%%%%%%%%%%%%%%%%%%%%%%%%%%%%%%%%%%%%%%%%%%%%%%%%%%%%%%%%%

\subsection{Averaging method}

The averaging method is one of the most traditional and famous singular perturbation methods based on an idea of 
averaging a given equation, which is periodic in time $t$, with respect to $t$ to obtain an autonomous equation
(see Sanders, Verhulst and Murdock [47]).
In most literature, only first and second order averaging equations are given and they coincide with first 
and second order RG equations, respectively.
Our $m$-th order RG equation gives a generalization of such lower order averaging equations.
In Chiba and Paz\'{o}[13], the third and fifth order RG equations are used to unfold a phase diagram
of the Kuramoto model of coupled oscillators.

Many authors formulated the averaging method for time-periodic differential equations.
However, the RG equation (the averaging equation) can be defined as long as limits in Eqs.(2.11, 13) exist
and $R_i(y)$ are well-defined even if a system (2.1) does not satisfy the assumption (A).
If $R_1 , \cdots , R_m$ are well-defined for a system (2.1), which need not satisfy the assumption (A),
we say (2.1) satisfies the \textit{KBM condition up to order $m$} 
(KBM stands for Krylov, Bogoliubov and Mitropolsky, see Bogoliubov and Mitropolsky [6]).
We showed that a system (2.1) with the assumption (A) satisfies the KBM condition up to all order (Lemma 2.1).
If the assumption (A) is violated,
the error estimate (2.27) for approximate solutions may get worse, 
such as $|| x(t) - \alpha ^{(m)}_t(y(t)) || < C\sqrt{\varepsilon }$,
even if the KBM condition is satisfied and thus the RG equation is defined.
See Sanders \textit{et al.} [47] and DeVille \textit{et al.} [32] for such examples.

Even if (2.1) satisfies the KBM condition, we need additional assumptions to prove Thm.2.9.
It is because to apply Fenichel's theorem, which was used in the proof of Thm.2.9,
we have to show that the error function $S(t,y, \varepsilon )$ in Eq.(2.23) is $C^1$ with respect to $t, x$
and bounded in $t \in \mathbf{R}$. See Chiba [10] for more detail.

%%%%%%%%%%%%%%%%%%%%%%%%%%%%%%%%%%%%%%%%%%%%%%%%%%%%%%%%%%%%%%%%%%%%%%%%%%%%%%%%%%%%%%%%%%%%%%%%%%%%%%%%%%

\subsection{Multiple time scale method}

The multiple time scale method [3,38,44] was perhaps introduced in the early 20th century and now it is one of the most
used perturbation techniques along with the averaging method.
In this subsection, we give a brief review of the method and show that it yields the same results as the RG method.

Consider the system (2.1) satisfying the assumption (A).
Let us introduce the new time scales $t_m$ as
\begin{equation}
t_m = \varepsilon ^m t,\,\, m= 0,1,2,\cdots ,
\end{equation}
and consider $t$ and $t_0, t_1, t_2 ,\cdots $ to be independent of each other.
Then, the ``total" derivative $d/dt$ is rewritten as
\begin{equation}
\frac{d}{dt} = \frac{\partial }{\partial t_0} + \varepsilon \frac{\partial }{\partial t_1}
 + \varepsilon ^2 \frac{\partial }{\partial t_2} + \cdots .
\end{equation}
Let us expand the dependent variable $x$ as
\begin{equation}
x = x_0 + \varepsilon x_1 + \varepsilon ^2 x_2 + \cdots ,
\end{equation}
where $x_i = x_i(t_0, t_1, t_2, \cdots )$. Substituting Eqs.(4.2) and (4.3) into (2.1), we obtain
\begin{equation}
\left( \frac{\partial }{\partial t_0} + \varepsilon \frac{\partial }{\partial t_1}
 + \varepsilon ^2 \frac{\partial }{\partial t_2} + \cdots \right)
\left(x_0 + \varepsilon x_1 + \varepsilon ^2 x_2 + \cdots \right)
 = \sum^\infty_{k=1} \varepsilon ^k g_k(t_0, x_0 + \varepsilon x_1 + \varepsilon ^2 x_2 + \cdots ).
\end{equation}
Expanding the both sides of the above in $\varepsilon $ and equating the coefficients of each $\varepsilon ^k$,
we obtain ODEs on $x_0 , x_1, x_2,\cdots $;
\begin{equation}
\left\{ \begin{array}{l}
\displaystyle \frac{\partial x_0}{\partial t_0} = 0,  \\[0.2cm]
\displaystyle \frac{\partial x_1}{\partial t_0} = G_1(t_0, x_0) - \frac{\partial x_0}{\partial t_1}, \\[0.2cm]
\displaystyle \frac{\partial x_2}{\partial t_0} = G_2(t_0, x_0, x_1)
                - \frac{\partial x_1}{\partial t_1} - \frac{\partial x_0}{\partial t_2}, \\[0.2cm]
\quad \quad \vdots \\
\displaystyle \frac{\partial x_m}{\partial t_0} = G_m(t_0, x_0, \cdots , x_{m-1})
                - \sum^m_{j=1}\frac{\partial x_{m-j}}{\partial t_j},
\end{array} \right.
\end{equation}
where the functions $G_k,\, k= 1,2, \cdots $ are defined through Eq.(2.6).
Let $x_0 = y = y(t_1, t_2, \cdots )$ be a solution of the zeroth order equation.
Then, a general solution of the first order equation of (4.5) is given as
\begin{equation}
x_1 = B_1 + \int^{t_0}\! \left( G_1(s, y) - \frac{\partial y}{\partial t_1} \right) ds
 = B_1 + u^{(1)}_{t_0}(y) + R_1(y) t_0 - \frac{\partial y}{\partial t_1}t_0, 
\end{equation}
where $B_1 = B_1(y; t_1, t_2,\cdots )$ is independent of $t_0$ and where
$R_1$ and $u^{(1)}_{t_0}$ are defined by Eqs.(2.11) and (2.12), respectively.
Now we define  $\partial y/ \partial t_1$ so that $x_1$ above is bounded in $t_0$ :
\begin{equation}
R_1(y) t_0 - \frac{\partial y}{\partial t_1} t_0 = 0.
\end{equation}
This condition is called the \textit{non-secularity condition} and it yields
\begin{equation}
\frac{\partial y}{\partial t_1} = R_1(y), \quad x_1 = x_1(y) = B_1 + u^{(1)}_{t_0}(y).
\end{equation}
Next thing to do is calculating $x_2$.
The equation on $x_2$ is written as
\begin{eqnarray}
\frac{\partial x_2}{\partial t_0} &=& G_2(t_0, y, B_1 + u^{(1)}_{t_0}(y))
 - \frac{\partial }{\partial t_1} (B_1 + u^{(1)}_{t_0}(y)) - \frac{\partial y}{\partial t_2} \nonumber \\
&=& \frac{\partial g_1}{\partial y}(t_0, y)(B_1 + u^{(1)}_{t_0}(y)) + g_2(t_0, y)
- \frac{\partial u^{(1)}_{t_0}}{\partial y}(y)R_1(y) - \frac{\partial B_1}{\partial t_1}
 - \frac{\partial B_1}{\partial y} R_1(y) - \frac{\partial y}{\partial t_2},
\end{eqnarray}
a general solution of which is given by
\begin{eqnarray}
x_2 &=& B_2 + \int^{t_0} \! \left( \frac{\partial g_1}{\partial y}(s, y)(B_1 + u^{(1)}_{s}(y)) + g_2(s, y)
- \frac{\partial u^{(1)}_{s}}{\partial y}(y)R_1(y) - \frac{\partial B_1}{\partial t_1}
 - \frac{\partial B_1}{\partial y} R_1(y) - \frac{\partial y}{\partial t_2}\right) ds \nonumber \\
&=& B_2 + u^{(2)}_{t_0}(y) + \frac{\partial u^{(1)}_{t_0}}{\partial y}(y)B_1 \nonumber \\
& & \quad + \frac{\partial R_1}{\partial y}(y)B_1 t_0 + R_2(y) t_0 - \frac{\partial B_1}{\partial t_1}t_0
 - \frac{\partial B_1}{\partial y}R_1(y)t_0 - \frac{\partial y}{\partial t_2} t_0,
\end{eqnarray}
where $B_2 = B_2(y; t_1, t_2,\cdots )$ is independent of $t_0$.
If we impose the non-secularity condition so that $x_2$ is bounded in $t_0$,
we obtain
\begin{eqnarray}
& & \frac{\partial y}{\partial t_2} = R_2(y) - \frac{\partial B_1}{\partial t_1} - [B_1, R_1](y), \\
& & x_2 = x_2(y) = B_2 + u^{(2)}_{t_0}(y) + \frac{\partial u^{(1)}_{t_0}}{\partial y}(y)B_1,
\end{eqnarray}
where the commutator $[\bm{\cdot},\bm{\cdot} ]$ is defined in Eq.(\ref{2_4_8}).

By proceeding in a similar manner, we can show that the non-resonance condition at each stage
yields equations on $y$ of the form
\begin{equation}
\frac{\partial y}{\partial t_i} = R_i(y) + Q_i(R_1, \cdots ,R_{i-1}, B_1 , \cdots ,B_{i-2})(y)
 - \frac{\partial B_{i-1}}{\partial t_1} - [B_{i-1}, R_1](y),
\end{equation}
for $i = 2,3, \cdots $, where $Q_i$ is a function of $R_1, \cdots ,R_{i-1}$ and $B_1 , \cdots ,B_{i-2}$.
Combining these equations on $y$, we obtain
\begin{eqnarray}
\frac{dy}{dt} &=& \frac{\partial y}{\partial t_0} + \varepsilon \frac{\partial y}{\partial t_1}
 + \varepsilon ^2 \frac{\partial y}{\partial t_2} + \cdots  \nonumber \\
&=& \varepsilon R_1(y) + \varepsilon ^2 \left( R_2(y) - \frac{\partial B_1}{\partial t_1} - [B_1, R_1](y) \right) + \cdots.
\end{eqnarray}
The functions $B_1, B_2, \cdots $ are still left undetermined.
We propose two different ways to determine them.
\\[0.2cm]
\textbf{(i)} \, If we define $B_i$'s as solutions of the linear inhomogeneous differential equations
\begin{equation}
\frac{\partial B_{i}}{\partial t_1}  =  -[B_{i}, R_1](y) + R_{i+1}(y) + Q_{i+1}(R_1, \cdots ,R_{i}, B_1 , \cdots ,B_{i-1})(y),
\end{equation}
then Eq.(4.14) is reduced to the equation
\begin{equation}
\frac{dy}{dt} = \varepsilon R_1(y).
\end{equation}
Let $y = y(t)$ be a solution of this equation.
Then, we can prove that a curve defined by
\begin{eqnarray}
\alpha ^{(m)}_t (y(t)) 
&:=& y(t) + \varepsilon  x_1(y(t)) + \varepsilon ^2x_2(y(t)) + \cdots + \varepsilon ^m x_m(y(t)) \nonumber \\
&=& y + \varepsilon \left( B_1 + u^{(1)}_t(y) \right)
 + \varepsilon ^2 \left( B_2 + u^{(2)}_{t}(y) + \frac{\partial u^{(1)}_{t}}{\partial y}(y)B_1 \right) + \cdots \Bigl|_{y = y(t)}
\end{eqnarray}
provides an approximate solution for Eq.(2.1) and satisfies the same inequality as (2.27) (see Murdock[38]).
\\
\textbf{(ii)} Otherwise, we suppose that $B_i$'s are independent of $t_1, t_2, \cdots $ so that
Eq.(4.14) is independent of $t_1, t_2, \cdots $.
Then Eq.(4.14) takes the form
\begin{equation}
\frac{dy}{dt} = \varepsilon R_1(y) + \varepsilon ^2 \left( R_2(y) - [B_1, R_1](y) \right) + \cdots, 
\end{equation}
which coincides with the (simplified) RG equation (see Sec.2.4).

%%%%%%%%%%%%%%%%%%%%%%%%%%%%%%%%%%%%%%%%%%%%%%%%%%%%%%%%%%%%%%%%%%%%%%%%%%%%%%%%%%%%%%%%%%%%%%%%%%%%%%%%%%

\subsection{Normal forms}

Let us consider the system (2.5) with assumptions (C1) and (C2).
The technique of normal forms is used to analyze local dynamics of such a system
and had been well developed in the last century (Murdock [39]).
Let us recall the definition of normal forms.
For a given system (2.5), there exist a \textit{time-independent} local coordinate transformation
$x = h(z)$ defined near the origin which brings Eq.(2.5) into the form
\begin{equation}
\dot{z} = Fz + \varepsilon \tilde{g}_1(z) + \cdots  + \varepsilon^m \tilde{g}_m(z) + \varepsilon ^{m+1}\tilde{S}(z,\varepsilon ) 
\end{equation}
with the properties that $\tilde{g}_i(z),\, i = 1, \cdots  ,m$ satisfy $\tilde{g}_i(e^{Ft}z) = e^{Ft}\tilde{g}_i(z)$
and that $\tilde{S}$ is $C^\infty$ in $z$ and $\varepsilon $.
Then the truncated system
\begin{equation}
\dot{z} = Fz + \varepsilon \tilde{g}_1(z) + \cdots  + \varepsilon^m \tilde{g}_m(z) 
\end{equation}
is called the \textit{normal form of Eq.(2.5) up to order $m$}.

Note that if the matrix $F$ is \textit{not} diagonalizable, different definitions are adopted for normal forms.

Now we consider the RG equation and the RG transformation for Eq.(2.5), respectively, defined by Eqs.(2.50) and (2.51)
with Eqs.(2.46) to (2.49).
The flow $\varphi _t(y)$ in Eqs.(2.46) to (2.49) is given as $\varphi _t(y) = e^{Ft}y$ by using the fundamental matrix $e^{Ft}$
in our situation. Thm.2.5 shows that the RG transformation $x = \alpha ^{(m)}_t(y)$ brings Eq.(2.5) into the form of 
Eq.(2.23). Because of Thm.2.15, further change of variables as $y = e^{-Ft}z$ transforms Eq.(2.23)
into the system
\begin{equation}
\dot{z} = Fz + \varepsilon R_1(z) + \cdots  + \varepsilon ^m R_m(z) + \varepsilon ^{m+1}e^{Ft}S(t, e^{-Ft}z, \varepsilon )
\end{equation}
with the property that $R_i(e^{Ft}z) = e^{Ft}R_i(z)$ for $i= 1, \cdots  ,m$.

Now we need the next lemma.
\\[0.2cm]
\textbf{Lemma 4.1.}\, For the autonomous system (2.5), the RG transformation satisfies the equality
$\alpha ^{(m)}_t(e^{Ft'}y) = \alpha ^{(m)}_{t+t'}(y)$ for all $t, t', y$.
\\[-0.2cm]

This lemma immediately follows from the equality $h^{(i)}_t(\varphi _{t'}(y)) = h^{(i)}_{t+t'}(y),\, i= 1,2, \cdots $
proved in the proof of Thm.2.15.
Thus it turns out that the composition of two transformations 
$x = \alpha ^{(m)}_t(y)$ and $y = e^{-Ft} z$ is independent of $t$:
\begin{equation}
x = \alpha ^{(m)}_t(e^{-Ft} z) = \alpha ^{(m)}_0 (z).
\end{equation}
This proves that Eq.(4.21) is obtained from Eq.(2.5) by the time-independent transformation,
and thus the truncated system
\begin{equation}
\dot{z} = Fz + \varepsilon R_1(z) + \cdots  + \varepsilon ^m R_m(z) 
\end{equation}
is a normal form of Eq.(2.5) up to order $m$.

Note that normal forms of a given system are not unique in general and the simplest form among them
is called a \textit{hyper-normal form} ([2,39,40]).
The RG method can also provide hyper-normal forms by choosing undetermined integral constants in Eqs.(2.47)
and (2.49) appropriately as was discussed in Sec.2.4. See also Chiba [11].

%%%%%%%%%%%%%%%%%%%%%%%%%%%%%%%%%%%%%%%%%%%%%%%%%%%%%%%%%%%%%%%%%%%%%%%%%%%%%%%%%%%%%%%%%%%%%%%%%%%%%%%%%%

\subsection{Regular perturbation and envelopes}

In this subsection, we give some properties of the regular perturbation method.
The RG equation is obtained from regular perturbation solutions through Kunihiro's idea based on envelopes [28,29].

Consider the system (2.1) on $\mathbf{R}^n$ with the assumption (A).
Let us construct a formal solution of it of the form
\begin{equation}
x = \hat{x} = x_0 + \varepsilon x_1 + \varepsilon ^2 x_2 + \cdots .
\end{equation}
Substituting Eq.(4.24) into Eq.(2.1) and equating the coefficients of each $\varepsilon ^k$, we obtain a system
of ODEs :
\begin{equation}
\left\{ \begin{array}{l}
\dot{x}_0 = 0,   \\
\dot{x}_1 = G_1(t, x_0), \\
\quad \vdots \\
\dot{x}_k = G_k(t, x_0, \cdots , x_{k-1}), \\
\quad \vdots
\end{array} \right.
\end{equation}
where the functions $G_k$ are defined through Eq.(2.6).
Solving these equations and substituting solutions $x_k = x_k(t)$ into Eq.(4.24),
we obtain a formal solution $x = \hat{x}(t)$.
If Eq.(2.1) is analytic in $\varepsilon \in \mathcal{E} \subset \mathbf{C}$, Eq.(4.24) converges for $\varepsilon \in \mathcal{E}$
and gives an exact solution of Eq.(2.1).
However, in this section, we need not such an assumption and regard Eq.(4.24) as a formal power series in $\varepsilon $.
To construct a formal solution $\hat{x}(t)$ as above is called the \textit{regular perturbation method}.

Let us define $R_i(y)$ and $u^{(i)}_t(y),\, i = 1,2,\cdots $ by Eqs.(2.11) to (2.14).
By using these functions, regular perturbation solutions (4.24) are given as follows:
\\[0.2cm]
\textbf{Proposition 4.2.}\, Solutions of the system of ODEs (4.25) are given as
\begin{equation}
x_k = x_k(t, y) = u^{(k)}_t(y) + p^{(k)}_1(t,y)t +p^{(k)}_2(t,y)t^2 + \cdots  + p^{(k)}_k(t,y)t^k,
\end{equation}
for $k = 1,2, \cdots $, where $y \in \mathbf{R}^n$ is an arbitrary constant (a solution of $\dot{x}_0 = 0$)
and where $p^{(j)}_i$ are defined by
\begin{equation}
\left\{ \begin{array}{ll}
p^{(1)}_1 (t, y) = p^{(1)}_1(y) = R_1(y) &  \\[0.2cm]
\displaystyle p^{(i)}_1 (t, y) 
   = R_i(y) + \sum^{i-1}_{k=1}\frac{\partial u^{(k)}_t}{\partial y}(y)R_{i-k}(y),  & i= 2,3, \cdots , \\[0.2cm]
\displaystyle p^{(i)}_j (t, y) 
   = \frac{1}{j} \sum^{i-1}_{k=1} \frac{\partial p^{(k)}_{j-1}}{\partial y}(t,y)R_{i-k}(y), &  j=2,3,\cdots ,i-1, \\[0.2cm]
\displaystyle p^{(i)}_i (t, y) = p^{(i)}_i (y) 
   = \frac{1}{i}\frac{\partial p^{(i-1)}_{i-1}}{\partial y}(y)R_{1}(y), & i= 2,3, \cdots, \\
p^{(i)}_j (t, y) = 0, & j > i.
\end{array} \right.
\end{equation}
In particular, $p^{(i)}_j(t,y)$ are almost periodic functions with respect to $t$ for all $i, j \in \mathbf{N}$.
\\[0.2cm]
\textbf{Proof.}\, This proposition is proved in Chiba [10]. \hfill $\blacksquare$
\\[-0.2cm]

Now we derive the RG equation (2.18) according to an idea of Kunihiro [28,29].
We can show that regular perturbation solutions are also expressed as
\begin{equation}
x = \hat{x}(t, \tau, y) = y + \varepsilon x_1(t, \tau, y)+ \varepsilon ^2 x_2(t, \tau, y)+ \cdots 
\end{equation}
with 
\begin{equation}
x_k(t, \tau, y) = u^{(k)}_t(y) + p^{(k)}_1(t,y)(t-\tau ) +p^{(k)}_2(t,y)(t-\tau )^2 + \cdots  + p^{(k)}_k(t,y)(t-\tau )^k,
\end{equation}
by choosing initial times and initial values appropriately when solving the system (4.25) so that
a formal solution $\hat{x}(t, \tau, y)$ passes through $y + \sum^{\infty}_{k=1} \varepsilon ^k u^{(k)}_\tau (y)$ at $t=\tau$,
where $\tau \in \mathbf{R}$ is an arbitrary constant.
Further, we regard $y = y(\tau)$ as a function of $\tau$ to be determined, and consider the family 
$\{\hat{x}(t, \tau, y(\tau)) \}_{\tau \in \mathbf{R}}$ of regular perturbation solutions parameterized by $\tau \in \mathbf{R}$.
It is known that regular perturbation solutions are close to exact solutions only for a short time interval
$|t| \sim O(1)$. However, if we move the initial value $y + \sum^{\infty}_{k=1} \varepsilon ^k u^{(k)}_\tau (y)$
along an exact solution by varying $\tau$ as Fig~.7, 
it seems that the envelope of the family $\{\hat{x}(t, \tau, y(\tau)) \}_{\tau \in \mathbf{R}}$
gives the exact solution of Eq.(2.1).

\begin{figure}[h]
\begin{center}
\includegraphics[scale=1.0]{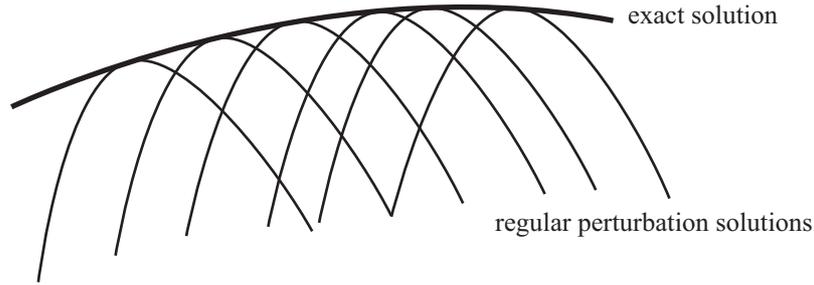}
\end{center}
\caption{Exact solution of Eq.(2.1) and a family of regular perturbation solutions $\hat{x}(t, \tau, y(\tau))$.}
\end{figure}

The envelope is calculated as follows: 
At first, we differentiate $\hat{x}(t, \tau, y(\tau))$ with respect to $\tau = t$ and equate it to zero :
\begin{equation}
\frac{d}{d\tau}\Bigl|_{\tau = t} \hat{x}(t, \tau, y(\tau)) = 0.
\end{equation}
Substituting Eqs.(4.28) and (4.29) into the above yields
\begin{eqnarray}
0 &=& \frac{dy}{dt} + \sum^\infty_{k=1}\varepsilon ^k \frac{d}{d\tau}\Bigl|_{\tau = t} x_k(t, \tau, y(\tau)) \nonumber \\
&=& \frac{dy}{dt} + \sum^\infty_{k=1}\varepsilon ^k 
    \left( \frac{\partial u^{(k)}_t}{\partial y}(y) \frac{dy}{dt} - p^{(k)}_1(t, y) \right)\nonumber \\
&=& \frac{dy}{dt} + \sum^\infty_{k=1}\varepsilon ^k 
   \left( \frac{\partial u^{(k)}_t}{\partial y}(y) \frac{dy}{dt}- R_k(y)
   - \sum^{k-1}_{j=1}\frac{\partial u^{(j)}_t}{\partial y}(y) R_{k-j}(y) \right) .
\end{eqnarray}
This equality is formally satisfied if we define $y = y(t)$ as a solution of the ODE
\begin{equation}
\frac{dy}{dt} = \sum^\infty_{k=1}\varepsilon ^k R_k(y).
\end{equation}
This recovers the $m$-th order RG equation if truncated at order $\varepsilon ^m$.
With this $y(t)$, the envelope for the family $\{\hat{x}(t, \tau, y(\tau)) \}_{\tau \in \mathbf{R}}$
is given by $\hat{x}(t,t, y(t))$.
Again Eqs.(4.28) and (4.29) prove that
\begin{eqnarray}
\hat{x}(t,t, y(t)) &=& y(t) + \sum^\infty_{k=1} \varepsilon ^k x_k(t,t, y(t)) \nonumber \\
&=& y(t) + \sum^\infty_{k=1} \varepsilon ^k u^{(k)}_t(y(t)) \nonumber \\
&=& \alpha ^{(m)}_t(y(t)) + O(\varepsilon ^{m+1}),
\end{eqnarray}
which gives the $m$-th order RG transformation if truncated at order $\varepsilon ^m$.
\\

In the rest of this section, we prove that a formal solution (4.24) with Eq.(4.26) 
obtained by the regular perturbation method
includes infinitely many convergent subseries even if Eq.(2.1) (and thus Eq.(4.24))
is not analytic in $\varepsilon $.
\\[0.2cm]
\textbf{Proposition 4.3.} \, For sufficiently small $|\varepsilon t|$, the series
$\sum^\infty_{l=1}\varepsilon ^l t^l p^{(l+k)}_l (t, y)$ are convergent for $k=0,1,2,\cdots $.
\\[-0.2cm]

To prove Prop.4.3, we need the next simple lemma.
\\[0.2cm]
\textbf{Lemma 4.4.}\, Let $h_1(t)$ and $h_2(t)$ be almost periodic functions with 
$\mathrm{Mod} (h_1) \supset \mathrm{Mod} (h_2) $ (the module of almost periodic functions is defined in Sec.2.1).
If 
\begin{equation}
\lim_{t\to \infty}(h_1(t) - h_2(t)) = 0,
\end{equation}
then $h_1(t) = h_2(t)$ for all $t\in \mathbf{R}$.
\\[0.2cm]
\textbf{Proof of Lemma 4.4.}\, Functions $h_i ,\, i=1,2$ are expanded in Fourier series as $h_i(t)
 = \sum_{\lambda _k \in \mathrm{Mod} (h_1)} c_i(\lambda _k)e^{i\lambda _k t}$.
Then, Eq.(4.34) is rewritten as
\begin{equation}
\lim_{t\to \infty} \sum_{\lambda _k \in \mathrm{Mod} (h_1)} \left( c_1(\lambda _k) - c_2(\lambda _k)\right) e^{i\lambda _k t}=0.
\end{equation}
This proves that $c_1(\lambda _k) = c_2(\lambda _k)$ for all $\lambda _k\in \mathrm{Mod} (h_1)$. \hfill $\blacksquare$
\\[0.2cm]
\textbf{Proof of Prop.4.3.} \, Let 
\begin{equation}
x(t) = y_0 + \sum^\infty_{l=1} \varepsilon ^l t^l p^{(l)}_l (y_0)
  + \sum^\infty_{k=1} \varepsilon ^k \left( u^{(k)}_t (y_0) + \sum^\infty_{l=1} \varepsilon ^l t^l p^{(l+k)}_l (t,y_0)\right)
\end{equation}
be a formal solution of Eq.(2.1) constructed by the regular perturbation method and let $y = y(\varepsilon t, \varepsilon )$
be a solution of the $m$-th order RG equation
\begin{equation}
\frac{dy}{d(\varepsilon t)} = R_1(y) + \varepsilon R_2(y) + \cdots  + \varepsilon ^{m-1}R_m(y)
\end{equation}
for Eq.(2.1). By Thm.2.7, there exist positive constants $C$ and $T$ such that the inequality
\begin{equation}
|| x(t) - \alpha ^{(m)}_t (y(\varepsilon t, \varepsilon )) || < C |\varepsilon |^m , \quad 0 \leq t\leq T/\varepsilon 
\end{equation}
holds if we choose initial values of $x(t)$ and $y(\varepsilon t, \varepsilon )$ appropriately.
Putting $t = T/\varepsilon $ yields
\begin{eqnarray}
& & \alpha ^{(m)}_{T/\varepsilon }(y (T, \varepsilon ))
 = y(T, \varepsilon ) + \sum^{m}_{k=1}\varepsilon ^k u^{(k)}_{T/\varepsilon } (y(T, \varepsilon )) \nonumber \\
 &=& y_0 + \sum^\infty_{l=1} T^l p^{(l)}_l (y_0)
  + \sum^{m-1}_{k=1} \varepsilon ^k \left( u^{(k)}_{T/\varepsilon } (y_0)
   + \sum^\infty_{l=1} T^l p^{(l+k)}_l (T/\varepsilon ,y_0)\right) + O(\varepsilon ^m).
\end{eqnarray}
Note that functions $u^{(k)}_t(y)$ and $p^{(l+k)}_l (t, y)$ are almost periodic functions in $t$
such that $\mathrm{Mod} (u^{(k)}_t)$ and $\mathrm{Mod} (p^{(l+k)}_l )$
are included in $\mathrm{Mod} (g)$.
In particular, they are bounded in $t\in \mathbf{R}$.
Thus taking the limit $\varepsilon  \to 0$ in Eq.(4.39) provides
\begin{equation}
y(T,0) = y_0 + \sum^\infty_{l=1} T^l p^{(l)}_l (y_0).
\end{equation}
This proves that the series $\sum^\infty_{l=1} T^l p^{(l)}_l (y_0)$ is convergent.

Next thing to do is to prove Prop.4.3 for $k=1$.
Using Eq.(4.40) and dividing the both sides of Eq.(4.39) by $\varepsilon $, we obtain
\begin{equation}
\sum^\infty_{k=1}\frac{\varepsilon ^{k-1}}{k!}\frac{\partial ^k y}{\partial \varepsilon ^k} (T,0)
 + \sum^{m}_{k=1}\varepsilon ^{k-1} u^{(k)}_{T/\varepsilon } (y(T, \varepsilon ))
 = \sum^{m-1}_{k=1} \varepsilon ^{k-1} \left( u^{(k)}_{T/\varepsilon } (y_0)
   + \sum^\infty_{l=1} T^l p^{(l+k)}_l (T/\varepsilon ,y_0)\right) + O(\varepsilon ^{m-1}).
\end{equation}
Taking the limit $\varepsilon \to 0$ yields
\begin{equation}
\frac{\partial y}{\partial \varepsilon }(T,0) + \lim_{\varepsilon \to 0}
\left( u^{(1)}_{T/\varepsilon } (y(T,0)) - u^{(1)}_{T/\varepsilon } (y_0)
 - \sum^\infty_{l=1} T^l p^{(l+1)}_l (T/\varepsilon , y_0)\right) = 0.
\end{equation}
Now Lemma 4.4 proves that
\begin{equation}
\sum^\infty_{l=1} T^l p^{(l+1)}_l (t , y_0)
 = \frac{\partial y}{\partial \varepsilon }(T,0) + u^{(1)}_{t} (y(T,0)) - u^{(1)}_{t} (y_0),
\end{equation}
and the series $\sum^\infty_{l=1} T^l p^{(l+1)}_l (t , y_0)$ proves to be convergent.

The proof of Prop.4.3 for $k \geq 2$ is done in a similar way and omitted. \hfill $\blacksquare$

Though it seems that Prop.4.3 has no relationship to the RG method,
it is instructive to notice the equality (4.40).
Since $y(T,0)$ in the left hand side is obtained from a solution of the first order RG equation 
$dy/d(\varepsilon t) = R_1(y)$, the series $\sum^\infty_{l=1} T^l p^{(l)}_l (y_0)$
proves to be determined by only the first order RG equation.
This is remarkable because if we want to obtain $\sum^\infty_{l=1} T^l p^{(l)}_l (y_0)$ by using
the regular perturbation method, we have to solve infinitely many ODEs (4,25).
A similar argument shows that a solution of the $m$-th order RG equation involves all terms
of the form $\varepsilon ^{k+l}t^{l}p^{(l+k)}_l (t, y_0)$ for $k= 0,\cdots ,m-1$ and $l = 1,2,\cdots $ in the 
formal solution (4.24).

%%%%%%%%%%%%%%%%%%%%%%%%%%%%%%%%%%%%%%%%%%%%%%%%%%%%%%%%%%%%%%%%%%%%%%%%%%%%%%%%%%%%%%%%%%%%%%%%%%%%%%%%%%
%%%%%%%%%%%%%%%%%%%%%%%%%%%%%%%%%%%%%%%%%%%%%%%%%%%%%%%%%%%%%%%%%%%%%%%%%%%%%%%%%%%%%%%%%%%%%%%%%%%%%%%%%%

\section{Infinite order RG equation}

Infinite order RG equations and transformations are not convergent series in general.
In this section, we give a necessary and sufficient condition for the convergence.
It will be proved in Sec.5.3 that infinite order RG equations for linear systems are convergent
and related to Floquet theory.

%%%%%%%%%%%%%%%%%%%%%%%%%%%%%%%%%%%%%%%%%%%%%%%%%%%%%%%%%%%%%%%%%%%%%%%%%%%%%%%%%%%%%%%%%%%%%%%%%%%%%%%%%%

\subsection{Convergence condition}

Let us consider the system
\begin{equation}
\dot{x} = \varepsilon g(t,x,\varepsilon ), \quad x\in M
\end{equation}
on a real analytic manifold $M$ with the following assumption :
\\[0.2cm]
\textbf{(A')} \, The vector field $g$ is analytic with respect to $t\in \mathbf{R}, x\in M$ and 
$\varepsilon \in I$, where $I \subset \mathbf{R}$ is an open interval containing $0$.
Further $g$ is $T$-periodic in $t$.
\\[-0.2cm]

Because of the periodicity, we can regard (5.1) as a system on $S^1 \times M$ :
\begin{equation}
\frac{d}{dt}\left(
\begin{array}{@{\,}c@{\,}}
s \\
x
\end{array}
\right) = \left(
\begin{array}{@{\,}c@{\,}}
1 \\
\varepsilon g(s, x, \varepsilon )
\end{array}
\right),
\end{equation}
where $S^1$ is a circle with a $C^\omega $ structure.

Theorem 2.5 states that the infinite order RG transformation defined by
\begin{equation}
x = \alpha _t (y) := y + \varepsilon u^{(1)}_t(y) + \varepsilon ^2 u^{(2)}_t(y) + \cdots 
\end{equation}
formally brings the system (5.1) into the infinite order RG equation
\begin{equation}
\dot{y} = \varepsilon R_1(y) + \varepsilon ^2 R_2(y) + \cdots ,
\end{equation}
where ``formally" means that Eqs.(5.3) and (5.4) are not convergent in general.
If Eq.(5.3) is convergent, so is Eq.(5.4).
A necessary and sufficient condition for the convergence of Eq.(5.3) is given as follows:
\\[0.2cm]
\textbf{Theorem 5.1.}\, For the system (5.1) with the assumption (A'),
there exist an open neighborhood $U = U_y$ of $S^1 \times \{y\} \times \{0\}$ in $S^1 \times M \times I$ for each $y \in M$
and an analytic infinite order RG transformation on $U$,
if and only if the system (5.1) is invariant under the $\mathbb{T}^1$ ($1$-torus) action
of the form
\begin{equation}
\mathbb{T}^1 : (t,x) \mapsto (t + k,\, x + \varepsilon \sigma_k (t,x, \varepsilon )), \quad k\in \mathbf{R},
\end{equation}
where $\sigma_k (t,x,\varepsilon )$ is analytic with respect to 
$k,t,x,\varepsilon$ and $T$-periodic in $k$ and $t$.

Recall that RG equations and RG transformations are not unique and not all of them are convergent
even if the condition of Thm.5.1 is satisfied.

The proof of this theorem involves Lie group theory and will be given in Sec.5.2, though
the idea of the proof is shown below.

Since the infinite order RG equation is an autonomous system, it is invariant under the translation of $t$,
$(t, y) \mapsto (t + k, y)$. If an infinite order RG equation and transformation are convergent and well-defined,
the system (5.1) is invariant under the action defined by pulling back the translation by the RG transformation :
\begin{equation}
(t, x) \mapsto (t + k, \alpha _{t+k} \circ \alpha _t^{-1}(x)).
\end{equation}
Since $\alpha _t$ is $T$-periodic (Lemma 2.1 (ii)), this defines the $\mathbb{T}^1$ action on the space $S^1 \times M$.
Conversely, if the system (5.1) is invariant under the action (5.5), then a simple extension of Bochner's
linearization theorem proves that there exists a $C^\omega $ coordinates transformation $x \mapsto y$ such that
the action (5.5) is written as $(t, y) \mapsto (t + k, y)$.
We can show that this transformation is just an RG transformation.
See Section 5.2 for the detail.
\\[-0.2cm]

In the rest of this subsection, we consider an autonomous system on $\mathbf{C}^n$ of the form
\begin{eqnarray}
\dot{x} &=& Fx + \varepsilon g(x, \varepsilon ) \nonumber \\
&=& Fx + \varepsilon g_1(x) + \varepsilon ^2 g_2(x) + \cdots ,\,\, x\in \mathbf{C}^n,
\end{eqnarray}
with the assumptions (C1),(C2) (see Sec.2.1) and (C3) below :
\\[0.2cm]
\textbf{(C3)} \, $g(x, \varepsilon )$ is analytic with respect to $x \in \mathbf{C}^n$ and $\varepsilon \in I \subset \mathbf{R}$. 
\\[-0.2cm]

For this system, RG transformations are defined by Eq.(2.51) with Eqs.(2.46) to (2.49).
The next corollary immediately follows from Thm.5.1.
\\[0.2cm]
\textbf{Corollary 5.2.} \, Suppose that all eigenvalues of $F$ have pairwise rational ratios.
Then, there exist an open neighborhood $U = U_y$ of $S^1 \times \{y\} \times \{0\}$
in $S^1 \times M \times I$ for each $y \in M$
and an analytic infinite order RG transformation on $U$,
if and only if the system (5.7) is invariant under the $\mathbb{T}^1$ action
of the form
\begin{equation}
\mathbb{T}^1 : x \mapsto e^{Fk}x + \varepsilon \sigma_k (t,x, \varepsilon ), \quad k\in \mathbf{R},
\end{equation}
where $\sigma_k (t,x,\varepsilon )$ is analytic with respect to $k,t,x,\varepsilon$ and periodic in $k$ and $t$.
\\[0.2cm]
\textbf{Proof.}\, By changing the coordinates as $x = e^{Ft}X$, Eq.(5.7) is rewritten as 
$\dot{X} = \varepsilon e^{-Ft}g(e^{Ft}X, \varepsilon )$.
Since $e^{Ft}$ is periodic because of the assumption of Corollary 5.2, we can apply Thm.5.1 to this system.
Note that we do not need the assumption (C2) to prove Corollary 5.2. \hfill $\blacksquare$
\\[-0.2cm]

Recall that RG equations for Eq.(5.7) are equivalent to normal forms (Sec.4.3).
If there are irrational ratios among eigenvalues of $F$, Thm.5.1 is no longer applicable.
For such a system, Zung [54] gives a necessary and sufficient condition for the convergence of normal forms
of infinite order, although he supposes that $g_i(x)$ in Eq.(5.7) is a homogeneous vector field of degree $i$
for $i = 2,3, \cdots $.
A necessary and sufficient condition for the convergence of infinite RG equations and transformations
for Eq.(5.7) is given in a similar way to Zung's theorem as follows:

Remember that an RG equation for Eq.(5.7) has the property that $R_i(e^{Fk}y) = e^{Fk}R_i(y),\, k\in \mathbf{R}$
for $i = 1,2, \cdots $, if integral constants in Eqs.(2.47) and (2.49) are appropriately chosen (Thm.2.15).
If a matrix $B$ satisfies the equalities $FB=BF$ and $Q(e^{Bk}y) = e^{Bk}Q(y),\, k\in \mathbf{R}$
for all polynomial vector fields $Q$ such that $Q(e^{Fk}y) = e^{Fk}Q(y),\, k\in \mathbf{R}$,
then $B$ is called \textit{subordinate} to $F$.
It is known that if $F$ is a diagonal matrix, we can take $B$ having the form 
$B = \mathrm{diag}\, (ib_1, \cdots ,ib_n )$, where $i = \sqrt{-1}$ and $b_j \in \mathbf{Z}$ for $j=1,\cdots ,n$
(see Murdock [39], Zung [54]).
Let $p$ be the maximum number of linearly independent such matrices $B_1, \cdots , B_p$ and call it
the \textit{toric degree of $F$}.
Then, the matrix $e^{B_1k_1 + \cdots + B_p k_p},\, (k_1, \cdots  ,k_p ) \in \mathbf{R}^p$
induces the $\mathbb{T}^p$ action on $\mathbf{C}^n$ and the RG equation is invariant under this action.
By a similar way to the proof of Thm.5.1, we can prove the next theorem, whose proof is omitted here.
\\[0.2cm]
\textbf{Theorem 5.3.} \, Let $p$ be the toric degree of $F$.
For the system (5.7), there exist an open neighborhood $U = U_y$ of 
$S^1 \times \{y\} \times \{0\}$ in $S^1 \times M \times I$ for each $y \in M$
and an analytic infinite order RG transformation on $U$,
if and only if (5.7) is invariant under the $\mathbb{T}^p$ action
of the form
\begin{equation}
\mathbb{T}^p : x \mapsto e^{B_1k_1 + \cdots + B_p k_p}x 
+ \varepsilon \sigma_{k_1, \cdots  ,k_p} (t,x, \varepsilon ), \quad (k_1, \cdots  ,k_p ) \in \mathbf{R}^p,
\end{equation}
where $\sigma_{k_1, \cdots  ,k_p} (t,x,\varepsilon )$ is analytic with respect to 
$k_1, \cdots ,k_p,t,x,\varepsilon $ and periodic in $k_1 , \cdots ,k_p$ and $t$.

%%%%%%%%%%%%%%%%%%%%%%%%%%%%%%%%%%%%%%%%%%%%%%%%%%%%%%%%%%%%%%%%%%%%%%%%%%%%%%%%%%%%%%%%%%%%%%%%%%%%%%%%%%

\subsection{Proof of Theorem 5.1}

In this subsection, we give a proof of Thm.5.1.
At first, we provide a few notations and facts from Lie group theory.

Let $K$ be a compact Lie group, $dk$ a bi-invariant measure on $K$,
$V$ a complete locally convex topological vector space, and $\pi$ a representation of $K$ in $V$.
We define the \textit{average} $\mathrm{av} (\pi) : V \to V$ of $\pi$ by
\begin{equation}
\mathrm{av} (\pi) (v) 
= \int_{K} \! \pi (k)(v) dk. 
\end{equation}
Then, the next theorem holds (see Duistermaat and Kolk [15] for the proof).
\\[0.2cm]
\textbf{Theorem 5.4.}\, The average $\mathrm{av} (\pi)$ is a linear projection from $V$ onto the space
$V^{\pi (K)}:= \{ v\in V \, | \, \pi (k)v = v, \,\, \forall k\in K\}$; that is, equalities
\begin{equation}
\pi (k) \circ \mathrm{av} (\pi) (v) = \mathrm{av} (\pi) (v)
\end{equation}
and 
\begin{equation}
\mathrm{av} (\pi)(v) = v, \quad v\in V^{\pi (K)}
\end{equation}
hold.

Let $S^1$ be a circle with a $C^\omega $ structure, $M$ a real analytic manifold,
and $I \subset \mathbf{R}$ an open interval containing $0$, as Sec.5.1.
The next theorem is a simple extension of Bochner's linearization theorem [5,15].
\\[0.2cm]
\textbf{Theorem 5.5.} \, Let $A : K \to \mathrm{Diff}^\omega (S^1 \times M \times I)$ be a $C^\omega $ action
of $K$ on $S^1 \times M \times I$ such that $A(k)$ is expressed as
\begin{equation}
A(k)(t,x, \varepsilon ) = (a_1(k)(t),\, a_2(k)(t,x,\varepsilon ),\, \varepsilon ),\quad k\in K,
\end{equation} 
where $a_1(k) : S^1 \to S^1$ and $a_2(k) : S^1 \times M \times I \to M$ are $C^\omega $ maps.
Suppose that there exists $x_0 \in M$ such that
\begin{equation}
a_2(k)(t, x_0, 0) = x_0.
\end{equation}
Then, there exist an open neighborhood $U$ of $S^1 \times \{x_0\} \times \{0\}$ in $S^1 \times M \times I$
and a $C^\omega $ diffeomorphism $\overline{\varphi }$ from $U$ into $S^1 \times T_{x_0}M \times \mathbf{R}$
such that
\begin{equation}
\overline{\varphi } \circ A(k) = (a_1(k) \times D_xa_2(k)(t, x_0, 0) \times id) \circ \overline{\varphi },
\end{equation}
where $D_x$ is the derivative with respect to $x$.
The $\overline{\varphi }$ is expressed as
\begin{equation}
\overline{\varphi }(t,x, \varepsilon ) = (t,\, a_3(t,x,\varepsilon ),\, \varepsilon )
\end{equation}
with a $C^\omega $ map $a_3$ and satisfies
\begin{equation}
\overline{\varphi }(t, x_0, 0) = (t,0,0), \quad D_x\overline{\varphi }(t, x_0,0) = id.
\end{equation}
Further if we suppose
\begin{equation}
D_x^n a_2(k)(t, x_0,0) = 0 \quad (n\geq 2),
\end{equation}
then $\overline{\varphi }$ satisfies
\begin{equation}
D_x^n \overline{\varphi }(t, x_0,0) = 0 \quad (n\geq 2).
\end{equation}
\\
\textbf{Proof.}\, Let $W$ be a $K$-invariant open neighborhood of $S^1 \times \{ x_0\} \times \{0\}$
and $V$ a space of $C^\omega $ maps $v : W \to S^1 \times T_{x_0}M \times \mathbf{R}$ such that
\begin{equation}
v(S^1 \times \{ x_0\} \times \{0\}) = S^1 \times \{ x_0\} \times \{0\}.
\end{equation}
Then $V$ is a complete locally convex topological vector space.

Take a $\varphi \in V$ expressed as
\begin{equation}
\varphi (t,x,\varepsilon ) = (t,\, \beta (t,x,\varepsilon ),\, \varepsilon )
\end{equation}
satisfying
\begin{equation}
D_x\beta (t, x_0,0) = id, \quad D^n_x\beta (t, x_0,0) = 0 \quad (n\geq 2),
\end{equation}
where $\beta $ is a $C^\omega $ map.
Define the representation $\pi$ of $K$ in $V$ to be 
\begin{equation}
\pi (k)(v) = (a_1(k) \times D_xa_2(k)(t, x_0, 0) \times id) \circ v \circ A(k)^{-1},
\end{equation}
and define $\overline{\varphi } \in V$ to be
\begin{equation}
\overline{\varphi } = \mathrm{av} (\pi) (\varphi )
 = \int_{K} \!  (a_1(k) \times D_xa_2(k)(t, x_0, 0) \times id) \circ \varphi \circ A(k)^{-1} dk. 
\end{equation}
Then Thm.5.4 implies that $\overline{\varphi } = \pi (k) (\overline{\varphi })$
and this proves Eq.(5.15).
Eqs.(5.16),(5.17) and (5.19) immediately follow from the definition (5.24) of $\overline{\varphi }$
with Eqs.(5.14),(5.18),(5.20),(5.21) and (5.22).
Since $D\overline{\varphi }(t, x_0,0) = id$, $\overline{\varphi }$ is a $C^\omega $ diffeomorphism
on a neighborhood of $S^1 \times \{x_0\} \times \{0\}$ by virtue of the inverse mapping theorem. \hfill $\blacksquare$
\\[0.2cm]
\textbf{Proof of Thm.5.1.}\, Suppose that an infinite order RG transformation $\alpha _t(y)$ is analytic
in $U$. Then, so is an infinite order RG equation. Since the RG equation is invariant under the translation
of $t$, Eq.(5.1) is invariant under the action defined by pulling back the translation by the RG transformation:
\begin{eqnarray}
(t, x) &\mapsto &(t + k,\, \alpha _{t+k} \circ \alpha _t^{-1}(x)) \nonumber \\
& =& (t + k,\, x + \varepsilon (u^{(1)}_{t+k}(x) - u^{(1)}_t(x)) + O(\varepsilon ^2)), \quad k\in \mathbf{R}.
\end{eqnarray}
Since $\alpha _t$ is $T$-periodic, this defines the $\mathbb{T}^1$ action on $S^1 \times M$
of the form of (5.5).

Conversely, suppose that Eq.(5.1) is invariant under the $\mathbb{T}^1$ action (5.5).
Let us rewrite Eq.(5.1) as
\begin{equation}
\frac{d}{dt}\left(
\begin{array}{@{\,}c@{\,}}
s \\
x \\
\varepsilon 
\end{array}
\right) = \left(
\begin{array}{@{\,}c@{\,}}
1 \\
\varepsilon g(s, x, \varepsilon ) \\
0
\end{array}
\right) .
\end{equation}
Then, the action (5.5) induces the action on $S^1 \times M \times I$ of the form
\begin{equation}
A(k) : (s, x, \varepsilon ) \mapsto (s + k,\, a_2(k)(s,x,\varepsilon ), \varepsilon ), \quad k\in \mathbf{R},
\end{equation}
where $a_2(k)(s,x,\varepsilon ) = x + \varepsilon \sigma_k (s,x,\varepsilon )$.
Since $a_2(k)$ satisfies Eqs.(5.14) and (5.18) for any $x\in M$,
Thm.5.5 applies to show that there exist an open neighborhood $U = U_{x_0}$ of $S^1 \times \{x_0\} \times \{0\}$
and a $C^\omega $ diffeomorphism $\overline{\varphi }$ satisfying Eqs.(5.15) to (5.17) and (5.19) for each $x_0\in M$.

By taking a local coordinate near $x_0$, we put $x_0 = 0$ and 
identify a neighborhood of $x_0$ with a neighborhood of $0$ in $T_{x_0}M$.
Then, Eqs.(5.16), (5.17) and (5.19) prove that $\overline{\varphi }$ is expressed as
\begin{equation}
\overline{\varphi }(t, x, \varepsilon ) = (t , x + \varepsilon \hat{a}_3 (t,x,\varepsilon ),\, \varepsilon ),
\end{equation}
where $\hat{a}_3$ is a $C^\omega $ map.
Thus $\overline{\varphi }$ defines the $C^\omega $ transformation $\psi_t$ by
\begin{equation}
x = \psi_t (y) = y + \varepsilon \hat{a}_3 (t,y,\varepsilon ).
\end{equation}
Now Eq.(5.15) proves that if we transform Eq.(5.1) by $x = \psi_t (y, \varepsilon )$, then
the resultant equation is invariant under the action
\begin{equation}
a_1(k) \times D_xa_2(k)(t, x_0, 0) \times id :  (t, x, \varepsilon ) \mapsto (t + k, x, \varepsilon ).
\end{equation}
Since this is the translation of $t$, the resultant equation has to be an autonomous system of the form
\begin{equation}
\frac{dy}{dt} = \varepsilon \mathcal{R}(y, \varepsilon ),
\end{equation}
where $\mathcal{R}$ is analytic in $y$ and $\varepsilon $.

Finally, we show that there exist an RG transformation and an RG equation which coincide with
Eqs.(5.29) and (5.31), respectively.
Let Eqs.(5.3) and (5.4) be one of the RG transformations and the RG equations for Eq.(5.1), respectively.
Since both of Eqs.(5.31) and (5.4) are obtained from Eq.(5.1) by $C^\infty$ transformations
$\psi_t$ and $\alpha _t$, respectively,
there exists an $C^\infty$ transformation $\phi$ which brings Eq.(5.31) into Eq.(5.4):
\begin{equation}
\alpha _t \circ \phi (y, \varepsilon ) = \psi_t (y, \varepsilon ).
\end{equation}
Because of Claim 2.16, there exists an RG transformation $\widetilde{\alpha }_t$ other than $\alpha _t$
such that $\widetilde{\alpha }_t = \alpha _t \circ \phi$.
This proves that there exists an analytic RG transformation $\widetilde{\alpha }_t = \psi_t$.\hfill $\blacksquare$

\subsection{Infinite order RG equation for linear systems}

In this subsection, we consider an $n$-dimensional linear system of the form
\begin{eqnarray}
\dot{x} &=& \varepsilon A(t, \varepsilon )x,, \nonumber \\
&=& \varepsilon A_1(t)x + \varepsilon ^2 A_2(t)x + \cdots ,\quad x\in \mathbf{C}^n,
\end{eqnarray}
where $A(t, \varepsilon )$ is an $n\times n$ matrix.
If this system satisfies the assumption (A), the RG method is applicable.
Since the RG equation and the RG transformation are linear in $y \in \mathbf{C}^n$,
it is convenient to define matrices $R_i$ and $u^{(i)}_t$ to be
\begin{eqnarray}
& & R_1 = \lim_{t\to \infty} \frac{1}{t}\int^t \! A_1(s)ds, \\
& & u^{(1)}_t = \int^t \! \left( A_1(s) - R_1 \right) ds, 
\end{eqnarray}
and
\begin{eqnarray}
& & R_i = \lim_{t\to \infty} \frac{1}{t}\!\int^t \! 
\left(\sum^{i-1}_{k=1} A_{i-k}(s)u^{(k)}_s + A_i(s) - \sum^{i-1}_{k=1} u^{(k)}_s R_{i-k} \right) ds , \\
& & u^{(i)}_t = \int^t \! 
\left( \sum^{i-1}_{k=1} A_{i-k}(s)u^{(k)}_s + A_i(s) - \sum^{i-1}_{k=1} u^{(k)}_s R_{i-k} - R_i \right) ds ,
\end{eqnarray}
for $i = 2,3, \cdots$, respectively.
With these matrices, the $m$-th order RG equation and the $m$-th order RG transformation for Eq.(5.33) are defined by
\begin{eqnarray}
& & \dot{y} = \varepsilon R_1y + \varepsilon ^2R_2y + \cdots + \varepsilon ^m R_my, \quad y\in \mathbf{C}^n, \\
& & \alpha ^{(m)}_ty = y + \varepsilon u^{(1)}_ty + \cdots + \varepsilon ^mu^{(m)}_ty,\quad y\in \mathbf{C}^n,
\end{eqnarray}
respectively.
Since the RG equation is a linear system with a constant coefficient, we can easily determine the stability
of the trivial solution $x = 0$ of Eq.(5.33) by using the RG equation. See Chiba [11] for the detail.

In what follows, we suppose the following assumption :
\\[0.2cm]
\textbf{(L)}\, The matrix $A(t, \varepsilon )$ is $T$-periodic in $t$ and analytic with respect to $\varepsilon \in D$,
where $D \subset \mathbf{C}$ is an open neighborhood of the origin.
\\[-0.2cm]

Let us consider the infinite order RG equation $\dot{y} = \sum^\infty_{k=1}\varepsilon ^kR_ky := \mathcal{R}(\varepsilon )y$
and the infinite order RG transformation $\alpha _t y = y + \sum^\infty_{k=1}\varepsilon ^k u^{(k)}_ty$.
Then, the fundamental matrix $X(t, \varepsilon )$ of Eq.(5.33) is given by
\begin{equation}
X(t, \varepsilon ) = \alpha _t \cdot e^{\mathcal{R}(\varepsilon )t}.
\end{equation}
An initial value is given by $X(0, \varepsilon ) = \alpha _0 = id + \sum^\infty_{k=1}\varepsilon ^k u^{(k)}_0$
and it can be taken arbitrarily by choosing integral constants in Eqs.(5.35) and (5.37) appropriately.
We suppose that the integral constants are chosen so that an initial value $\alpha _0$ is analytic in $\varepsilon \in D$.
Then, $X(t, \varepsilon )$ is analytic in $\varepsilon \in D$.
Now the question arises whether $\alpha _t$ and $\mathcal{R}(\varepsilon )$ are analytic in $\varepsilon $.
\\[0.2cm]
\textbf{Theorem 5.6.}\, Suppose that Eq.(5.33) satisfies the assumption (L) and $\alpha _0$ is chosen to be
analytic in $\varepsilon \in D$.
Then, there exists a positive number $r_0$ such that $\mathcal{R}(\varepsilon ),\, e^{\mathcal{R}(\varepsilon )t}$
and $\alpha _t$ are analytic on the disk $|\varepsilon | < r_0$.
\\[0.2cm]
\textbf{Proof.}\, Since $\alpha _t$ is $T$-periodic (Lemma 2.1 (ii)), the equality
\begin{eqnarray}
X(t + T, \varepsilon ) &=& \alpha _{t+ T} \cdot e^{\mathcal{R}(\varepsilon )(t + T)} \nonumber \\
&=& \alpha _t \cdot e^{\mathcal{R}(\varepsilon )t} \cdot e^{\mathcal{R}(\varepsilon )T} \nonumber \\
&=& X(t, \varepsilon ) \cdot e^{\mathcal{R}(\varepsilon )T}
\end{eqnarray}
holds. Putting $t=0$ yields
\begin{equation}
e^{\mathcal{R}(\varepsilon )T} = X(0, \varepsilon )^{-1} \cdot X(T, \varepsilon ).
\end{equation}
Since $X(t, \varepsilon )$ is nonsingular and analytic in $\varepsilon \in D$ for all $t\in \mathbf{R}$,
$e^{\mathcal{R}(\varepsilon )T}$ is also analytic in $\varepsilon \in D$.
Since $\det e^{\mathcal{R}(\varepsilon )T} \neq 0$, the theory of analytic matrix-functions concludes that
$\mathcal{R}(\varepsilon )T$ is analytic on some disk $\{ \varepsilon \in \mathbf{C} \, | \, |\varepsilon | < r_0\}
\subset D$ (see Yakubovich and Starzhinskii [52], Erugin [17]).
This proves that $\mathcal{R}(\varepsilon ), e^{\mathcal{R}(\varepsilon )t}$ and 
$\alpha _t = X(t, \varepsilon ) e^{-\mathcal{R}(\varepsilon )t}$ are also analytic on the disk. \hfill $\blacksquare$
\\[-0.2cm]

A few remarks are in order.
Thm.5.6 with Thm.5.1 shows that Eq.(5.33) with the assumption (L) admits a Lie group action other than
the scalar multiple : $x \mapsto kx ,\, k\in \mathbf{C}$.

In general, analyticity of $e^{\mathcal{R}(\varepsilon )T}$ on $D$ does not conclude analyticity of 
$e^{\mathcal{R}(\varepsilon )t}$ on $D$ for all $t\in \mathbf{R}$.
For example, consider the case $\mathcal{R}(\varepsilon ) =\log (1 + \varepsilon ),\, T=1$.

The convergence radius $r_0$ of $\mathcal{R}(\varepsilon )$ is given as follows:
Fix $\varepsilon _0 \in D$ and a pass $l$ in $D$ from the origin to $\varepsilon _0$.
Let $\lambda _1(\varepsilon ) , \cdots  ,\lambda _n(\varepsilon )$ be eigenvalues of $e^{\mathcal{R}(\varepsilon )T}$.
Suppose that there are $i \neq j$ such that $\lambda _i(\varepsilon _0) = \lambda _j (\varepsilon _0)$
and passes of $\lambda _i(\varepsilon )$ and $\lambda _j(\varepsilon )$ along $l$ surround the origin (see Fig.8).
Then, $\log \lambda _i (\varepsilon _0)$ and $\log \lambda _j (\varepsilon _0)$ are located in different sheets
of the Riemann surface.
The smallest absolute value $r_0 = |\varepsilon _0|$ among such $\varepsilon _0$'s gives the convergence radius.

\begin{figure}[h]
\begin{center}
\includegraphics[scale=1.0]{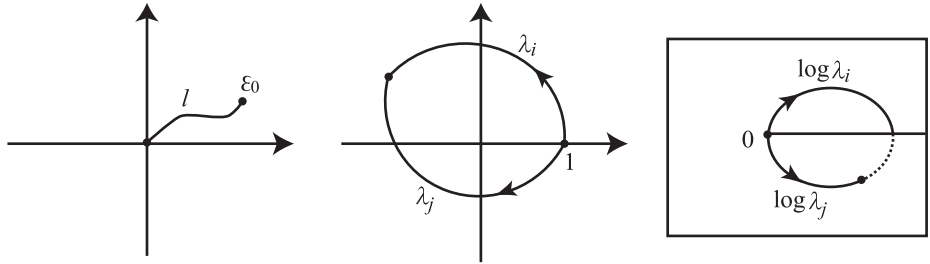}
\end{center}
\caption{Passes of eigenvalues $\lambda _i,\,  \lambda _j$ on $\mathbf{C}$
and  passes of $\log \lambda _i,\, \log \lambda _j$ on the Riemann surface.}
\end{figure}

Floquet theorem states that for a given linear system with a periodic coefficient,
there exist a periodic matrix $Q(t)$ and a constant matrix $B$ such that the fundamental matrix $X(t)$
is written as $X(t) = Q(t)e^{Bt}$.
The RG method just gives the matrix $B$; $B = \mathcal{R}(\varepsilon ) = \sum^\infty_{k=1}\varepsilon ^kR_k$.
Since $Q(t) = \alpha _t$ is periodic, the stability of the trivial solution $x=0$ is determined by eigenvalues of 
$B = \mathcal{R}(\varepsilon )$, called \textit{Floquet exponents}.

Let $T= 2\pi$ for simplicity.
The matrix $A(t, \varepsilon )$ is expanded in a Fourier series as 
$A(t, \varepsilon ) = \sum^{+\infty}_{-\infty} c_n(\varepsilon ) e^{int}$.
By changing the independent variables as $z = e^{int}$, Eq.(5.33) is written as
\begin{equation}
\frac{dx}{dz} = -i \varepsilon \sum^{+\infty}_{-\infty} c_n(\varepsilon ) z^{n-1} x.
\end{equation}
This is a linear system on a complex domain having the singularity at the origin.
It is known that the fundamental matrix of this system is expressed as
\begin{equation}
X(z) = S(z) \cdot e^{M \cdot \log z},
\end{equation}
where $S(z)$ is a single-valued matrix and $M = M(\varepsilon )$ is a constant matrix called
\textit{monodromy matrix}.
It is easy to verify that $\mathcal{R}(\varepsilon ) = iM$ and the RG method provides the monodromy matrix.
\\

\textbf{Acknowledgments}

The author would like to thank Professor Toshihiro Iwai for critical reading of the
manuscript and for useful comments.
\\

\setlength{\baselineskip}{12pt}
{\large\textbf{References}}
\\[-0.3cm]
\begin{enumerate}
\renewcommand{\labelenumi}{[\theenumi]}

\item
L. Arnold et al,
Dynamical systems, 
Lecture notes in mathematics 1609, Springer Verlag, 1995

\item
G. Belitskii,
$C^ \infty$-normal forms of local vector fields,
Acta Appl. Math. 70 (2002), no. 1-3, 23--41

\item
C. M. Bender, S. A. Orszag,
Advanced mathematical methods for scientists and engineers,
McGraw-Hill Book Co., New York, 1978

\item
I. I. Blekhman,
Synchronization in Science and Technology,
ASME press, 1988

\item
S. Bochner,
Compact groups of differentiable transformations,
Ann. of Math. (2) 46, (1945). 372--381

\item
N. N. Bogoliubov, Y. A. Mitropolsky,
Asymptotic methods in the theory of non-linear oscillations,
Hindustan Publishing Corp., Delhi, Gordon and Breach Science Publishers, New York 1961

\item
J. Carr,
Applications of Centre Manifold Theory,
Springer-Verlag, 1981

\item 
L. Y. Chen, N. Goldenfeld, Y. Oono, 
Renormalization group theory for global asymptotic analysis, 
Phys. Rev. Lett. 73 (1994), no. 10, 1311-15

\item 
L. Y. Chen, N .Goldenfeld, Y. Oono, 
Renormalization group and singular perturbations: Multiple scales, boundary layers, and reductive perturbation theory,
Phys. Rev. E 54, (1996), 376-394 

\item
H. Chiba,
$C^1$ approximation of vector fields based on the renormalization group method,
SIAM J. Appl. Dym. Syst., Vol.7, No.3, (2008) 895--932 

\item
H. Chiba,
Simplified Renormalization Group Equations for Ordinary Differential Equations,
J. Differ. Equ. 246, (2009) 1991-2019

\item
H. Chiba,
Approximation of Center Manifolds on the Renormalization Group Method,
J. Math. Phys. Vol.49, 102703 (2008).

\item
H. Chiba, D. Paz\'{o},
Stability of an $[N/2]$-dimensional invariant torus
in the Kuramoto model at small coupling,
Physica D, (in press)

\item
S. N. Chow, C. Li, D. Wang,
Normal forms and bifurcation of planar vector fields, 
Cambridge University Press, 1994

\item
J. J. Duistermaat, J. A. C. Kolk,
Lie groups,
Universitext. Springer-Verlag, Berlin, 2000

\item 
S. Ei, K. Fujii, T. Kunihiro, 
Renormalization-group method for reduction of evolution equations; invariant manifolds and envelopes,
Ann. Physics 280 (2000), no. 2, 236-298 

\item
N. P. Erugin,
Linear systems of ordinary differential equations with periodic and quasi-periodic coefficients,
Academic Press, New York-London, 1966

\item 
N. Fenichel,
Persistence and smoothness of invariant manifolds for flows,
Indiana Univ. Math. J., 21(1971), pp. 193-226

\item
N. Fenichel,
Geometric singular perturbation theory for ordinary differential equations,
J. Differential Equations 31 (1979), no.1, pp.53-98

\item
A. M. Fink,
Almost periodic differential equations,
Lecture Notes in Mathematics, Vol. 377. Springer-Verlag, Berlin-New York, 1974

\item
S. Goto,
Renormalization Reductions for Systems with Delay,
Prog. Theor. Phys. 118. No.2 (2007), 211--227

\item
S. Goto, Y. Masutomi, K. Nozaki,
Lie-Group Approach to Perturbative Renormalization Group Method,
Progr. Theoret. Phys., vol.102, No.3.(1999), pp. 471-497.

\item
W. Hahn,
Stability of motion,
Springer-Verlag New York, 1967

\item 
M.W. Hirsch, C.C. Pugh, M. Shub,
Invariant manifolds,
Springer-Verlag, 1977,
Lec. Notes in Math., 583

\item
F.C. Hoppensteadt, E.M. Izhikevich,
Weakly connected neural networks,
Springer-Verlag, New York, 1997

\item
E. Kirkinis,
On the Reduction of Amplitude Equations by the Renormalization Group Approach,
Phys. Rev. E, (1) 77, 2008

\item
E. Kirkinis,
The Renormalization Group and the Implicit Function Theorem for Amplitude Equations,
J. Math. Phys. (7), 49, 2008

\item 
T. Kunihiro,
A geometrical formulation of the renormalization group method for global analysis,
Progr. Theoret. Phys. 94 (1995), no. 4, 503-514

\item 
T. Kunihiro,
The renormalization-group method applied to asymptotic analysis of vector fields,
Progr. Theoret. Phys. 97 (1997), no. 2, 179-200 

\item
T. Kunihiro, J. Matsukidaira,
Dynamical reduction of discrete systems based on the renormalization-group method,
Phys. Rev. E. Vol.57, No.4, 4817--4820, 1998

\item
T. Kunihiro, K. Tsumura,
Application of the renormalization-group method to the reduction of transport equations,
J. Phys. A 39 (2006), no.25, 8089--8104

\item
R. E. Lee DeVille, A. Harkin, M. Holzer, K. Josi\'{c}, T. Kaper,
Analysis of a renormalization group method and normal form theory for perturbed ordinary differential equations,
Physica D, (2008)

\item
I.G. Malkin,
The Methods of Lyapunov and Poincare in the Theory of Nonlinear Oscillations,
Moscow-Leningrad, 1949

\item
I.G. Malkin,
Some problems of the theory of nonlinear oscillations,
Gosudarstv. Izdat. Tehn.-Teor. Lit., Moscow, 1956

\item
B. Mudavanhu, R.E. O'Malley Jr.,
A new renormalization method for the asymptotic solution of weakly nonlinear vector systems,
SIAM J. Appl. Math. 63 (2002), no.2, 373--397

\item
B. Mudavanhu, R.E. O'Malley Jr.,
A renormalization group method for nonlinear oscillators,
Stud. Appl. Math. 107 (2001), no. 1, 63--79

\item
J. Murdock,
Frequency entrainment for almost periodic evolution equations,
Proc. Amer. Math. Soc. 96 (1986), no. 4, 626--628

\item
J. Murdock,
Perturbations. Theory and methods,
John Wiley \& Sons, Inc., New York, 1991

\item 
J. Murdock,
Normal forms and unfoldings for local dynamical systems,
Springer-Verlag, New York, (2003)

\item 
J. Murdock,
Hypernormal form theory: foundations and algorithms,
J. Differential Equations, 205 (2004), no. 2, 424-465

\item
J. Murdock, L.C. Wang,
Validity of the multiple scale method for very long intervals,
Z. Angew. Math. Phys. 47 (1996), no. 5, 760--789

\item
S. Murata, K. Nozaki,
Renormalization group symmetry method and gas dynamics,
Internat. J. Non-Linear Mech. 39 (2004), no. 6, 963--967

\item
A. H. Nayfeh, 
Method of Normal Forms,
John Wiley \& Sons, INC., 1993

\item
A. H. Nayfeh,
Introduction to perturbation techniques,
John Wiley \& Sons, New York, 1981

\item 
K. Nozaki, Y. Oono,
Renormalization-group theoretical reduction,
Phys. Rev. E, 63(2001), 046101

\item
K. Nozaki, Y. Oono, Y. Shiwa,
Reductive use of renormalization group,
Phys. Rev. E (3) 62 (2000), no. 4, 4501--4504

\item
J.A. Sanders, F. Verhulst, J. Murdock,
Averaging methods in nonlinear dynamical systems,
Springer, New York, 2007

\item
T. Tsumura, T. Kunihiro, K. Ohnishi,
Derivation of covariant dissipative fluid dynamics in the renormalization-group method,
Phys. Lett. B 646 (2007), no. 2-3, 134--140.

\item
S.I. Tzenov, R.C. Davidson,
Renormalization group reduction of the Henon map and application to the transverse betatron motion in cyclic accelerators,
New J. Phys. 5 (2003), 67.1--67.13

\item
F. Verhulst,
Methods and applications of singular perturbations.
Boundary layers and multiple timescale dynamics,
Texts in Applied Mathematics, 50. Springer, New York, (2005)

\item 
S. Wiggins, 
Normally hyperbolic invariant manifolds in dynamical systems ,
Springer-Verlag, 1994

\item
V. A. Yakubovich, V. M. Starzhinskii,
Linear differential equations with periodic coefficients 1,
John Wiley \& Sons, New York-Toronto, 1975

\item
M. Ziane,
On a certain renormalization group method, 
J. Math. Phys. 41 (2000), no. 5, 3290-3299

\item
N. T. Zung,
Convergence versus integrability in Poincare-Dulac normal form,
Math. Res. Lett. 9 (2002), no. 2-3, 217--228

\item
I. Moise, R. Temam,
Renormalization group method. Application to Navier-Stokes equation,
Discrete Contin. Dynam. Systems 6 (2000), no. 1, 191--210

\item
I. Moise, M. Ziane,
Renormalization group method. Applications to partial differential equations,
J. Dynam. Differential Equations 13 (2001), no. 2, 275--321

\item
M. Petcu, R. Temam, D. Wirosoetisno,
Renormalization group method applied to the primitive equations,
J. Differential Equations 208 (2005), no. 1, 215--257

\item
Y. Shiwa,
Renormalization-group theoretical reduction of the Swift-Hohenberg model,
Phys. Rev. E 63, 016119 (2000)

\item
Tao. Tu, G. Cheng,
Renormalization group theory for perturbed evolution equations,
Phys. Rev. E (3) 66 (2002), no. 4, 046625

\end{enumerate}

\end{document}